\documentclass[3p,times]{elsarticle}

\usepackage{graphicx}
\usepackage{psfrag}
\usepackage{color}
\usepackage{amssymb}
\usepackage{amsmath}
\usepackage{dsfont}
\usepackage{rotating}

\newcommand{\K}{\mathbf{K}}
\newcommand{\M}{\mathbf{M}}

\newcommand{\U}{\mathcal{U}}

\renewcommand{\S}{\mathbf{S}} 
\newcommand{\Q}{\mathbf{Q}} 
\newcommand{\q}{\mathbf{q}} 
\newcommand{\f}{\mathbf{f}} 
\newcommand{\g}{\mathbf{g}}

\newcommand{\F}{\mathbf{F}} 
 
\newcommand{\x}{\mathbf{x}} 
\newcommand{\xxi}{\boldsymbol{\xi}} 
\newcommand{\X}{\mathbf{X}} 
\renewcommand{\P}{\mathbf{P}} 
\newcommand{\w}{\mathbf{w}}

\renewcommand{\v}{\mathbf{v}}

\begin{document}

\begin{frontmatter}



\title{High Order Cell-Centered Lagrangian-Type Finite Volume Schemes with Time-Accurate Local Time Stepping on Unstructured Triangular Meshes}

\author[UNITN]{Walter Boscheri}
\ead{walter.boscheri@unitn.it}
\author[UNITN]{Michael Dumbser$^{*}$}
\ead{michael.dumbser@unitn.it}
\cortext[cor1]{Corresponding author}
\author[UNITN]{Olindo Zanotti}
\ead{olindo.zanotti@unitn.it}

\address[UNITN]{Laboratory of Applied Mathematics, Department of Civil, Environmental and Mechanical Engineering \\ University of Trento, Via Mesiano 77, I-38123 Trento, Italy}

\begin{abstract}
% Text of abstract
We present a novel cell-centered direct Arbitrary-Lagrangian-Eulerian (ALE) finite volume scheme on unstructured triangular meshes that is high order accurate in space and 
time and that also allows for \textit{time-accurate local time stepping} (LTS). It extends our previous investigations on high order Lagrangian finite volume schemes with 
LTS carried out in \cite{ALELTS1D} in one space dimension. 
The new scheme uses the following basic ingredients: a high order WENO reconstruction in space on unstructured meshes, an element-local high-order accurate space-time Galerkin 
predictor that performs the time evolution of the reconstructed polynomials within each element, the computation of numerical ALE fluxes at the moving element interfaces through 
approximate Riemann solvers, and a one-step finite volume scheme for the time update which is directly based on the integral form of the conservation equations in space-time. The 
inclusion of the LTS algorithm requires a number of crucial extensions, such as a proper scheduling criterion for the time update of each element and for each node; a virtual 
projection of the elements contained in the reconstruction stencils of the element that has to perform the WENO reconstruction; and the proper computation of the fluxes through 
the space-time boundary surfaces that will inevitably contain hanging nodes in time due to the LTS algorithm. 

We have validated our new unstructured Lagrangian LTS approach over a wide sample of test cases solving the Euler equations of compressible gasdynamics in two space dimensions, 
including shock tube problems, cylindrical explosion problems, as well as specific tests typically adopted in Lagrangian calculations, such as the Kidder and the Saltzman  
problem. When compared to the traditional global time stepping (GTS) method, the newly proposed LTS algorithm allows to reduce the number of element updates in a given simulation by 
a factor that may depend on the complexity of the dynamics, but which can be as large as $\sim 4.7$. 
\end{abstract}
\end{frontmatter}

% main text

%
\section{Introduction}
\label{sec.introduction}

In the last few years there has been a renewed interest in the development of novel accurate and robust cell-centered Lagrangian finite volume schemes for hydrodynamics.  
Since in a Lagrangian method the computational mesh moves with the local fluid velocity, such schemes are regarded as the first choice in all problems presenting moving 
material interfaces appearing in compressible multi-phase and multi-material flows, for instance, the numerical simulation of inertial confinement fusion (ICF). 
%The use of a Lagrangian approach allows, for instance, to identify and locate with very good accuracy the position of 
%material interfaces, 
The vast majority of modern Lagrangian schemes adopts a \textit{cell-centered} finite-volume approach, see for example 
\cite{Despres2009,Depres2012,ShashkovCellCentered,MaireCyl1,Maire2007,Maire2008,Maire2009b,Maire2010}, where all flow variables are defined as cell-averaged quantities 
inside a control volume. However, also \textit{staggered} Lagrangian schemes are possible, see e.g. \cite{StagLag}, where the velocity is defined at the cell interfaces, 
while the other flow variables are still defined at the cell centers. 

In \cite{munz94,Despres2009} Godunov-type finite volume schemes have been presented for Lagrangian hydrodynamics, while in \cite{DepresMazeran2003,Despres2005} the governing equations 
have been coupled with the equations for the evolution of the geometry and the resulting weakly hyperbolic system has been solved using a node-based finite volume solver. Unstructured 
multidimensional meshes have been considered by Maire in \cite{Maire2009,Maire2010,Maire2011}, who developed up to second order accurate cell-centered Lagrangian schemes where the time 
derivatives of the fluxes have been computed with a node-centered solver. This approach may be regarded as a multi-dimensional Lagrangian extension of the Generalized Riemann problem 
methodology used for example in the ADER approach of Titarev and Toro \cite{titarevtoro,Toro:2006a} in the Eulerian context. 
Arbitrary-Lagrangian-Eulerian (ALE) methods based on remeshing and remapping have also been investigated very recently for single and multi-material flows in 
\cite{ShashkovMultiMat1,ShashkovMultiMat2,ShashkovMultiMat3,ShashkovMultiMat4}. 

In \cite{chengshu1,chengshu2} Cheng and Shu presented the first better than second order accurate Lagrangian schemes for hydrodynamics on structured meshes, where the use of a high order 
Essentially Non-Oscillatory (ENO) reconstruction operator yielded high order of accuracy in space, while high order of accuracy in time was guaranteed using either a Runge-Kutta or a 
Lax-Wendroff-type time stepping. Arbitrary high order accurate cell-centered Lagrangian-type finite volume schemes for conservative and non-conservative hyperbolic PDE on moving unstructured 
triangular and tetrahedral meshes have been considered for the first time by Boscheri et al. in a very recent series of papers \cite{Lagrange2D,LagrangeNC,LagrangeMHD,LagrangeMDRS,Lagrange3D}. 
A new class of \textit{meshless} Lagrangian particle methods based on a high order accurate moving least-squares WENO reconstruction has been forwarded in \cite{WENOSPH}. 

High order accurate Lagrangian algorithms using the classical continuous finite element method (FEM) can be found, for example, in the work of Scovazzi et al. \cite{scovazzi1,scovazzi2} 
and Dobrev et al. \cite{Dobrev1,Dobrev2,Dobrev3}, while Lagrangian discontinuous Galerkin finite elements have been recently proposed by Vilar et al. and Yu et al. in 
\cite{Vilar1,Vilar2,Vilar3,Yuetal}. Arbitrary-Lagrangian-Eulerian DG schemes have been developed and applied, for example, in \cite{Feistauer4,Feistauer3}. 

Almost all of the above mentioned algorithms use an explicit \textit{global} time stepping scheme in which the timestep is computed under a classical \textit{global} CFL stability condition, 
so that the timestep is essentially determined by the smallest control volume appearing in the mesh. In Lagrangian hydrodynamics, where the mesh follows as closely as possible the local fluid 
motion, very severe deformations and distortions may occur in the computational cells, especially at shocks and shear waves. As a consequence, the computational efficiency of the 
algorithm drastically decreases, because the smallest timestep imposed by the most deformed control volumes dictates the timestep for the entire computational grid, including those elements 
which are much bigger or which lie in a zone where the fluid is moving uniformly. 
In the Eulerian framework such a problem can be partially avoided controlling the mesh quality \textit{a priori} and designing a high quality mesh once in a pre-processing step, 
since the grid will not change anymore during the simulation. 
Of course, the CFL condition can be circumvented by using implicit or semi-implicit schemes, see for example \cite{Casulli1999,CasulliStelling2011,CasulliVOF,Dolejsi1,Dolejsi2,Dolejsi3}, but 
this approach does not yet seem to be very popular in the context of cell-centered Lagrangian-type finite volume methods. 
An alternative to overcome the global CFL condition consists in the development of numerical schemes that allow for time-accurate \textit{local} time stepping (LTS), where each 
element has to obey only a less restrictive \textit{local} CFL stability condition, hence using its own optimal local timestep. Therefore, many efforts have been devoted to the construction of 
high order accurate Eulerian schemes with time-accurate LTS, developing either discontinuous Galerkin finite element methods \cite{FlahertyLTS,dumbserkaeser06d,TaubeMaxwell,stedg1,stedg2,KrivodonovaLTS,FerrariLevelSet} 
or high order accurate finite volume schemes with LTS \cite{Berger-Oliger1984,Berger-Colella1989,DomainDecomp,CastroLTS,MuletAMR1,MuletAMR2,Burger2012,FerrariLevelSet,Dumbser2012b,Dumbser2014}. 
The finite volume schemes with LTS adopt mainly classical adaptive mesh refinement (AMR) techniques in space and time or block-clustered local time stepping algorithms. 
In \cite{GroteLTS1,GroteLTS2} also high order accurate Runge-Kutta time integrators with local time stepping (so-called multi-rate integrators) can be found. 
To the knowledge of the authors, the first high order accurate \textit{Lagrangian} algorithm with \textit{time accurate local time stepping} on moving grids has been proposed 
very recently in \cite{ALELTS1D}, where the equations of hydrodynamics and of classical magnetohydrodynamics (MHD) have been solved in one spatial dimension. 
In the present paper we extend the algorithm presented in \cite{ALELTS1D} to \textit{moving unstructured triangular meshes}. 

The rest of the paper is structured as follows: in Section \ref{sec.numethod} the numerical scheme is described, including the details of the local time stepping algorithm on 
moving unstructured meshes, while numerical convergence studies as well as some classical numerical test problems for hydrodynamics are presented in Section 
\ref{sec.validation}. We conclude the paper giving an outlook to future research and  developments in Section \ref{sec.concl}. 

\section{Numerical method}
\label{sec.numethod}
\subsection{Formulation of the equations and basic set-up}
In this article we consider nonlinear hyperbolic conservation laws,  cast in the form 
\begin{equation}
\label{PDE}
  \frac{\partial \Q}{\partial t} + \nabla \cdot \F(\Q) = \S(\Q), \qquad \x = (x,y) \in \Omega(t) \subset \mathbb{R}^2, \quad t \in \mathbb{R}_0^+, \quad \Q \in \Omega_{\Q} \subset \mathbb{R}^\nu,     
\end{equation} 
where $\Q=(q_1,q_2,...,q_\nu)$ is the vector of conserved variables defined in the space of the admissible states $\Omega_{\Q} \subset \mathbb{R}^\nu$, $\F(\Q)=\left( \f(\Q),\g(\Q) \right)$ denotes the nonlinear flux tensor and $\S(\Q)$ represents a nonlinear algebraic source term which is not allowed to be stiff. The system of equations \eqref{PDE} is defined in two space dimensions, hence addressing the space coordinate vector and the time with $\mathbf{x}=(x,y)$ and $t$, respectively. The two-dimensional computational domain $\Omega(t)$ is time-dependent since in the Lagrangian framework the mesh is moving, hence changing its configuration at each time step. The domain is discretized using a total number of $N_E$ conforming triangles $T^n_i$, therefore at a general time $t^n$ the \textit{current triangulation} $\mathcal{T}^n_{\Omega}$ of the domain $\Omega(t^n)=\Omega^n$ is given by the union of all elements, i.e.
\begin{equation}
\mathcal{T}^n_{\Omega} = \bigcup \limits_{i=1}^{N_E}{T^n_i}. 
\label{trian}
\end{equation}
Within the Lagrangian LTS algorithm that is going to be presented in this paper each element moves in time independently from the others and following its own local timestep, 
hence the triangulation $\mathcal{T}^n_{\Omega}$ will in general \textit{never} be assembled at a common time level. 
In the LTS framework hanging nodes in time are naturally produced and one is in general not able to define the configuration of the computational mesh at a certain time level $t^n$,  
unless we force the computation to reach the same time $\overline{t}$, which could be typically the case either at the final time or at an intermediate output time. For this reason 
in the rest of the paper each time level $t^n$ will be addressed also with the element number it refers to, i.e. $t^n_i$, with the subscript $i$ denoting the number of the physical 
triangle $T_i$. 

As suggested in \cite{Lagrange2D}, we adopt a spatial reference system $\xi-\eta$ defined by the coordinate vector $\boldsymbol{\xi}=(\xi,\eta)$ where the unit reference triangle $T_e$ is composed of the nodes $\boldsymbol{\xi}^e_{1}=(\xi^e_{1},\eta^e_{1})=(0,0)$, $\boldsymbol{\xi}^e_{2}=(\xi^e_{2},\eta^e_{2})=(1,0)$ and $\boldsymbol{\xi}^e_{3}=(\xi^e_{3},\eta^e_{3})=(0,1)$. The physical element $T^n_i$ defined in the physical system $x-y$ is mapped to the reference element $T_e$ using the transformation
\begin{equation} 
 \x = \x(\xxi,t^n) = \X^n_{1,i} + \left( \X^n_{2,i} - \X^n_{1,i} \right) \xi + \left( \X^n_{3,i} - \X^n_{1,i} \right) \eta, 
 \label{xietaTransf} 
\end{equation} 
where $\mathbf{X}^n_{k,i} = (X^n_{k,i},Y^n_{k,i})$ represents the vector of physical coordinates of the $k$-th vertex of triangle $T^n_i$ at time $t^n_i$. In the Lagrangian framework the use of the reference system, which does not change in time, is much more convenient rather than carrying on the computation in the physical system, where elements are moving and deforming in time.

As usual for cell-centered finite volume schemes, data are represented and evolved in time within each control volume as piecewise constant cell averages 
\begin{equation}
  \Q_i^n = \frac{1}{|T_i^n|} \int_{T^n_i} \Q(\x,t^n_i) dV,     
 \label{eqn.cellaverage}
\end{equation}  
where the volume of element $T_i^n$ is denoted by $|T_i^n|$ at the current element time $t^n_i$. In the time-accurate LTS algorithm a cell $T_i^n$ is allowed to evolve the solution in time only if the so-called \textit{update criterion} \cite{DumbserKaeser07,stedg1,ALELTS1D} is satisfied, namely if
\begin{equation}
\max \limits_{j \in \mathcal{N}_i} \left(t_j^n\right) \leq \left(t_i^n + \Delta t_i^n \right) \leq \min \limits_{j \in \mathcal{N}_i} \left(t_j^n + \Delta t_j^n \right), 
\label{eqn.updateFV}
\end{equation} 
where $\mathcal{N}_i$ denotes the \textit{Neumann neighborhood} of element $T_i$, i.e. the three direct side neighbors $T_j$ of the cell, while $t_i^n$ and $\Delta t_i^n$ represent the current 
local time and the local timestep of triangle $T_i$, respectively. Hence, $\left(t_i^n + \Delta t_i^n \right)$ is the future time of element $T_i$ and to make notation easier it will be addressed 
with $t_i^{n+1}$.

There are two important issues that need to be clarified:
\begin{enumerate}
	\item 
	In order to develop a numerical scheme that evolves the cell averages \eqref{eqn.cellaverage} with high order 
	of accuracy in space and in time in one single step, two strategies are followed. For the accuracy in space we implement a 
	suitable Weighted Essentially Non-Oscillatory (WENO) reconstruction technique that is able to deal with LTS and which is presented in detail in the next 
	Section \ref{sec.weno}, while for the accuracy in time we use an element-local space-time Galerkin predictor approach, as illustrated in Section \ref{sec.lst}. 
	\item In a time-accurate LTS finite volume scheme, each element $T^n_i$ evolves the solution $\Q_i^n$ in time with a local timestep $\Delta t_i^n$ that is computed according to 
	a local CFL stability condition. As a result, the WENO reconstruction will be carried out \textit{locally}, i.e. considering only the element $T^n_i$ which is 
	currently updating the solution to its new time level $t^{n+1}_i$, as well as an appropriate neighborhood of $T_i^n$ that is necessary to carry out the reconstruction,
	the so-called reconstruction stencil $\mathcal{S}_i^W$. Since the neighbor elements of $T_i$ in general have a different local time, the reconstruction needs to get time-accurate
	\textit{virtual} cell averages from the neighbor cells as input. These virtual cell averages are readily available from the local space-time Galerkin predictor solution 
  inside the neighbors. 
\end{enumerate}

\subsection{High order WENO reconstruction for local time stepping} 
\label{sec.weno}
In order to obtain high order of accuracy in space a nonlinear WENO reconstruction algorithm is used. As done in \cite{Lagrange2D,LagrangeNC,LagrangeMHD,LagrangeMDRS,Lagrange3D} 
we adopt the \textit{polynomial} formulation presented in \cite{DumbserKaeser07,DumbserKaeser06b,MixedWENO2D,MixedWENO3D}, instead of the original \textit{pointwise} approach 
proposed by Shu et al. in \cite{shu_efficient_weno,balsarashu,HuShuVortex1999,ZhangShu3D}. Other high order accurate reconstruction algorithms on unstructured meshes can be found, 
e.g. in \cite{AboiyarIske,MOODorg,MOODhighorder,ADERMOOD}.  
While all the details of high order WENO reconstruction are contained in the above-mentioned references, we present here only a brief summary of the main features of the scheme, 
highlighting the modifications that are necessary to handle a time accurate local time stepping formulation. 

The reconstructed solution $\mathbf{w}_h(\x,t_i^n)$ is given in terms of piecewise polynomials of degree $M$ and is computed \textit{locally} for each control volume $T_i^n$. First, one has to construct a set of reconstruction stencils $S_i^s$ relative to the element $T_i$, namely
\begin{equation}
\mathcal{S}_i^s = \bigcup \limits_{j=1}^{n_e} T_{m(j)}, 
\label{stencil}
\end{equation}
where $1\leq j \leq n_e$ denotes a local index which counts the elements belonging to the stencil, while $m(j)$ maps the local counter $j$ to the global element number used in the triangulation 
\eqref{trian}. As explained in \cite{StencilRec1990,Olliver2002,KaeserIske2005,DumbserKaeser06b}, in two space dimensions on unstructured meshes one has to take a total number of elements $n_e$ 
for each stencil that is bigger than the smallest number $\mathcal{M} = (M+1)(M+2)/2$ needed to reach the formal order of accuracy $M+1$, hence we typically set $n_e = 2 \mathcal{M}$. 
Furthermore, according to \cite{KaeserIske2005,DumbserKaeser06b}, we need a total number of stencils $s=7$ in order to perform the polynomial WENO reconstruction, namely one central stencil $s=1$, three 
primary sector stencils $s \in \{2,3,4\}$ and  three reverse sector stencils $s \in \{5,6,7\}$. As a consequence, the update criterion \eqref{eqn.updateFV} must be extended to the {\em total} 
WENO stencil $\mathcal{S}_i^W$ given by 
\begin{equation}
\mathcal{S}_i^W= \bigcup \limits_{s=1}^{7} \mathcal{S}_i^s,
\label{eqn.WENOstencil}
\end{equation}
hence obtaining
\begin{equation}
\max \left(t_j^n\right) \leq t_i^{n+1} \leq \min \left(t_j^{n+1} \right), \qquad \forall T_j \in \mathcal{S}_i^W.
\label{eqn.WENOupdate}
\end{equation}
In order to guarantee that at least one element in the entire mesh satisfies condition \eqref{eqn.WENOupdate}, the total stencils $\mathcal{S}_i^W$ need to be constructed in such a way that 
they are \textit{symmetric}, i.e. each element $T_j \in \mathcal{S}_i^W$ inside the stencil of $T_i$ \textit{must} contain in its own WENO stencil $\mathcal{S}_j^W$ the element $T_i$. In other 
words, if $T_j \in \mathcal{S}_i^W$ then $T_i \in \mathcal{S}_j^W$. It is always possible to construct such symmetric stencils by adding elements to the stencils until the condition of 
symmetry is satisfied for all elements. 

For the sake of clarity we give a simple example of what could happen if we take \textit{non-symmetric} stencils. Let element $T_j$ be \textit{not} contained in the stencil of $T_i$ and let $T_i$ belong to the stencil $\mathcal{S}_j^W$ of element $T_j$. Let furthermore the current time level of $T_i$ and $T_j$ be $t_i^n$ and $t_j^n$, respectively, with the corresponding future times $t_i^{n+1}$ and $t_j^{n+1}$. Without loss of generality we assume $t_i^n=t_j^n$, while the future time levels are chosen such that $t_i^{n+1}>t_j^{n+1}$. If the update criterion on the \textit{non-symmetric} stencil $\mathcal{S}_i^W$ is supposed to be satisfied, then element $T_i$ is allowed to update the numerical solution to its future time, which will subsequently become the \textit{current} time of $T_i$, i.e. $t_i^n \rightarrow t_i^{n+1}$. The resulting situation will lead to a \textit{dead lock} in the algorithm, where element $T_j$ will never obey condition \eqref{eqn.WENOupdate} since $t_j^{n+1}<t_i^n$. A simple solution is to always build a \textit{symmetric} stencil. In this case element $T_j$ performs the update \textit{first} and does not prevent element $T_i$ from updating its solution. The drawback of this approach is that slightly larger stencils are required.

Due to \eqref{eqn.WENOupdate}, the current time $t_j^n$ of the neighbor elements belonging to the WENO stencil $\mathcal{S}_i^W$ must be lower than the current time level $t_i^n$ of the triangle $T_i$ 
for which the reconstruction has to be performed. Moreover, in Lagrangian algorithms the mesh is moving in time, therefore the local WENO reconstruction is carried out on a \textit{virtual geometry}  
with \textit{virtual cell averages}, as suggested in \cite{ALELTS1D}. 
These virtual cell averages, which are needed for the reconstruction, are obtained from the local space-time predictor solution $\q_h(\x,t_i^n)$ inside the neighbor elements $T_j$ using 
a simple integral \textit{projection} (averaging). 
The way how this predictor solution is computed will be described in the next Section \ref{sec.lst}. A similar projection is used also for the virtual geometry of the elements inside the 
total WENO stencil, where all elements $T_j^n \in \mathcal{S}_i^W$ are moved \textit{virtually} until time $t_i^n$ is reached. We emphasize that the projection of the stencil geometry and 
of the cell averages is done only virtually, just for the purpose of reconstruction, because the real mesh motion and the real conservative update of the cell averages will be performed 
individually by each element at its scheduled time according to the update criterion \eqref{eqn.WENOupdate}.  
The geometry of each stencil element $T_j^n$, i.e. the vertex coordinates, are projected and also all the other geometric quantities used for the computation, e.g. normal vectors, volumes, side 
lengths, \textit{etc.}. For the sake of clarity, the projected quantities will be denoted by a tilde symbol in the following, hence 
\begin{equation}
 \mathbf{\tilde{X}}^{n+1}_{k,j}	= \mathbf{X}^{n}_{k,j}	+ \left( t_i^n - t_j^n\right) \, \overline{\mathbf{V}}_{k,j}^n, \qquad \forall T_j^n \in \mathcal{S}_i^W, \quad k=1,2,3
\label{eqn.virtGeom}
\end{equation}      
and
\begin{equation}
  \tilde{\Q}_j^n = \left\{ \begin{array}{ccc} 
  \Q_i^n, & \textnormal{ if } & j = i, \\ 
	\frac{1}{|\tilde{T}_j^n|} \int_{\tilde{T}^n_j} \q_h(\x,t^n_i) dV, & \textnormal{ if } & j \neq i, 
	\end{array} \right.  \qquad \forall T_j \in \mathcal{S}_i^W. 
 \label{eqn.virtCellaverage}
\end{equation} 
In \eqref{eqn.virtGeom} the time-averaged node velocity $\overline{\mathbf{V}}_{k,j}^n$ is computed according to the node solver algorithm, see also \cite{Lagrange2D,Lagrange3D,LagrangeMHD,LagrangeMDRS}, 
which will be briefly described in Section \ref{sec.meshMot}, while in \eqref{eqn.virtCellaverage} the virtual cell averages $\tilde{\Q}_j^n$ of the neighbor elements are 
given as the spatial integral of the predicted solution at  time $t_i^n$ over the virtual control volumes $\tilde{T}^n_j$. 

Once the virtual geometry and the virtual cell averages have been computed for the entire stencil $\mathcal{S}_i^W$, we are in the position to carry out the \textit{local} high order WENO reconstruction procedure. To obtain the reconstruction polynomial $\mathbf{w}_h(\x,t_i^n)$, integral conservation of the projected cell averages $\tilde{\Q}^n_j$ in each reconstruction stencil $\mathcal{S}_i^s$ is required, i.e.
\begin{equation}
\label{intConsRec}
\frac{1}{|\tilde{T}^n_j|} \int \limits_{\tilde{T}^n_j} \mathbf{w}^s_h(\x,t_i^n) dV = 
\frac{1}{|\tilde{T}^n_j|} \int \limits_{\tilde{T}^n_j} \psi_l(\xi,\eta) \hat \w^{n,s}_{l,i} = 
\tilde{\Q}^n_j, \qquad \forall T^n_j \in \mathcal{S}_i^s,     
\end{equation} 
where the integrals are evaluated using Gaussian quadrature formulae of suitable order (see \cite{stroud} for details).
For simplicity, in the above equation, as well as in the rest of the paper, 
we have adopted the Einstein summation convention over repeated indices.
The reconstruction polynomial on each stencil is expressed in terms of a set of orthogonal spatial basis functions $\psi_l(\xi,\eta)$ on the reference element, see \cite{Dubiner,orth-basis,CBS-book}, 
and $\mathcal{M}$ unknown degrees of freedom $\hat \w^{n,s}_{l,i}$.  
Since each stencil contains a total number of elements $n_e > \mathcal{M}$, system \eqref{intConsRec} results in an overdetermined linear algebraic system that is solved by a constrained least-squares  technique \cite{DumbserKaeser06b}. In the Lagrangian framework the geometry evolves in time. Hence, the reconstruction matrix, which is given by the multidimensional integrals in \eqref{intConsRec},  continuously changes in time. As a consequence, the system \eqref{intConsRec} must be solved whenever element $T_i$ performs its WENO reconstruction. To maintain the scheme as simple as possible and  reasonably cost efficient, the stencil topology is fixed once and forall in a preprocessing stage and is not dynamically recomputed. 

In order to avoid spurious oscillations at discontinuities, the reconstruction operator must be nonlinear. Therefore the polynomials defined on each stencil are combined with each other and weighted 
in a nonlinear way, where the non-linearity is introduced in the WENO weights $\omega_s$ 
\begin{equation}
\tilde{\omega}_s = \frac{\lambda_s}{\left(\sigma_s + \epsilon \right)^r}, \qquad 
\omega_s = \frac{\tilde{\omega}_s}{\sum_q \tilde{\omega}_q},  
\end{equation}
through the oscillation indicators $\sigma_s$, which are computed according to \cite{shu_efficient_weno,DumbserKaeser07,DumbserKaeser06b} as
\begin{equation}
\sigma_s = \Sigma_{lm} \hat w^{n,s}_{l,i} \hat w^{n,s}_{m,i},
\end{equation}
with 
\begin{equation}
\Sigma_{lm} = \sum \limits_{ \alpha + \beta \leq M}  \, \, \int \limits_{T_e} \frac{\partial^{\alpha+\beta} \psi_l(\xi,\eta)}{\partial \xi^\alpha \partial \eta^\beta} \cdot 
                                                                         \frac{\partial^{\alpha+\beta} \psi_m(\xi,\eta)}{\partial \xi^\alpha \partial \eta^\beta} d\xi d\eta.
\end{equation}

As done in \cite{DumbserKaeser06b,DumbserKaeser07}, we set $\epsilon=10^{-14}$, $r=8$, $\lambda_s=1$ for the one-sided stencils ($s>1$) and $\lambda_1=10^5$ for the central stencil. 
The final nonlinear WENO reconstruction polynomial and its coefficients are then given by  
\begin{equation}
\label{eqn.weno} 
 \w_h(x,y,t_i^n) =  \psi_l(\xi,\eta) \hat \w^{n}_{l,i}, \qquad \textnormal{ with } \qquad  
 \hat \w^{n}_{l,i} = \sum \limits_{s=1}^{7} \omega_s \hat \w^{n,s}_{l,i}.   
\end{equation}

\subsection{Local space-time Galerkin predictor on moving triangles} 
\label{sec.lst}
In order to achieve high order of accuracy in time we use the local space-time continuous Galerkin method, where the reconstructed polynomial $\w_h$ obtained at the current element time $t_i^n$ are \textit{evolved} locally within element $T_i(t)$ until the future time $t_i^{n+1}$. This method was first introduced for the Eulerian framework in \cite{Dumbser20088209} and then extended to moving meshes in \cite{Lagrange1D,Lagrange2D,LagrangeNC,Lagrange3D}. In all the above-mentioned references the space-time continuous Galerkin procedure has been proposed \textit{locally}, i.e. the high order evolution of the reconstructed polynomial has always been carried out within each control volume and considering separately all the elements of the entire mesh. As a consequence, such a procedure automatically fits the construction of a time-accurate local time stepping algorithm.  

As previously done for the WENO reconstruction, we use again the spatial reference system $\xi-\eta$, where now the relative time $\tau$ is also considered. Therefore the physical element 
can be mapped to the reference space-time element $T_E \times [0,1]$ using the local space transformation \eqref{xietaTransf} and the following mapping in time: 
\begin{equation}
t = t_i^n + \tau \, \Delta t_i^n, \qquad  \tau = \frac{t - t_i^n}{\Delta t_i^n}. 
\label{timeTransf}
\end{equation} 
The spatial coordinate vector in physical and reference coordinates are given by $\mathbf{x}=(x,y)$ and $\boldsymbol{\xi}=(\xi,\eta)$, respectively, while $\mathbf{\tilde{x}}=(x,y,t)$ and $\boldsymbol{\tilde{\xi}}=(\xi,\eta,\tau)$ are the corresponding space-time coordinate vectors. According to \cite{Dumbser20088209}, we adopt a \textit{weak integral formulation} of the governing PDE \eqref{PDE}, which is rewritten in the space-time reference system using the relations \eqref{xietaTransf}-\eqref{timeTransf}:
\begin{equation}
\frac{\partial \Q}{\partial \tau}\tau_t + \frac{\partial \Q}{\partial \xi}\xi_t + \frac{\partial \Q}{\partial \eta}\eta_t + \frac{\partial \f}{\partial \tau}\tau_x + \frac{\partial \f}{\partial \xi}\xi_x + \frac{\partial \f}{\partial \eta}\eta_x + \frac{\partial \g}{\partial \tau}\tau_y + \frac{\partial \g}{\partial \xi}\xi_y + \frac{\partial \g}{\partial \eta}\eta_y = \mathbf{S}(\Q).  
\label{PDEweak}
\end{equation}
The Jacobian of the spatial and temporal transformation and its inverse read
\begin{equation}
J_{st} = \frac{\partial \mathbf{\tilde{x}}}{\partial \boldsymbol{\tilde{\xi}}} = \left( \begin{array}{ccc} x_{\xi} & x_{\eta} & x_{\tau} \\ y_{\xi} & y_{\eta} & y_{\tau} \\ 0 & 0 & \Delta_t \\ \end{array} \right), \quad J_{st}^{-1} = \frac{\partial \boldsymbol{\tilde{\xi}}}{\partial \mathbf{\tilde{x}}} = \left( \begin{array}{ccc} \xi_{x} & \xi_{y} & \xi_{t} \\ \eta_{x} & \eta_{y} & \eta_{t} \\ 0 & 0 & \frac{1}{\Delta t} \\ \end{array} \right), 
\label{Jac}
\end{equation}
where we used the properties $\tau_x = \tau_y = 0$ and $\tau_t = \frac{1}{\Delta t}$, according to the definition \eqref{timeTransf}. We rely on the inverse of the Jacobian matrix for reducing Eqn. \eqref{PDEweak} to
\begin{equation}
\frac{\partial \Q}{\partial \tau} + \Delta t \left( \frac{\partial \Q}{\partial \xi}\xi_t + \frac{\partial \Q}{\partial \eta}\eta_t + \frac{\partial \f}{\partial \xi}\xi_x + \frac{\partial \f}{\partial \eta}\eta_x + \frac{\partial \g}{\partial \xi}\xi_y + \frac{\partial \g}{\partial \eta}\eta_y  \right) = \Delta t \mathbf{S}(\Q),
\label{PDECG}
\end{equation}
which can be simply reformulated as
\begin{equation}
\Q_\tau =  \Delta t \P,
\label{PCG}
\end{equation}
with the aid of the term $\P$ defined as 
\begin{equation}
 \P := \S(\Q) - \left( \frac{\partial \Q}{\partial \xi}\xi_t + \frac{\partial \Q}{\partial \eta}\eta_t + \frac{\partial \f}{\partial \xi}\xi_x + \frac{\partial \f}{\partial \eta}\eta_x + \frac{\partial \g}{\partial \xi}\xi_y + \frac{\partial \g}{\partial \eta}\eta_y  \right). 
 \label{eqn.pdef}
\end{equation}    

As done in \cite{Dumbser20088209}, the solution vector $\Q$, the flux tensor $\F$, the source term $\S$ as well as the term $\P$ are discretized using a nodal finite element approach. 
The discrete solutions are denoted by $\q_h$, $\F_h$, $\S_h$ and $\P_h$, respectively, and are given by
\begin{eqnarray}
\q_h=\q_h(\xi,\eta,\tau) = \theta_{l}(\xi,\eta,\tau) \widehat{\q}_{l,i}, \qquad & \S_h=\S_h(\xi,\eta,\tau) = \theta_{l}(\xi,\eta,\tau) \widehat{\S}_{l,i},  \nonumber\\
\F_h=\F_h(\xi,\eta,\tau) = \theta_{l}(\xi,\eta,\tau) \widehat{\F}_{l,i}, \qquad & \P_h=\P_h(\xi,\eta,\tau) = \theta_{l}(\xi,\eta,\tau) \widehat{\P}_{l,i}.  
\label{thetaSol}
\end{eqnarray}
Here, $\theta_l=\theta_l(\boldsymbol{\tilde{\xi}})=\theta_l(\xi,\eta,\tau)$ are a set of space-time nodal basis functions defined by the Lagrange interpolation polynomials passing through a set of space-time nodes $\boldsymbol{\tilde{\xi}}_m=(\xi_m,\eta_m,\tau_m)$, see \cite{Dumbser20088209} for details. The same approximation also applies to the mapping from the physical space-time coordinate vector $\mathbf{\tilde{x}}$ to the reference space-time coordinate vector $\boldsymbol{\tilde{\xi}}$, hence
\begin{equation}
 \x(\xi,\eta,\tau) = \theta_l(\xi,\eta,\tau) \widehat{\x}_{l,i}, \qquad t(\xi,\eta,\tau) = \theta_l(\xi,\eta,\tau) \widehat{t}_l, 
 \label{eqn.isoparametric} 
\end{equation} 
where the use of the \textit{same} basis functions $\theta_l$ is allowed by the adoption of an \textit{isoparametric} approach. $\widehat{\mathbf{x}}_{l,i} = (\widehat{x}_{l,i},\widehat{y}_{l,i})$ are the degrees of freedom of the vector of physical coordinates in space, that are partially unknown, while $\widehat{t}_l$ denote the \textit{known} degrees of freedom of the physical time at each space-time node $\tilde{\x}_{l,i} = (\widehat{x}_{l,i}, \widehat{y}_{l,i}, \widehat{t}_l)$ according to \eqref{timeTransf}.  

In order to obtain the weak formulation of the governing PDE \eqref{PDE}, we first multiply \eqref{PDECG} with a test function which is given by the same space-time basis functions $\theta_k(\xi,\eta,\tau)$ and then we integrate it over the unit reference space-time element $T_e \times [0,1]$, i.e.  
\begin{equation}
\left\langle \theta_k,\frac{\partial \theta_l}{\partial \tau} \right\rangle \widehat{\q}_{l,i} = \Delta  t 
 \left\langle \theta_k,\theta_l \right\rangle \widehat{\P}_{l,i},
\label{LagrSTPDECG}
\end{equation} 
where the approximations given by \eqref{thetaSol} have been used as well as the following integral operator
\begin{equation}
 \left\langle f,g \right\rangle = \int \limits_{0}^{1} \int \limits_{T_e} f(\xi,\eta,\tau)g(\xi,\eta,\tau) d\xi d\eta d\tau,  
\label{intOperators}
\end{equation}
which denotes the scalar product of two functions $f$ and $g$ over the space-time reference element $T_e\times \left[0,1\right]$. Moreover the universal matrices
\begin{equation}
\K_{\tau} = \left\langle \theta_k,\frac{\partial \theta_l}{\partial \tau} \right\rangle \qquad \textnormal{and} \qquad \M = \left\langle \theta_k,\theta_l \right\rangle 
\label{Ktau}
\end{equation}
proposed in \cite{Lagrange2D,LagrangeNC} are used to write expression \eqref{LagrSTPDECG} in a more compact matrix notation, namely
\begin{equation}
\K_{\tau}\widehat{\q}_{l,i} = \Delta t \M \widehat{\P}_{l,i}. 
\label{LagrSTPDECGmatrix}
\end{equation}
Let $\widehat{\q}_{l,i}^{0}$ be the part of the degrees of freedom of vector $\widehat{\q}_{l,i}$ that are known from the initial condition $\w_h$ by setting the corresponding degrees of freedom to the known values (see \cite{Dumbser20088209} for details) and let $\widehat{\q}_{l,i}^{1}$ represent the unknown degrees of freedom for $\tau>0$. Since $\widehat{\q}_{l,i}^{0}$ are known, they can be moved onto the right-hand side of \eqref{LagrSTPDECGmatrix}, hence obtaining the following nonlinear algebraic equation system \eqref{PCG}, which can be solved by an iterative procedure, i.e.
\begin{equation}
  \K_{\tau} \widehat{\q}_{l,i}^{r+1} = \Delta t \M \widehat{\P}_{l,i}^r,  
\label{CGfinal}
\end{equation}
with the superscript $r$ denoting the iteration number. The initial guess ($r=0$) can be simply given by the reconstruction polynomial $\w_h$ at the initial time level, otherwise a more efficient initial condition based on a second order MUSCL-type scheme can be used (see \cite{HidalgoDumbser}). 

Due to the Lagrangian formulation, which implies mesh motion, we have also to consider the evolution of the vertex coordinates of the local space-time element. The motion is governed by the following ODE system 
\begin{equation}
\frac{d \mathbf{x}}{dt} = \mathbf{V}(x,y,t),
\label{ODEmesh}
\end{equation}
with the local mesh velocity $\mathbf{V}=\mathbf{V}(x,y,t)=(U,V)$ approximated again with a nodal approach as
\begin{equation}
\mathbf{V}_h=\mathbf{V}_h(\xi,\eta,\tau) = \theta_{l}(\xi,\eta,\tau) \widehat{\mathbf{V}}_{l,i}, \quad  \widehat{\mathbf{V}}_{l,i} = \mathbf{V}(\mathbf{\tilde{x}}_{l,i}).
\label{Vdof}
\end{equation}
Our algorithm belongs to the family of the so-called Arbitrary-Lagrangian-Eulerian (ALE) schemes, hence we allow the mesh velocity to be potentially different from the local fluid velocity. In this way Eulerian algorithms are reproduced by setting the mesh velocity to zero, while almost pure Lagrangian methods can be obtained when the mesh velocity coincides with the local fluid velocity. As suggested in \cite{Lagrange1D,Lagrange2D} the system \eqref{ODEmesh} can be conveniently solved for the unknown coordinate vector $\widehat{\mathbf{x}}_{l,i}$ using the same local space-time Galerkin method:
\begin{equation}
\left\langle \theta_k,\frac{\partial \theta_l}{\partial \tau} \right\rangle \widehat{\mathbf{x}}_{l,i} = \Delta t \left\langle \theta_k,\theta_l \right\rangle \widehat{\mathbf{V}}_{l,i},
\label{VCG}
\end{equation}
which yields the iteration scheme
\begin{equation}
\K_{\tau} \widehat{\mathbf{x}}^{r+1}_{l,i} = \Delta t \M \widehat{\mathbf{V}}^r_{l,i}.
\label{newVertPos}
\end{equation}
Since the physical triangle $T_i^n$ at time $t_i^n$ is known, the initial condition of the ODE system is simply given by the nodal degrees of freedom $\widehat{\x}_l$ at relative time $\tau=0$.

In practice, the ODE system \eqref{ODEmesh} is solved at each iteration of the PDE solver \eqref{CGfinal}
and the procedure is repeated until convergence is reached. 
At the end of the local space-time Galerkin procedure we obtain an \textit{element-local predictor} for the numerical solution $\q_h$, for the fluxes $\mathbf{F}_h=(\f_h,\g_h)$, for the source term $\S_h$ and also for the mesh velocity $\mathbf{V}_h$.

In a Lagrangian scheme with LTS we are dealing with hanging nodes in time and we generally do not have a matching in time of the geometry, as already explained before, but discontinuities in the 
geometry configuration are not admitted. In cell-centered Lagrangian schemes a unique node velocity is obtained by a so-called no-solver algorithm that takes as input all 
vertex-extrapolated states from the triangles in the Voronoi neighborhood surrounding the vertex. In the Lagrangian ADER-WENO schemes with global time stepping presented in \cite{Lagrange2D,LagrangeNC} 
we used a suitable node solver algorithm to update the mesh \textit{globally}, since the future time was the same for all the elements. Here, in the context of LTS, we adopt again the node solver  algorithm with the aim to fix a \textit{unique} node velocity, but the vertex will be \textit{physically} (and not virtually) moved only when an element of the Voronoi neighborhood of the vertex 
fulfills the update criterion \eqref{eqn.WENOupdate}. To handle this situation in practice, each node $k$ is also equipped with a local node time $t_k^n$.  

\subsection{Mesh motion with local time stepping}
\label{sec.meshMot}
As explained at the end of Section \ref{sec.lst}, each node $k$ of the computational mesh needs to be assigned a \textit{uniquely} defined velocity vector. The Voronoi neighborhood $\mathcal{V}_k$ 
of node $k$ is composed by all those elements $T_j$ which share the node $k$. The node $k$ will be moved each time the update criterion \eqref{eqn.WENOupdate} is satisfied by one element 
$T_i\in \mathcal{V}_k$. Therefore the future time to which node $k$ moves will coincide with the future time $t_i^{n+1}$ of that element $T_i$. 

In \cite{LagrangeMHD} three different node solver algorithms have been presented and here we consider the node solver denoted as $\mathcal{NS}_{cs}$, which adopts the idea of Cheng and Shu 
\cite{chengshu1,chengshu2}. However, rather than taking a simple arithmetic average of the velocity, the node velocity $\overline{\mathbf{V}}_k$ is computed as a \textit{mass weighted} average 
velocity among the neighborhood $\mathcal{V}_k$ of node $k$, i.e. 
\begin{equation}
\overline{\mathbf{V}}_k = \frac{1}{\mu_k}\sum \limits_{T_j \in \mathcal{V}_k}{\mu_{k,j}\overline{\mathbf{V}}_{k,j}}, 
\label{eqnNScs}
\end{equation}
with
\begin{equation}
\mu_k = \sum \limits_{T_j \in \mathcal{V}_k}{\mu_{k,j}}, \qquad \mu_{k,j}=\rho^{n}_j |T_j^{n}|.
\label{eqn.NScs.weights}
\end{equation}
The local weights $\mu_{k,j}$ are the masses of the elements $T_j$, obtained by multiplying the cell averages of the density $\rho_j$ with the cell area $|T_j|$ at the current neighbor time level $t_j^n$.

The mesh motion plays an important role in Lagrangian schemes, because it allows interfaces and shear waves to be precisely identified.  For this reason an accurate computation of  the node velocity represents a crucial step, and
in our approach the local velocity contributions $\overline{\mathbf{V}}_{k,j}$ are taken to be the time integrals of the high order vertex-extrapolated velocities at node $k$. We can use the space-time reference system $\xi-\eta-\tau$ and the velocity approximation given by \eqref{Vdof} to evaluate the time integral. Since each node $k$ can be moved by any of the Voronoi neighbors $T_j$, the vertex time level of node $k$ is not known \textit{a priori} when an element $T_i$ satisfies \eqref{eqn.WENOupdate} and is ready to update the geometry. Therefore, it is much more convenient to define a \textit{node time} variable $t_k^n$, that is independent of the time evolution of the elements and advances in time whenever the node is moved by any of its Voronoi neighbors $T_j$. As a result, the high order velocity integration for each element $T_j \in \mathcal{V}_k$ must be done within the time interval $\Delta t_{k} = [t_k^n,t_k^{n+1}]$, that has to be \textit{rescaled} to the corresponding reference time interval $\Delta \tau_k=[\tau_{k,j}^0,\tau_{k,j}^1]$ as 
\begin{equation}
\tau_{k,j}^0 = \frac{t_k^n-t_j^n}{\Delta t_j^n} \qquad \tau_{k,j}^1 = \frac{t_k^{n+1}-t_j^n}{\Delta t_j^n}, \qquad \forall T_j \in \mathcal{V}_k, 
\label{eqn.dtau}
\end{equation}
where $\Delta t_j^n$ is the local timestep of element $T_j$. Recall that $t_k^{n+1}=t_i^{n+1}$, if the node is moved by element $T_i$ which is supposed to satisfy the update criterion. 
Finally the local velocity contributions $\overline{\mathbf{V}}_{k,j}$ are given by 
\begin{equation}
\overline{\mathbf{V}}_{k,j} = \left( \int \limits_{\tau_{k,j}^0}^{\tau_{k,j}^1} \theta_l(\xi^e_{m(k)}, \eta^e_{m(k)}, \tau) d \tau \right) \widehat{\mathbf{V}}_{l,j}, 
\label{NodesVel}
\end{equation} 
where $m(k)$ is a mapping from the global node number $k$ to the local node number in element $T_j$, while $\xi^e_{m}$ and $\eta^e_{m}$ represent the coordinates of the vertices of the reference triangle in space. $\widehat{\mathbf{V}}_{l,j}$ are the space-time degrees of freedom which are \textit{known} from the local space-time predictor solution $\q_{h,j}$. 
Each node $k$ belonging to element $T_i$ is finally moved to the new position $\mathbf{X}^{n+1}_k$ with
\begin{equation}
\mathbf{X}^{n+1}_{k}	= \mathbf{X}^{n}_{k}	+ \Delta t_k \, \overline{\mathbf{V}}_k.
\label{eqn.newVertex}
\end{equation}

\subsection{Finite volume scheme}
\label{sec.SolAlg}

\begin{figure}[!htbp]
\begin{center}
\includegraphics[width=0.80\textwidth]{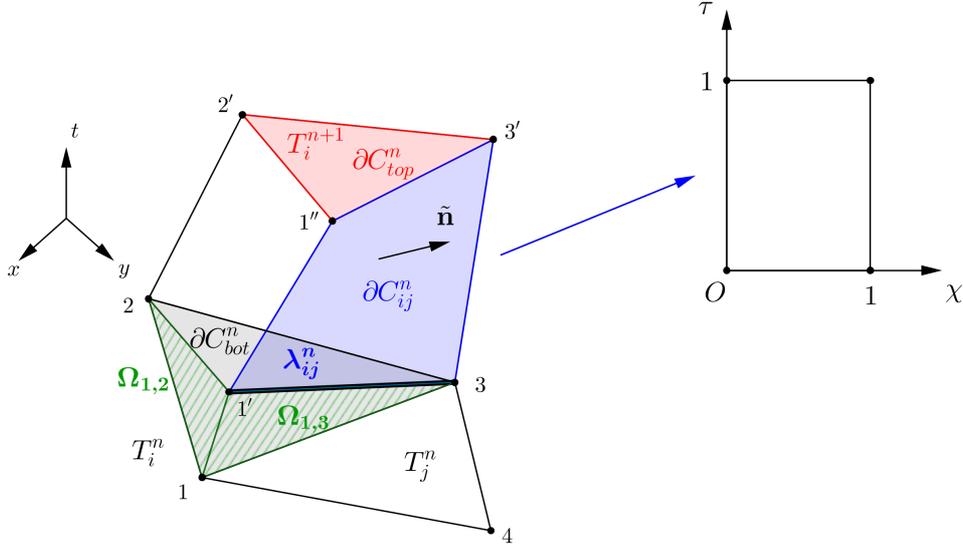} 
\caption{Space-time evolution of element $T_i$ from time $t_i^n$ (black triangle) to time $t_i^{n+1}$ (red triangle). The triangular sub-surfaces $\Omega_{1,2}$ and $\Omega_{2,3}$ (already 
computed in the past by some Voronoi neighbors of the vertices of $T_i$) are highlighted in green, while the trapezoidal space-time sub-surfaces $\partial C_{ij}^n$ computed with the current
element update are highlighted in blue.} 
\label{fig.FVscheme}
\end{center}
\end{figure}

The vector of conserved variables $\Q_i^n$ is evolved to the next time level $t_i^{n+1}$ only when element $T_i$ obeys the update criterion \eqref{eqn.WENOupdate}. As proposed in 
\cite{Lagrange2D,Lagrange3D} the governing PDE \eqref{PDE} can be rewritten in a more compact space-time divergence form, which reads 
\begin{equation}
\tilde \nabla \cdot \tilde{\F} = \mathbf{S}(\Q), 
\label{PDEdiv3D}
\end{equation} 
with the space-time nabla operator and the tensor $\tilde{\F}$ defined as
\begin{equation}
\tilde \nabla  = \left( \frac{\partial}{\partial x}, \, \frac{\partial}{\partial y}, \, \frac{\partial}{\partial t} \right)^T,  \qquad 
\tilde{\F}  = \left( \mathbf{F}, \, \Q \right) = \left( \mathbf{f}, \, \mathbf{g}, \, \Q \right).
\end{equation} 
The conservation law \eqref{PDEdiv3D} is then integrated in space and time over the space-time control volume $C^n_i = T_i(t) \times \left[t_i^{n}; t_i^{n+1}\right]$ generated by the time evolution of element $T_i$ and depicted in Figure \ref{fig.FVscheme}, hence yielding
\begin{equation}
\int\limits_{t_i^{n}}^{t_i^{n+1}} \int \limits_{T_i(t)} \tilde \nabla \cdot \tilde{\F} \, d\mathbf{x}\, dt = \int\limits_{t_i^{n}}^{t_i^{n+1}} \int \limits_{T_i(t)} \S \, d\mathbf{x}\, dt,   
\label{STPDE}
\end{equation}
which, after application of Gauss' theorem, reads 
\begin{equation}
\int \limits_{\partial C^n_i} \tilde{\F} \cdot \ \mathbf{\tilde n} \, \, dS = 
\int\limits_{t_i^{n}}^{t_i^{n+1}} \int \limits_{T_i(t)} \S \, d\mathbf{x}\, dt,   
\label{I1}
\end{equation}   
The vector $\mathbf{\tilde n} = (\tilde n_x,\tilde n_y,\tilde n_t)$ is the outward pointing space-time unit normal vector defined on the space-time surface $\partial C^n_i$, which is composed of five space-time sub-surfaces, as shown in Figure \ref{fig.FVscheme}: the first one $\partial C_{bot}^n$ is given by the element configuration $T_i^{n}$ at the current time level, while $\partial C_{top}^n$ represents the control volume $T_i^{n+1}$ evolved to the future time level. The remaining three lateral space-time sub-surfaces $\partial C_{ij}^n$ are usually shared with the so-called \textit{Neumann neighbors} $\mathcal{N}_i$ of $T_i$, i.e. with the direct side neighbors. As explained in \cite{Lagrange2D,LagrangeNC,LagrangeMHD} a set of bilinear basis functions are used to parametrize the lateral 
sub-surfaces, which are mapped onto a side-aligned local reference system $(\chi,\tau)$. The unit normal vector $\mathbf{\tilde n}$ can be computed from the parametrization of the lateral sub-surfaces,  while for $\partial C_{bot}^n$ and $\partial C_{top}^n$ it simply reads $\mathbf{\tilde n}=(0,0,-1)$ and $\mathbf{\tilde n}=(0,0,1)$, respectively. 
%Note that element $T_i^n$ may be defined by three verteces which are \textit{not} at the same time level, as shown in Figure \ref{fig.FVscheme}. Such a configuration might happen if one of the Voronoi neighbor elements of the verteces of $T_i$ has previuosly done the update timestep: in this case its nodes have moved according to \eqref{eqn.newVertex}, hence leading to a gap between the time level of the nodes and the time level of those elments that has not done the update yet. 

In the time-accurate local time stepping (LTS) algorithm, when the element $T_i$ is ready to update its numerical solution $\Q_i^{n}$, it might well be the case that the vertices of $T_i$ have already been moved by another element $T_j$ sharing one or more nodes with $T_i$. This situation generates hanging nodes in time, as shown in Figure \ref{fig.FVscheme}, where vertex $1$ has changed its position to $1^\prime$. In order to design a suitable finite volume scheme on moving meshes with LTS, some parts of the flux integral appearing in \eqref{I1} will be computed using a \textit{memory variable} $\Q_i^M$, according to \cite{ALELTS1D}. The memory variable contains all fluxes through the element space-time sub-surfaces $\partial C_{ij}^n$ in the \textit{past}, e.g. the fluxes through the space-time  triangular surfaces $\Omega_{1,2}$ and $\Omega_{1,3}$ depicted in Figure \ref{fig.FVscheme}. Therefore, from \eqref{I1} the following high order ALE one-step finite volume scheme with LTS is obtained: 
\begin{equation}
|T_i^{n+1}| \, \Q_i^{n+1} = |T_i^n| \, \Q_i^n - \sum \limits_{T_j \in \mathcal{N}_i} \,\, {\int \limits_0^1 
\int \limits_{0}^{1} | \partial C_{ij}^n| \tilde{\F}_{ij} \cdot \mathbf{\tilde n}_{ij} \, d \tau d\chi}
+ \int\limits_{t_i^{n}}^{t_i^{n+1}} \int \limits_{T_i(t)} \S(\mathbf{q}_h) \, d\mathbf{x} dt 
 + \Q_i^M\,,
\label{PDEfinal}
\end{equation}
with $|T_i^{n}|$ and $|T_i^{n+1}|$ representing the surface of triangle $T_i$ at the current and at the future time level, i.e. $t_i^n$ and $t_i^{n+1}$, and $|\partial C_{ij}^n|$ denoting the 
determinant of the coordinate transformation of each lateral sub-surface $\partial C^n_{ij}$. Furthermore $\tilde{\F}_{ij} \cdot \mathbf{\tilde n}_{ij}$ is the numerical flux used to resolve 
the discontinuity of the predictor solution $\mathbf{q}_h$ at the space-time sub-face $\partial C_{ij}^n$. In the finite volume scheme \eqref{PDEfinal} the flux integral across the 
quadrilateral sub-surface $\partial C_{ij}^n$ is computed in an edge-based unit reference system $(\chi,\tau) \in [0,1]^2$ that is linked to the physical coordinates of the four space-time 
nodes that define $\partial C_{ij}^n$. Note that in the edge-aligned system the relative time coordinate $\tau$ is in general \textit{different} from the ones in the adjacent left and right
elements $T_i$ and $T_j$, respectively, since the two nodes that define the edge may have already been moved before the update of element $T_i$. 
Let us denote the common edge between element $T_i$ and $T_j \in \mathcal{N}_i$ with $\lambda_{ij}$ and the global number of the first node on $\lambda_{ij}$ with $L$ and the one of the 
second node on the same edge with $R$, then the space-time coordiantes of the four space-time nodes defining the sub-surface $\partial C_{ij}^n$ in \eqref{PDEfinal} are given by 
\begin{equation} 
  \tilde{\mathbf{x}}_{ij}^1 = \left( \mathbf{X}^n_{L}, t_{L}^n \right), \qquad 
  \tilde{\mathbf{x}}_{ij}^2 = \left( \mathbf{X}^n_{R}, t_{R}^n \right), \qquad 	
  \tilde{\mathbf{x}}_{ij}^3 = \left( \mathbf{X}^{n+1}_{R}, t_{R}^{n+1} \right), \qquad 
  \tilde{\mathbf{x}}_{ij}^4 = \left( \mathbf{X}^{n+1}_{L}, t_{L}^{n+1} \right).  	
	\label{eqn.spacetimenodes} 
\end{equation} 
Note that $L=L(i,j)$ and $R=R(i,j)$ are functions of the numbers of element $T_i$ and the neighbor $T_j$, respectively, but to ease notation this explicit dependency is dropped. 
The associated space-time integral of the numerical flux over $\partial C_{ij}^n$ is also called \textit{edge flux} and denoted by $\mathbf{G}_{ij}^n$ in the following. The physical 
times of the four space-time nodes \eqref{eqn.spacetimenodes} have then to be rescaled to each individual reference space-time coordinate system associated with element $T_i$ and 
its neighbor $T_j$, respectively, using the time transformation \eqref{timeTransf}.  

In order to obtain a conservative scheme, the task of the memory variable $\Q_i^M$ in \eqref{PDEfinal} is to accumulate (sum) all past fluxes through the lateral space-time 
sub-surfaces, from the current element time $t_i^n$ to the current local node times $t_{L}^n$ and $t_R^n$, respectively, see \cite{ALELTS1D}. The edge flux $\mathbf{G}_{ij}^n$ 
through the sub-surface $\partial C_{ij}^n$ is given by  
\begin{equation}
  \mathbf{G}_{ij}^n = \int \limits_{\partial C_{ij}^n} \tilde{\F}_{ij} \cdot \mathbf{\tilde n}_{ij} \, d{S} 
	                  = {\int \limits_0^1 \int \limits_{0}^{1} | \partial C_{ij}^n| \tilde{\F}_{ij} \cdot \mathbf{\tilde n}_{ij} \, d \tau d\chi}. 
\end{equation} 
Then, if element $T_i$ is updated according to \eqref{PDEfinal}, the memory variable of the element itself is reset to zero and the memory variables of the \textit{neighbor} elements 
$T_j \in \mathcal{N}_i$ are updated by summing (accumulating) the contribution of the edge-flux $\mathbf{G}_{ij}^n$ to $\Q_j^M$. Note that for element $T_i$ the contribution 
$\mathbf{G}_{ij}^n$ has negative sign. Like in the 1D case presented in \cite{ALELTS1D} we therefore have after each update of element $T_i$:  
\begin{equation}
 \Q_i^M:=0, \qquad \Q_j^M:=\Q_j^M + \mathbf{G}_{ij}^n, \qquad \forall T_j \in \mathcal{N}_i. 
\label{eqn.Qm2}
\end{equation} 
%
%Following this strategy, when one neighbor element $T_j$  satisfies the update criterion and is allowed to advance it time with its local timestep, the flux contribution through the sub-surface $\Omega_{i,j}$ up to time $t_i^{n+1}$ will already have been computed by cell $T_i$ with \eqref{eqn.Qm2}. Thanks to this procedure, which was first proposed in \cite{dumbserkaeser06d}, the numerical flux across the surface given by the space-time evolution of the element edge $\lambda_{ij}$ is computed only \textit{once} by the first element, between $T_i$ and $T_j$, that is scheduled for the update. 
%In classical GTS schemes we typically compute the flux twice, when both element $T_i$ and element $T_j$ update the numerical solution according to the numerical scheme \eqref{PDEfinal}, where the memory variable disappears. 

The implementation of the finite volume scheme \eqref{PDEfinal} requires that a numerical flux is specified through an approximate Riemann solver. A possible simple formulation for the numerical 
flux is given by the Rusanov-type ALE flux, which, according to \cite{Lagrange2D}, reads 
\begin{equation}
  \tilde{\F}_{ij} \cdot \mathbf{\tilde n}_{ij} =  
  \frac{1}{2} \left( \tilde{\F}(\q_h^+) + \tilde{\F}(\q_h^-)  \right) \cdot \mathbf{\tilde n}_{ij}  - 
  \frac{1}{2} s_{\max} \left( \q_h^+ - \q_h^- \right),  
  \label{eqn.rusanov} 
\end{equation} 
where $s_{\max}$ is the maximum eigenvalue of the ALE Jacobian matrix w.r.t. the normal direction in space, which is 
\begin{equation} 
\mathbf{A}^{\!\! \mathbf{V}}_{\mathbf{n}}(\Q)=\left(\sqrt{\tilde n_x^2 + \tilde n_y^2}\right)\left[\frac{\partial \mathbf{F}}{\partial \Q} \cdot \mathbf{n}  - 
(\mathbf{V} \cdot \mathbf{n}) \,  \mathbf{I}\right], \qquad    
\mathbf{n} = \frac{(\tilde n_x, \tilde n_y)^T}{\sqrt{\tilde n_x^2 + \tilde n_y^2}},  
\end{equation} 
with $\mathbf{I}$ representing the identity matrix and $\mathbf{V} \cdot \mathbf{n}$ denoting the local normal mesh velocity. 

A more sophisticated alternative is given by  the Osher-type numerical flux, which guarantees a less dissipative numerical scheme if compared with the Rusanov flux. It has been presented in \cite{OsherUniversal} for the Eulerian case and then extended to moving meshes in multiple space dimensions in \cite{Lagrange1D,Lagrange2D,Lagrange3D}. The corresponding numerical flux is given by
\begin{equation}
  \tilde{\F}_{ij} \cdot \mathbf{\tilde n}_{ij} =  
  \frac{1}{2} \left( \tilde{\F}(\q_h^+) + \tilde{\F}(\q_h^-)  \right) \cdot \mathbf{\tilde n}_{ij}  - 
  \frac{1}{2} \left( \int \limits_0^1 \left| \mathbf{A}^{\!\! \mathbf{V}}_{\mathbf{n}}(\boldsymbol{\Psi}(s)) \right| ds \right) \left( \q_h^+ - \q_h^- \right),  
  \label{eqn.osher} 
\end{equation} 
where a simple straight-line segment path is used to connect the left and the right state across the discontinuity, i.e.   
\begin{equation}
\boldsymbol{\Psi}(s) = \q_h^- + s \left( \q_h^+ - \q_h^- \right), \qquad 0 \leq s \leq 1.  
\label{eqn.path} 
\end{equation} 
According to \cite{OsherUniversal} the integral in \eqref{eqn.osher} is evaluated numerically using Gaussian quadrature. 
The absolute value of the dissipation matrix in \eqref{eqn.osher} is evaluated as usual as 
\begin{equation}
 |\mathbf{A}| = \mathbf{R} |\boldsymbol{\Lambda}| \mathbf{R}^{-1},  \qquad |\boldsymbol{\Lambda}| = \textnormal{diag}\left( |\lambda_1|, |\lambda_2|, ..., |\lambda_\nu| \right),  
\end{equation}
where $\mathbf{R}$ and $\mathbf{R}^{-1}$ denote the right eigenvector matrix and its inverse, respectively.

\begin{figure}[!htbp]
\begin{center}
\includegraphics[width=0.70\textwidth]{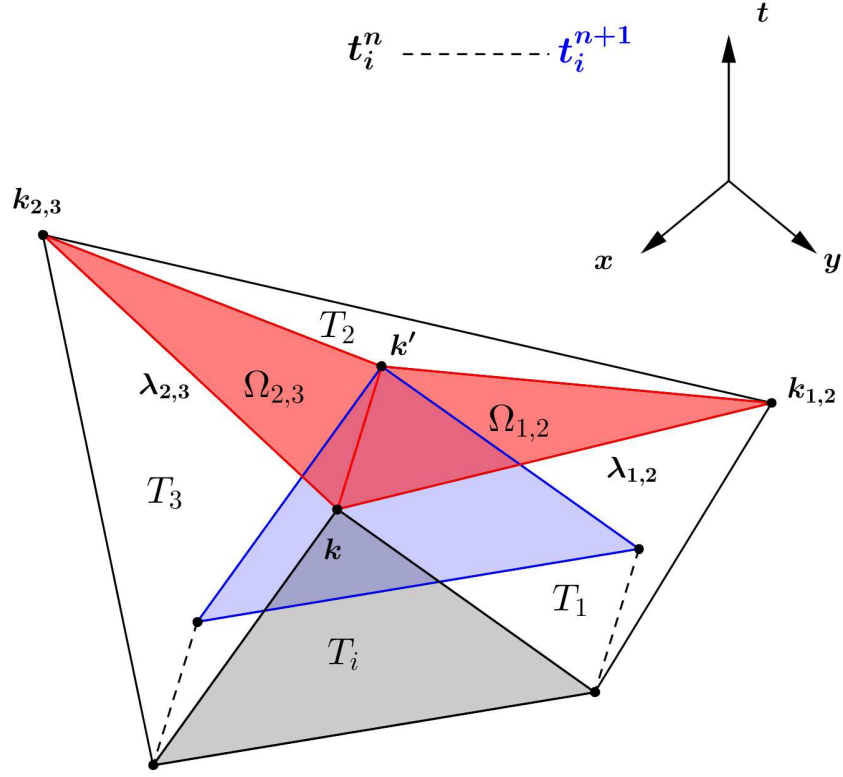} 
\caption{Space-time evolution of element $T_i$ from time $t_i^n$ (black triangle) to time $t_i^{n+1}$ (blue triangle). The triangular sub-surfaces $\Omega_{1,2}$ and $\Omega_{2,3}$ are highlighted in red.} 
\label{fig.STint}
\end{center}
\end{figure}

When element $T_i$ performs its local time update, the geometry of cell $T_i$ is also updated, because all three vertices of $T_i$ are moved according to \eqref{eqn.newVertex}. Using the 
memory variable $\Q_i^M$ we ensure conservation of the edge-fluxes, since the numerical fluxes computed over the space-time sub-surfaces $\partial C_{ij}$ are immediately saved 
(with opposite sign) in the memory variables of the neighbor elements $T_j \in \mathcal{N}_i$. While the consideration of edge fluxes is sufficient for the Lagrangian LTS algorithm 
presented in \cite{ALELTS1D}, its extension to moving unstructured triangular meshes requires an important modification due to the increased topological complexity of a two-dimensional 
mesh. 
As shown in Figure \ref{fig.STint}, each vertex $k$ of element $T_i$ is shared among the Voronoi neighbors $T_j \in \mathcal{V}_k$. Hence, we must also compute a numerical 
flux $\mathbf{G}_{k,m}$ across each edge defined by the vertices $k$ and $m$ which does \textit{not} belong to element $T_i$, i.e. 
\begin{equation}
 \mathbf{G}_{k,m} = {\int \limits_{\partial \Omega_{k,m}} \tilde{\F}_{{l},{r}} \cdot \mathbf{\tilde n}_{{l},{r}} \, d\tilde{\x}}. 
\label{eqn.STint}
\end{equation}
This \textit{vertex flux} will also be stored (with the proper sign) in the corresponding memory variables $\Q_{{l}}^M$ and $\Q_{{r}}^M$ of elements 
$T_{{l}}$ and $T_{{r}}$, where ${{l}}$ denotes the left element and ${{l}}$ denotes the right element on the corresponding edge composed of vertices 
$k-m$, respectively. 
As shown in Figure \ref{fig.STint},  the numerical flux is integrated over the \textit{triangular} space-time surfaces $\Omega_{j,j+1}$, defined by 
vertices $\left( \tilde{\x}(k),\tilde{\x}(k'),\tilde{\x}(m) \right)$, that represent the space-time coordinates of vertex $k$ at the old and at the 
new time level, and the space-time location of vertex $k_{i,j}$, respectively. 

In our finite volume formulation we are carrying out an integration over the closed space-time control volume $C_{i}^n$, which automatically guarantees the compliance with the geometric conservation law 
(GCL), see the appendix of \cite{Lagrange3D} for more details. From the Gauss theorem one has indeed    
\begin{equation}
 \int_{\partial \mathcal{C}_i^n} \mathbf{\tilde n} \, dS = 0. 
 \label{eqn.gcl} 
\end{equation} 
In order to verify whether the GCL is also satisfied in the practical implementation of our Lagrangian LTS algorithm, we need to compute the integral above whenever element $T_i$ performs an update.  For this purpose, we also compute a variable $H_i^M$ that behaves like the memory variable $\Q_i^M$, but for the GCL. All past contributions to the integral \eqref{eqn.gcl} relative to the cell 
$T_i$ are recorded in the \textit{geometrical memory variable} $H_i^M$, which is reset to zero when the local timestep procedure has been completed by element $T_i$. Strictly speaking this this is 
not needed, since Eqn. \eqref{eqn.gcl} is always satisfied at the end of a local time step because the final space-time control volume is always closed! In all test problems reported in Section 
\ref{sec.validation}, property \eqref{eqn.gcl} has always been explicitly verified for each element and for each local time step up to machine precision. 

\subsection{Description of the high order Lagrangian LTS algorithm in multiple space dimensions}
The aim of this Section is to give an overall overview of the entire LTS algorithm that has been previously described in all its parts. By placing each portion of the algorithm in a context, this  presentation should clarify how the numerical scheme can be practically implemented.
Due to the LTS approach, where elements are updated in the order given by the update criterion \eqref{eqn.WENOupdate}, we can no longer speak of \textit{timesteps} but we have to consider \textit{cycles}, as done in \cite{ALELTS1D}. In each cycle the scheme runs over all elements and only those which obey condition \eqref{eqn.WENOupdate} are allowed to update the numerical solution, while the others are simply skipped to the next cycle.
\begin{figure}[!htbp]
\begin{center}
\begin{tabular}{cc} 
\includegraphics[width=0.50\textwidth]{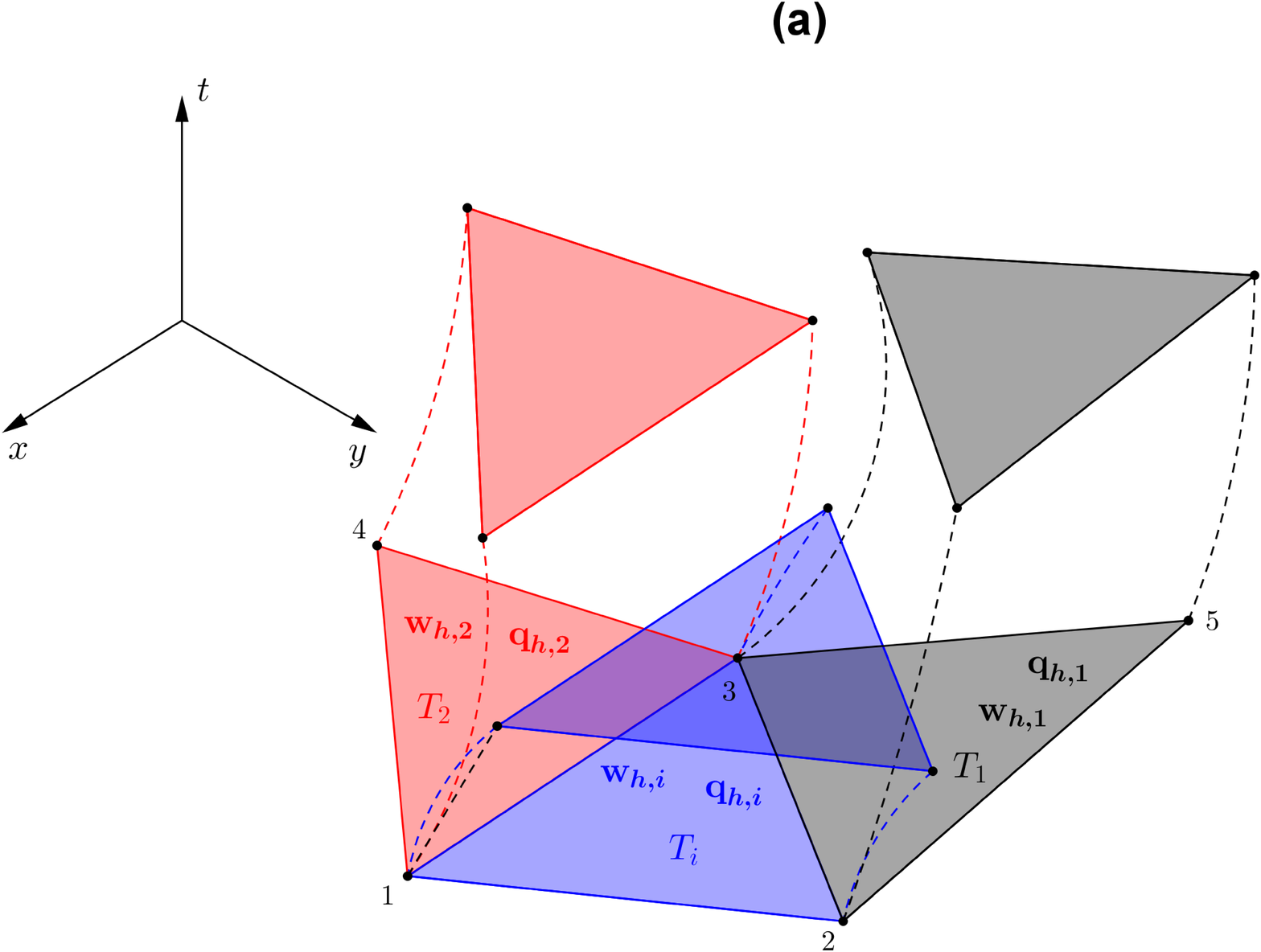}  &           
\includegraphics[width=0.50\textwidth]{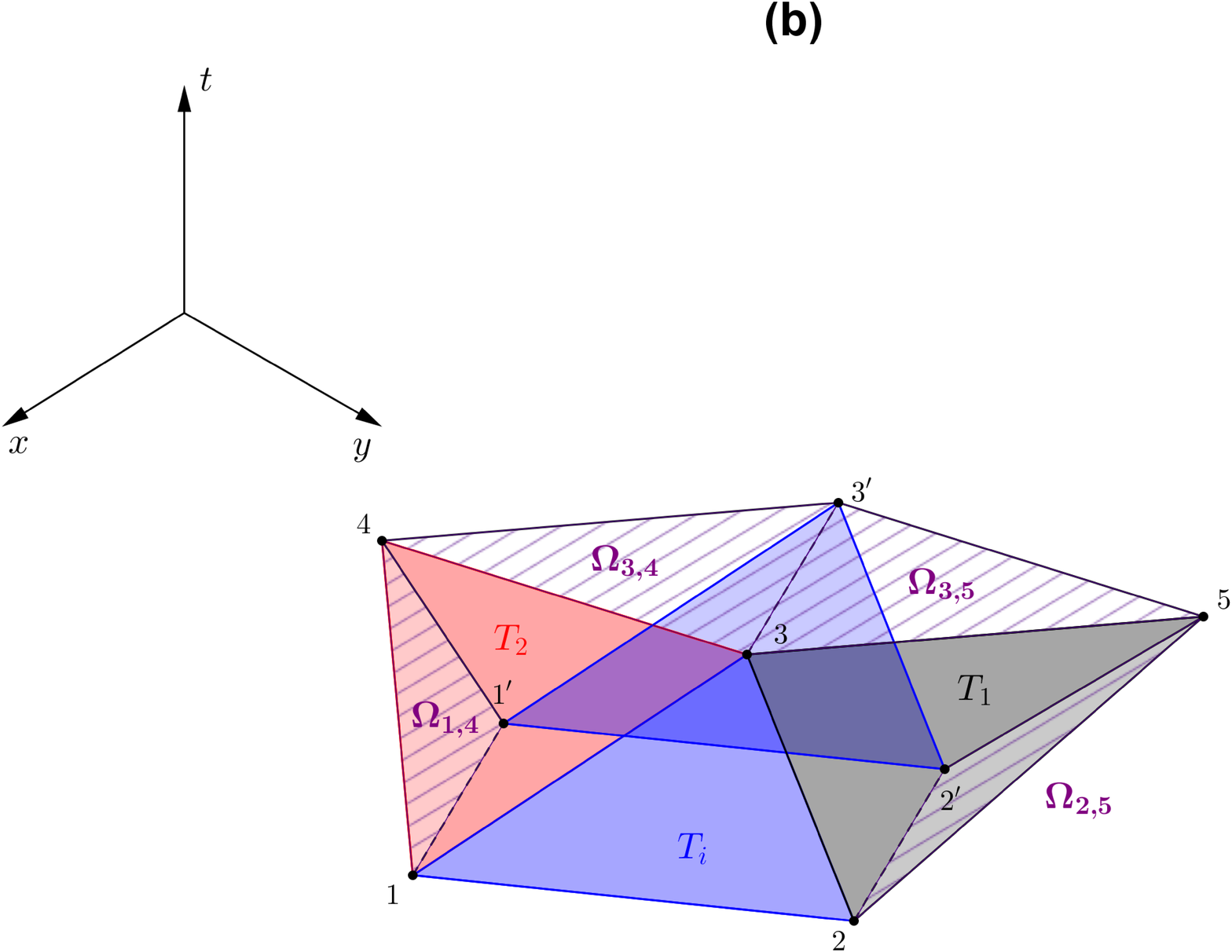}  \\
\includegraphics[width=0.50\textwidth]{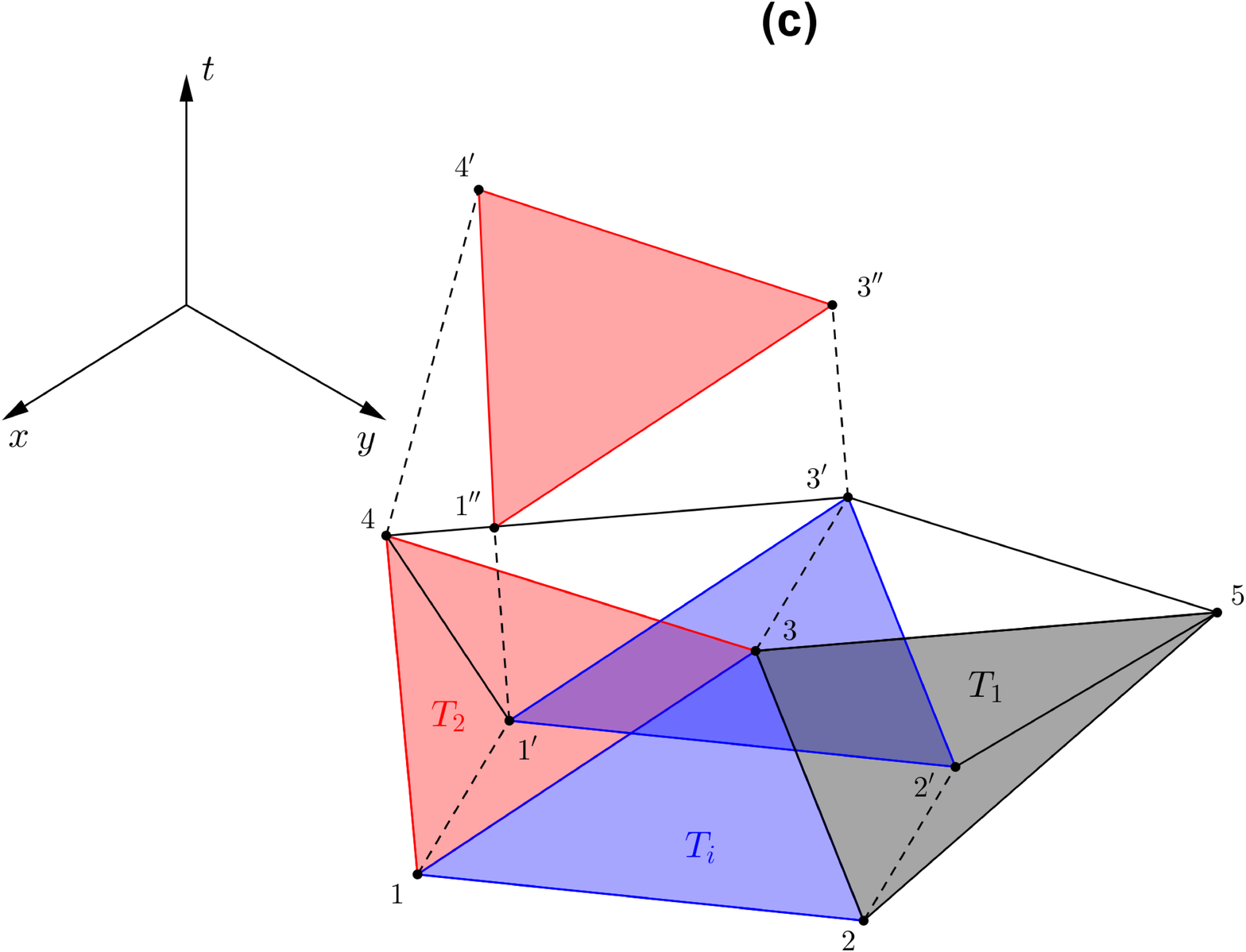}  &           
\includegraphics[width=0.50\textwidth]{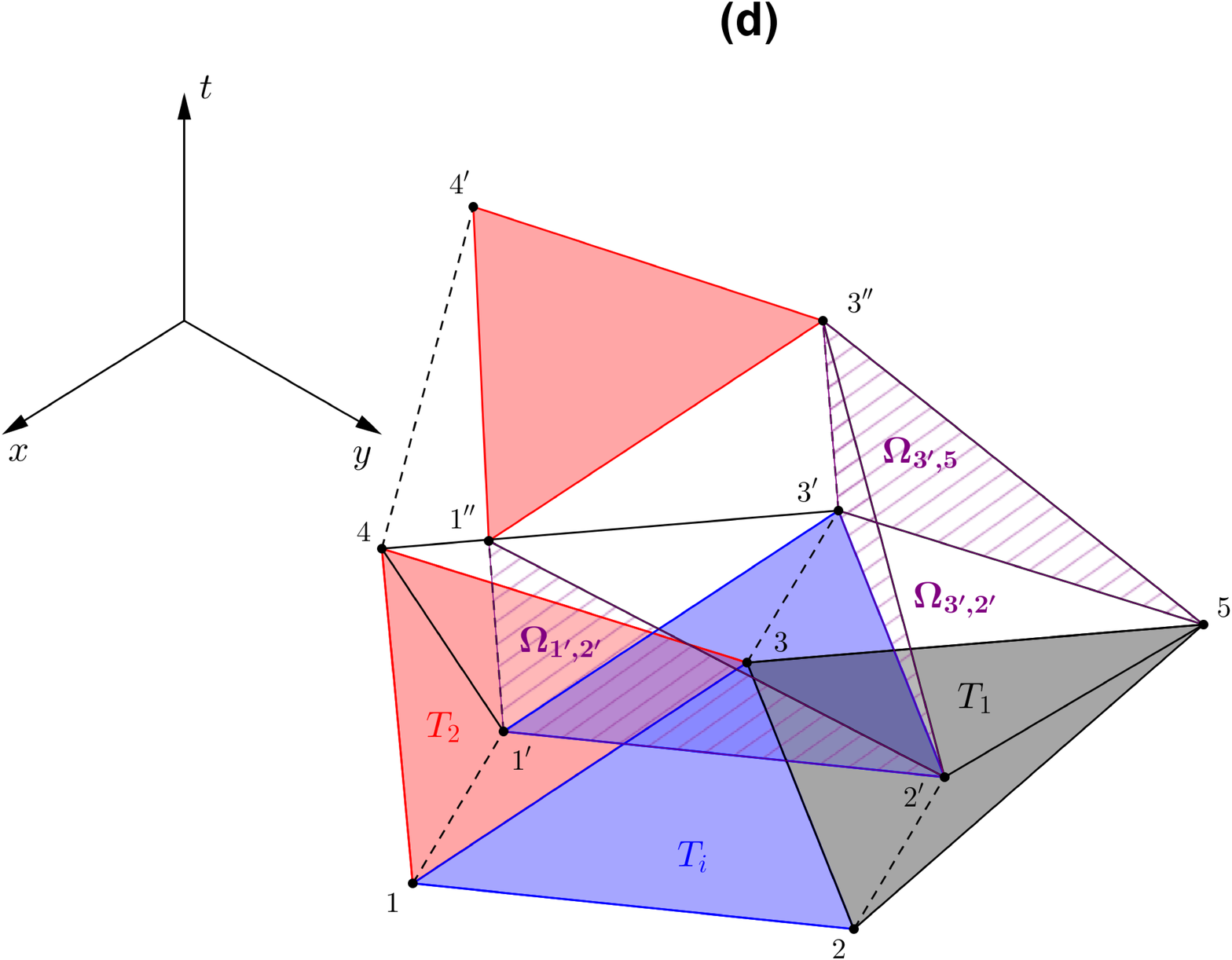}  \\
\end{tabular} 
\caption{Update of element $T_i$ and $T_2$ according to the high order Lagrangian LTS algorithm presented in this paper. At the beginning we assume the \textit{same} current time for each element, i.e. 
$t_i^n=t_1^n=t_2^n=t$. $\mathbf{(a)}$ At the current time level $t$ each element is given its own reconstruction and predictor solution $\w_h$ and $\q_h$, respectively. $\mathbf{(b)}$ Update of element $T_i$ to the new time level $t_i^{n+1}$. Computation of the necessary edge fluxes with the direct neighbors \textit{and} computation of the associated vertex fluxes $\Omega_{1,4},\Omega_{3,4},\Omega_{3,5},\Omega_{2,5}$. $\mathbf{(c)}$ Update of element $T_2$, where the edge fluxes are evaluated only over the space-time surfaces that exceeds the vertex fluxes previously calculated and stored in the memory variable $\Q_2^M$. $\mathbf{(d)}$ Computation of the vertex fluxes related to the update of element $T_2$.} 
\label{fig.LTSalgorithm}
\end{center}
\end{figure}
In the \textit{pre-processing phase} all elements of the mesh are assigned with the initial condition of the problem at the common time level $t=0$, i.e. the cell averages $\Q_i^n$ are defined 
according to \eqref{eqn.cellaverage} from the known initial condition. For each element the \textit{first} WENO reconstruction procedure presented in Section \ref{sec.weno} is carried out. Since
all elements are at the same time $t=0$, for this first reconstruction no virtual geometry or virtual cell averages $\tilde{\Q}$ are needed. As a result, we obtain the high order spatial polynomial 
$\w_h$ for each element. Then, the \textit{element-local timestep} $\Delta t_i^n$ is computed for each cell $T_i$ according to a classical CFL stability condition, considering only cell number $i$ 
and its Neumann neighborhood $\mathcal{N}_i$, i.e.    
\begin{equation}
\Delta t_i^n = \min \left( \textnormal{CFL} \,\frac{\tilde{d}_i}{|\tilde{\lambda}_{\max,i}|},\textnormal{CFL} \,\frac{\tilde{d}_j}{|\tilde{\lambda}_{\max,j}|} \right), \qquad \forall T_j \in \mathcal{N}_i, 
\label{eqn.timestep}
\end{equation}
with $\tilde{d}_j=d_j^0$ denoting the incircle diameter of element $T_j$ and $|\tilde{\lambda}_{\max,j}|=|\lambda_{\max,j}|^0$ representing the maximum absolute value of the eigenvalues computed 
from the initial condition $\tilde{\Q}_j=\Q_j^0$ in $T_j$. $\textnormal{CFL}$ is the Courant-Friedrichs-Levy number that must satisfy the inequality $\textnormal{CFL} \leq 0.5$ in the two-dimensional  case, as stated in \cite{ToroBook}. In the last part of the pre-processing stage, since the local element timestep $\Delta t^n_i$ as well as the local reconstruction polynomial $\w_h$ have already 
been computed, we are able to carry out the local space-time Galerkin predictor procedure described in Section \ref{sec.lst}, which gives the high order local space-time predictor solution $\q_h$. 
All cells are now at the same current time level $t=0$ and for each element $T_i$ the local predictor solution $\q_{h}$, the local reconstruction polynomial $\w_{h}$ and the cell average $\Q_i^n$ 
are given (Figure \ref{fig.LTSalgorithm} (a)). We underline that also each node $k$ of the entire computational mesh is assigned the initial time level $t_k^0=0$. 

The algorithm proceeds with the \textit{computational phase}, during which each element $T_i$ will reach the imposed final time of the simulation $t=t_f$ in a certain number of necessary cycles, 
according to its own optimal timestep. The first cycle starts by looping over all elements to check in which elements the update criterion \eqref{eqn.WENOupdate} is satisfied. If an element $T_i$  
obeys condition \eqref{eqn.WENOupdate}, then it performs the local timestep until its future time $t_i^{n+1} = t_i^n+\Delta t _i^n$ (Figure \ref{fig.LTSalgorithm} (b)) through the following 
sub-steps: 
\begin{itemize}
	\item \textit{mesh motion}: each vertex $k$ of element $T_i$ is moved to the new position at time $t_k^{n+1}=t_i^{n+1}$ using the node solver algorithm illustrated in Section \ref{sec.meshMot} and all other geometric quantities of element $T_i$ are also updated;
	\item \textit{edge flux computation}: we compute the numerical fluxes $\mathbf{G}_{ij}^n$ through the quadrilateral space-time sub-surfaces and using the high order Lagrangian finite volume 
	scheme \eqref{PDEfinal} we obtain the numerical solution $\Q_i^{n+1}$. Subsequently, we reset the memory variable of element $T_i$ to zero, i.e. $\Q_i^M:=0$ and accumulate the edge-fluxes into
	the memory variables of the neighbor elements to maintain conservation ($\Q_j^M:=\Q_j^M+\mathbf{G}_{ij}^n$). Also the geometry variable $H_i^M$ is reset to zero, after assuring that condition 
	\eqref{eqn.gcl} is satisfied; 
	\item \textit{vertex flux computation}: as explained in Section \ref{sec.SolAlg}, for each vertex $k$ of the element $T_i$ we also need to evaluate for each edge $k-k_{i,j}$ the additional fluxes $\mathbf{G}_{k,j}$ using \eqref{eqn.STint} (Figure \ref{fig.LTSalgorithm} (b)). The numerical fluxes evaluated over the space-time triangular sub-surface $\Omega_{j,j+1}$ (see Figure \ref{fig.STint}) are immediately stored into the memory variable of the adjacent elements $T_j,T_{j+1}$, while the part of the geometry integral \eqref{eqn.gcl} is stored into $H_j^M$ and $H_{j+1}^M$. In this way we ensure that the numerical scheme is fully conservative; 
	\item \textit{virtual projection}: all the elements $T_j$ belonging to the entire reconstruction stencil $\S_i^W$ of element $T_i$ are now moved \textit{virtually} to the future time level of cell $i$, i.e. $t_i^{n+1}$, and also the virtual cell averages $\tilde{\Q}_j$ are \textit{estimated} from the local predictor solution $\q_{h}$ in the neighbors $T_j$;
	\item \textit{local WENO reconstruction}: once the \textit{virtual} geometry and cell averages have been projected to the future time $t_i^{n+1}$, the local WENO reconstruction technique described 
	in Section \ref{sec.weno} can be carried out for element $T_i$, hence obtaining the new reconstruction polynomial $\w_{h}$ at time $t_i^{n+1}$; 
	\item \textit{local timestep computation}: using the virtual geometry and the virtual solution of the Neumann neighbors, the next local timestep $\Delta t_i^{n+1}$ is evaluated according to \eqref{eqn.timestep};
	\item \textit{local space-time predictor}: finally we compute the high order space-time predictor solution $\q_{h}$ valid within the next timestep of element $T_i$. 
\end{itemize}
This procedure is repeated for all elements, until all of them reach the final time of the simulation $t_f$. As soon as an element $T_i$ has finished its own computation because it has reached 
the final time $t_f$, it is automatically skipped at the beginning of each cycle, waiting for the remaining elements to reach the final time, too.  

This brief description summarizes how our high order Lagrangian LTS algorithm is organized. During the simulation hanging nodes in time appear because each node is moved \textit{physically} only 
by the updating element $T_i$ which the vertex belongs to. As a consequence, the resulting space-time mesh is computed \textit{dynamically}, producing a \textit{non-conforming} space-time mesh. 
Due to our high order approach, the edge and vertex fluxes have to be evaluated using higher order Gaussian quadrature rules, hence increasing the computational cost. In practical applications, 
for which first or second order accurate finite volume schemes are considered adequate, one could rely on the fast and simple mid-point rule that would significantly improve the computational 
efficiency of our LTS algorithm. 
\section{Test problems}
\label{sec.validation} 
\vspace{-2pt}

In the following we solve some numerical test problems in order to validate the high order Lagrangian ADER-WENO algorithm with time accurate local time stepping (LTS) presented so far. We consider the two-dimensional Euler equations of compressible gas dynamics, which can be cast into form \eqref{PDE} with
\begin{equation}
\label{eulerTerms}
\Q = \left( \begin{array}{c} \rho \\ \rho u \\ \rho v \\ \rho E \end{array} \right), \quad \f = \left( \begin{array}{c} \rho u \\ \rho u^2 + p \\ \rho uv \\ u(\rho E + p) \end{array} \right), \quad \g = \left( \begin{array}{c} \rho v \\ \rho uv \\ \rho v^2 + p  \\ v(\rho E + p) \end{array} \right),  
\end{equation}
where the vector of conserved variables is denoted by $\Q$ and the flux tensor is addressed with $\F=(\f,\g)$. Furthermore let $\rho$ and $\rho E$ denote the mass density and the total energy density, respectively, while $\mathbf{v}=(u,v)$ represents the velocity vector and $p$ is the fluid pressure. The source term $\S(\Q)$ is zero for the homogeneous Euler equations. The system is closed using the equation of state (EOS) for an ideal gas, namely
\begin{equation}
\label{eqn.eos} 
p = (\gamma-1)\left(\rho E - \frac{1}{2} \rho (u^2+v^2) \right), 
\end{equation}
where $\gamma$ is the ratio of specific heats.

In the next sections the governing PDE \eqref{PDE}, with the definitions provided by \eqref{eulerTerms}, will be assigned with different initial conditions, that may be given either in terms of the vector of conserved variables $\Q=(\rho,\rho u, \rho v, \rho E)$ or of the primitive variables $\U=(\rho,u,v,p)$. The system will be solved applying the Lagrangian ADER-WENO finite volume schemes illustrated in Section \ref{sec.SolAlg}, choosing among the Rusanov-type \eqref{eqn.rusanov} and the Osher-type \eqref{eqn.osher} numerical fluxes. In all the proposed test problems the local mesh velocity is chosen to be equal to the local fluid velocity ($\mathbf{V}=\mathbf{v}$), hence a formulation of our ALE algorithm has been chosen that comes as close as possible to a truly Lagrangian scheme.

\subsection{Numerical convergence studies}
\label{sec.conv}
In order to carry out the numerical convergence studies for the high order LTS  Lagrangian schemes we consider the classical smooth convected isentropic vortex proposed on triangular grids by Hu and Shu  \cite{HuShuVortex1999}. The initial computational domain is the square $\Omega(0)=[0;10]\times[0;10]$ defined on the $\x=(x,y)$ plane with periodic boundary conditions imposed on each side. The initial condition is given in terms of primitive variables as a linear superposition of a homogeneous background field and a perturbation:
\begin{equation}
\label{ShuVortIC}
\U = (\rho, u, v, p) = (1+\delta \rho, 1+\delta u, 1+\delta v, 1+\delta p).
\end{equation}
The flow is assumed to be isentropic, hence with no perturbation in the entropy,
while the perturbations for velocity $\mathbf{v}=(u,v)$ and temperature $T$ are given by
\begin{equation}
\label{ShuVortDelta}
\left(\begin{array}{c} \delta u \\ \delta v \end{array}\right) = \frac{\epsilon}{2\pi}e^{\frac{1-r^2}{2}} \left(\begin{array}{c} -(y-5) \\ \phantom{-}(x-5) \end{array}\right), \qquad \delta T = -\frac{(\gamma-1)\epsilon^2}{8\gamma\pi^2}e^{1-r^2},
\end{equation} 
where $r^2=(x-5)^2+(y-5)^2$ is the vortex radius, $\epsilon=5$ denotes the vortex strength and the ratio of specific heats is set to $\gamma=1.4$. The perturbations for density and pressure are then expressed as
\begin{equation}
\label{rhopressDelta}
\delta \rho = (1+\delta T)^{\frac{1}{\gamma-1}}-1, \quad \delta p = (1+\delta T)^{\frac{\gamma}{\gamma-1}}-1. 
\end{equation}

\begin{figure}[!htbp]
\begin{center}
\begin{tabular}{ccc} 
\includegraphics[width=0.33\textwidth]{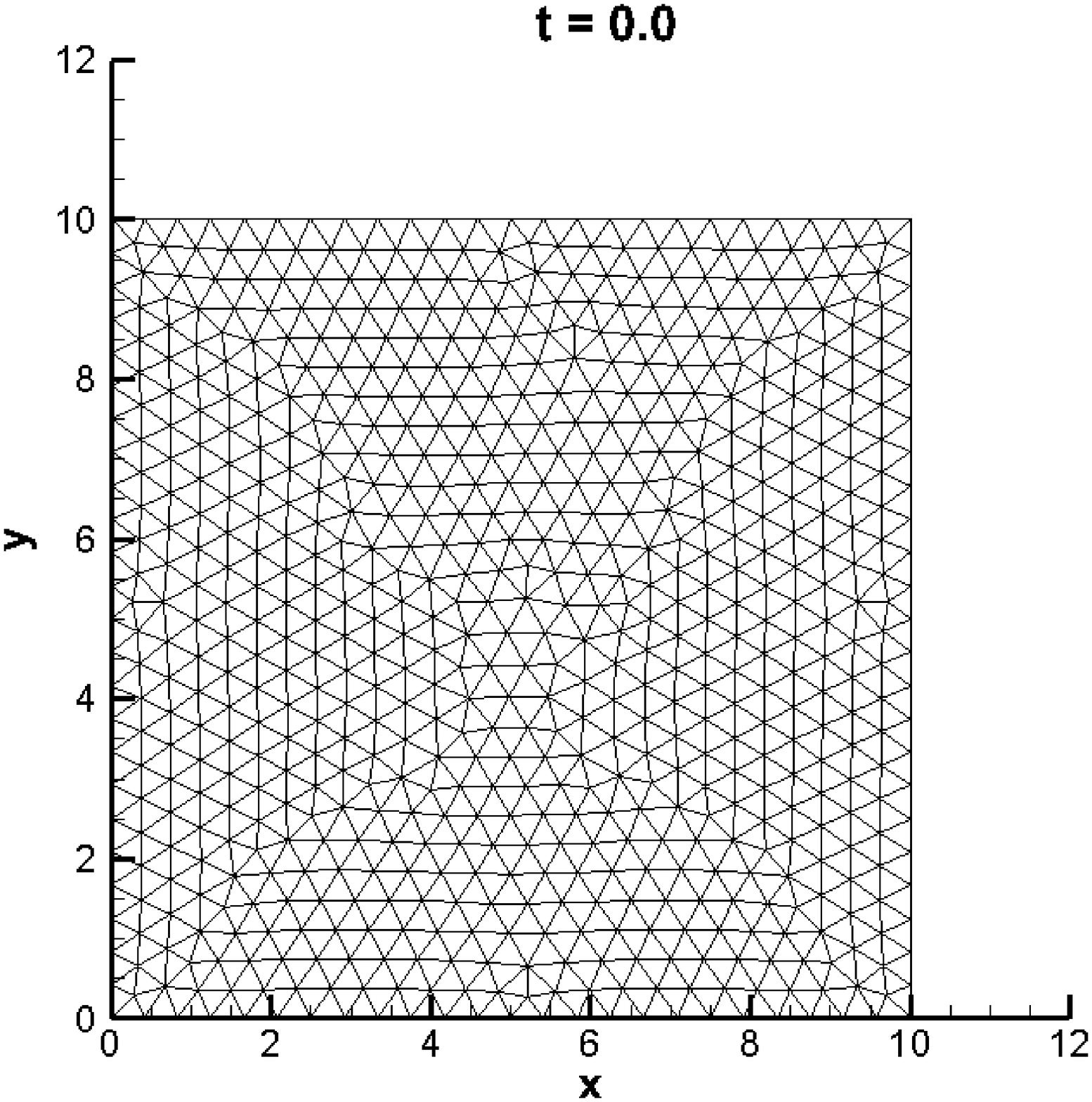}  &           
\includegraphics[width=0.33\textwidth]{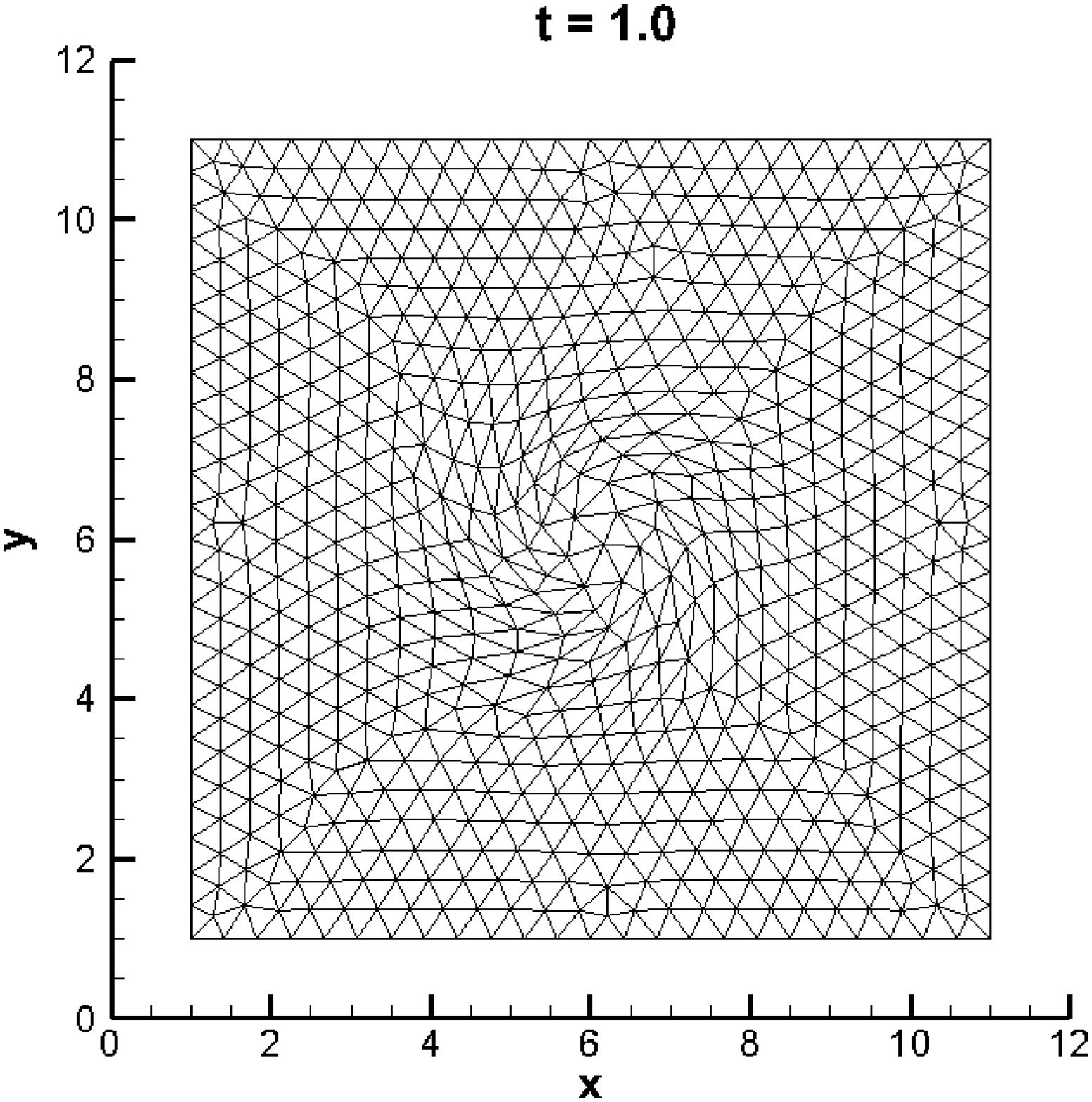}    &
\includegraphics[width=0.33\textwidth]{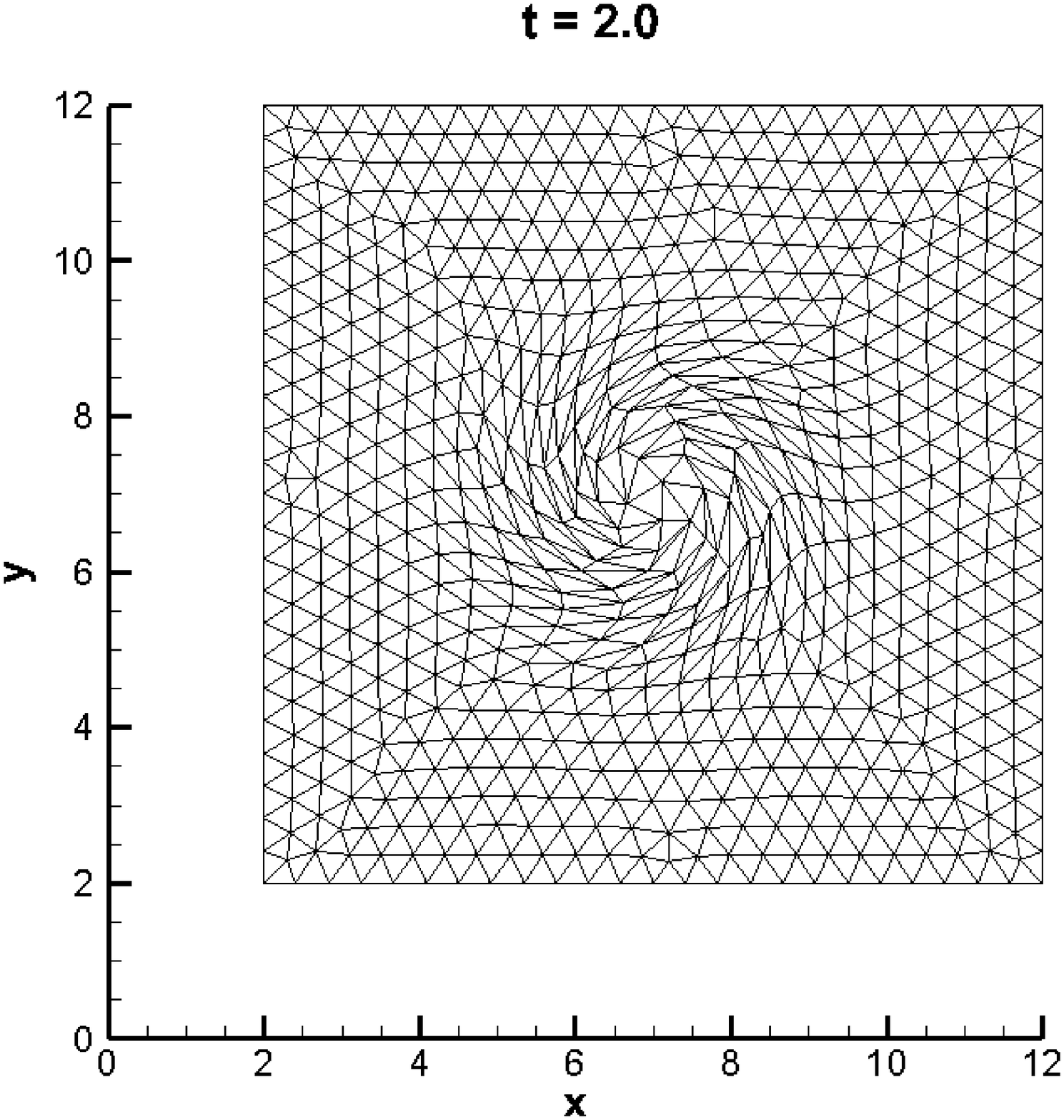}      \\
\end{tabular} 
\caption{Mesh configuration at three different output times for the smooth isentropic vortex test problem. The mesh is highly twisted in the center of the computational domain, which is furthermore convected with velocity $\v_c=(1,1)$.} 
\label{fig.SVgrid}
\end{center}
\end{figure}

The vortex is convected with velocity $\v_c=(1,1)$, so that at the final time $t_f$ of the simulation the exact solution $\Q_e(\x,t_f)$ is simply given by the time-shifted initial condition, e.g. $\Q_e(\x,t_f)=\Q(\x-\v_c t_f,0)$, with the averaged convection velocity of the vortex $\v_c=(1,1)$. As depicted in Figure \ref{fig.SVgrid}, the mesh is highly distorted and twisted by the vortex motion and no rezoning algorithm \cite{LagrangeMHD,LagrangeMDRS} is adopted here because we want to validate the new LTS algorithm inside an almost fully Lagrangian approach. Therefore the final time of the 
simulation is chosen to be $t_f=1.0$, which allows the computational mesh to remain reasonably well-shaped. We run this test case on successive refined meshes and for each mesh the corresponding 
error is expressed in the continuous $L_2$ norm as  
\begin{equation}
  \epsilon_{L_2} = \sqrt{ \int \limits_{\Omega(t_f)} \left( \Q_e(x,y,t_f) - \w_h(x,y,t_f) \right)^2 dxdy },  
	\label{eqnL2error}
\end{equation}
where $\w_h(x,y,t_f)$ represents the high order reconstructed solution at the final time, while the mesh size $h(\Omega(t_f))$ is evaluated as the maximum diameter of the 
circumcircles of the triangles in the final computational domain $\Omega(t_f)$. We use the Rusanov-type numerical flux \eqref{eqn.rusanov} to obtain the convergence results 
listed in Table \ref{tab.convEul}, achieving the designed order of accuracy of the scheme very well.  

\begin{table}  
\caption{Numerical convergence results for the compressible Euler equations using second to fourth order Lagrangian ADER-WENO finite volume schemes with time accurate local time stepping (LTS). The error norms refer to the variable $\rho$ (density) at time $t=1.0$.}  
\begin{center} 
\begin{small}
\renewcommand{\arraystretch}{1.0}
\begin{tabular}{ccccccccc} 
\hline
  \multicolumn{3}{c}{$\mathcal{O}2$} & \multicolumn{3}{c}{$\mathcal{O}3$}  & \multicolumn{3}{c}{$\mathcal{O}4$}  \\
\hline
  $h(\Omega(t_f))$ & $\epsilon_{L_2}$ & $\mathcal{O}(L_2)$ & $h(\Omega,t_f)$ & $\epsilon_{L_2}$ & $\mathcal{O}(L_2)$ & $h(\Omega(t_f))$ & $\epsilon_{L_2}$ & $\mathcal{O}(L_2)$ \\ 
\hline
3.58E-01 & 5.286E-02 & -   & 3.32E-01  & 3.471E-02 & -   & 7.00E-01 & 6.419E-02 & -    \\ 
2.48E-01 & 3.558E-02 & 1.1 & 2.51E-01  & 1.789E-02 & 2.4 & 3.28E-01 & 1.030E-02 & 2.4  \\ 
1.70E-01 & 1.514E-02 & 2.3 & 1.68E-01  & 6.346E-03 & 2.6 & 2.51E-01 & 3.598E-03 & 3.9  \\ 
1.28E-01 & 8.193E-03 & 2.1 & 1.28E-01  & 2.935E-03 & 2.8 & 1.68E-01 & 7.706E-04 & 3.8  \\ 
\hline 
\end{tabular}
\end{small}
\end{center}
\label{tab.convEul}
\end{table}

\subsection{Riemann problems} 
\label{sec.RP2D}
Here we solve two classical Riemann problems, namely the shock tube problems of Sod and the Lax, which are in the following addressed as RP1 and RP2, respectively, 
and which are widely adopted to validate numerical algorithms for the solution of the compressible Euler equations. 
They both include the formation of a left-propagating rarefaction wave, an intermediate contact discontinuity and a right-propagating shock wave. 
Though intrinsically one-dimensional, these tests become non-trivial and multidimensional when applied to unstructured meshes, where in general the element 
edges are not aligned with the fluid motion. Since a contact wave is present in the solution, we can also check how well it is resolved by our Lagrangian LTS 
scheme. 

The initial computational domain is given by the box $\Omega(0)=[-0.5;0.5]\times[-0.05;0.05]$ that is discretized with a characteristic mesh size of $h=1/200$, 
leading to a total number of $N_E=8862$ elements, while the initial conditions are given in terms of the primitive variables $\U=(\rho,u,v,p)$. Table~\ref{tab.iniRP2D}  
reports the relevant data for the setup of the two tests, where $t_f$ represents the final time of the simulation while $x_d$ gives the position of the initial discontinuity 
which splits the computational domain, as well as the initial conditions, in the two left and right states  $\U_L$ and $\U_R$. 
We set periodic boundary conditions in the $y$ direction, while transmissive boundaries are imposed along the $x$ direction. The ratio of specific heats is assumed to be 
$\gamma=1.4$ for both Riemann problems.

\begin{table}  
\caption{Initial condition for the Sod (RP1) and the Lax (RP2) shock tube problem. $t_f$ is the final time of the simulation and $x_d$ denotes the position of the initial discontinuity.}  
\begin{center} 
\begin{small}
\renewcommand{\arraystretch}{1.0}
\begin{tabular}{ccccccccccc} 
\hline
 Case   &  $\rho_L$ &  $u_L$   &  $v_L$  &  $p_L$  &  $\rho_R$  &  $u_R$   &  $v_R$  &  $p_R$ &  $t_f$  &  $x_d$  \\
\hline		
RP1 &  1.0     &  0.0   &  0.0  &  1.0   &  0.125   &  0.0   &  0.0  &  0.1   &  0.2  &  0.0  \\
RP2 &  0.445   &  0.698 &  0.0  &  3.528 &  0.5     &  0.0   &  0.0  &  0.571 &  0.1  &  0.0  \\
\hline   
\end{tabular}
\end{small}
\end{center}
\label{tab.iniRP2D}
\end{table}

\begin{figure}[!htbp]
\begin{center}
\begin{tabular}{cc} 
\includegraphics[width=0.47\textwidth]{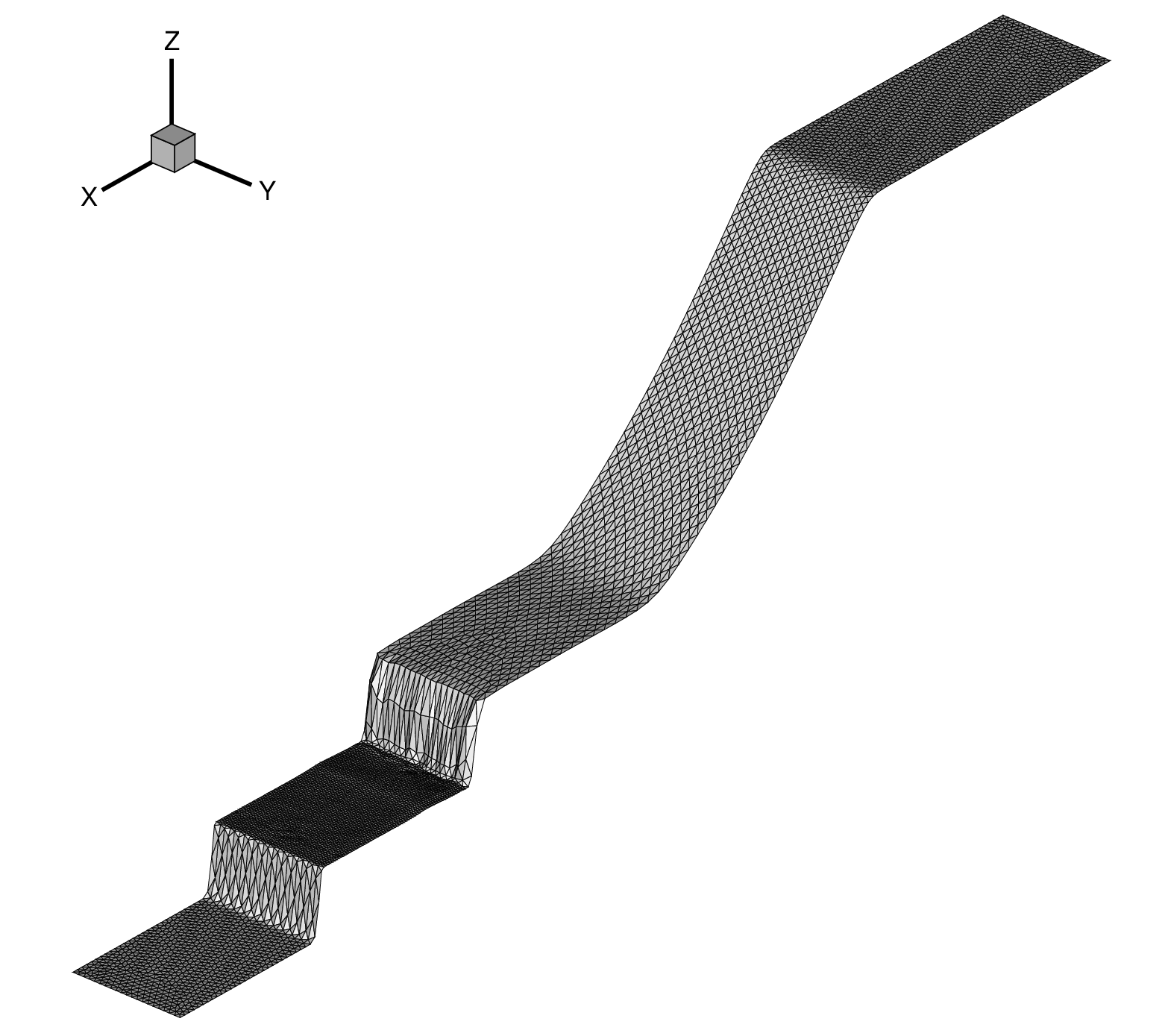}  &           
\includegraphics[width=0.47\textwidth]{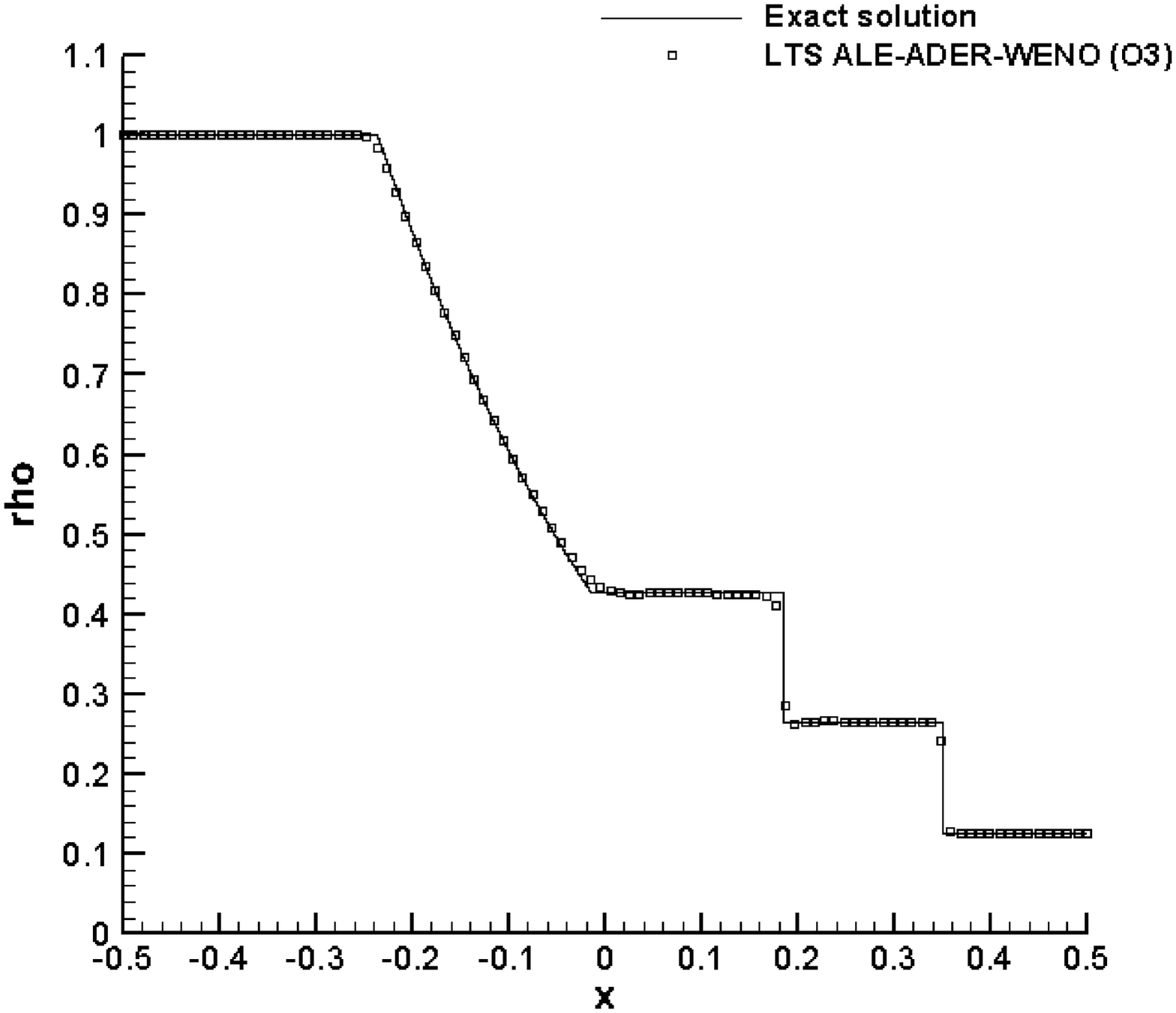} \\
\includegraphics[width=0.47\textwidth]{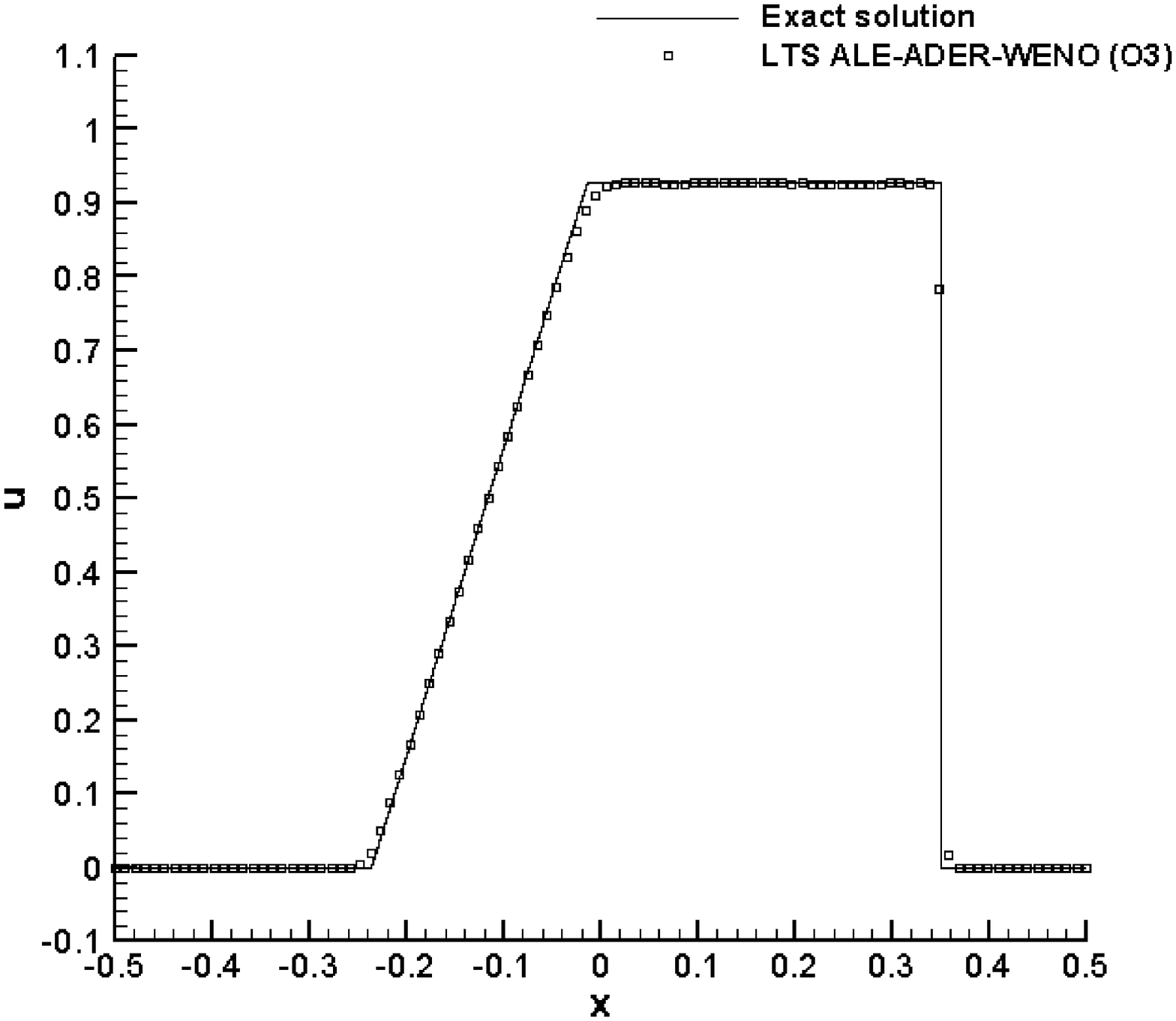}  &           
\includegraphics[width=0.47\textwidth]{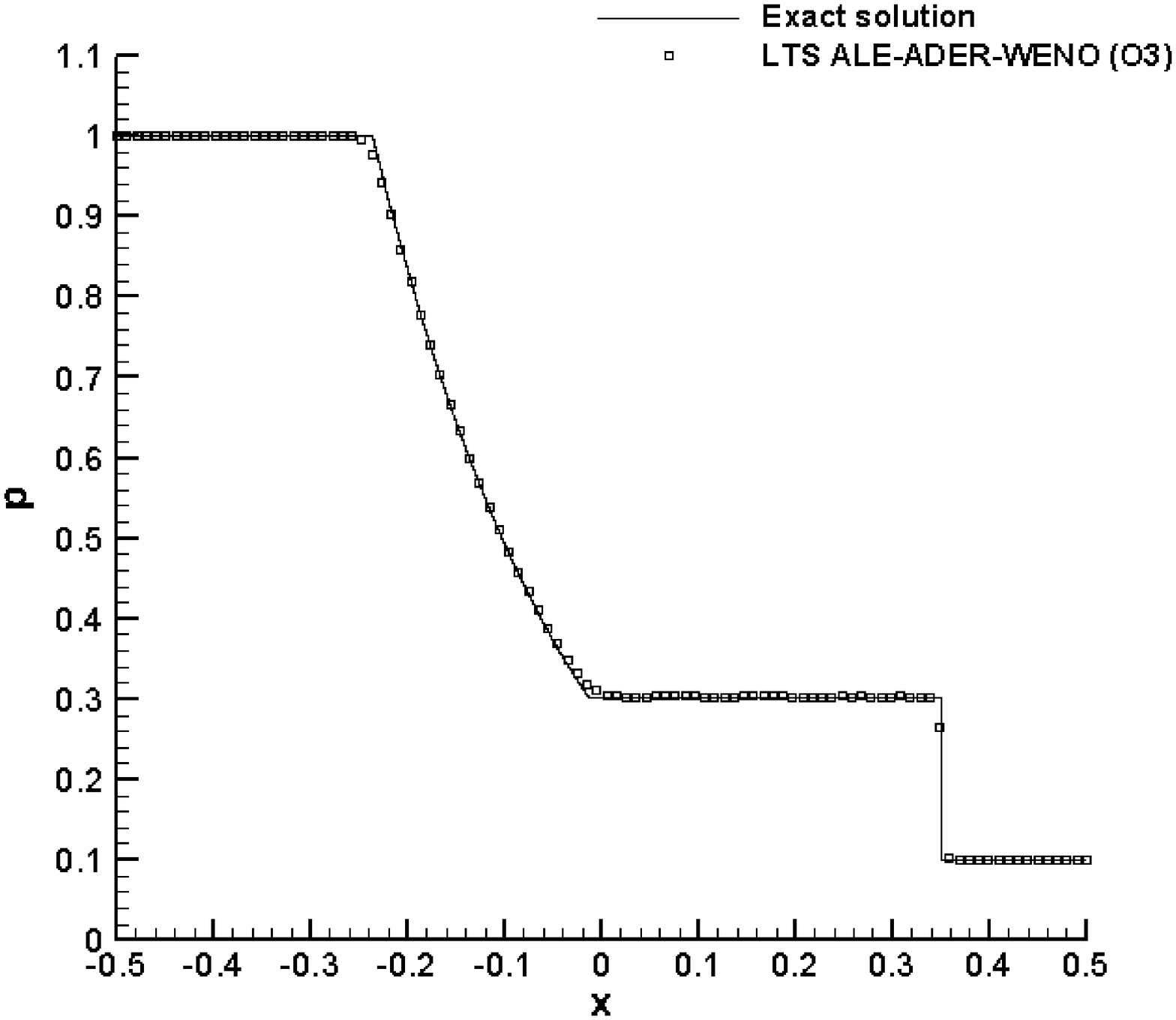} \\    
\end{tabular} 
\caption{Comparison between exact and third order accurate numerical solution for the Sod shock tube problem RP1. Density (top right), velocity (bottom left) and pressure (bottom right) distribution are shown as well as a 3D view of the density solution at the final time $t_f=0.2$ (top left).} 
\label{fig.RP1}
\end{center}
\end{figure}

\begin{figure}[!htbp]
\begin{center}
\begin{tabular}{cc} 
\includegraphics[width=0.47\textwidth]{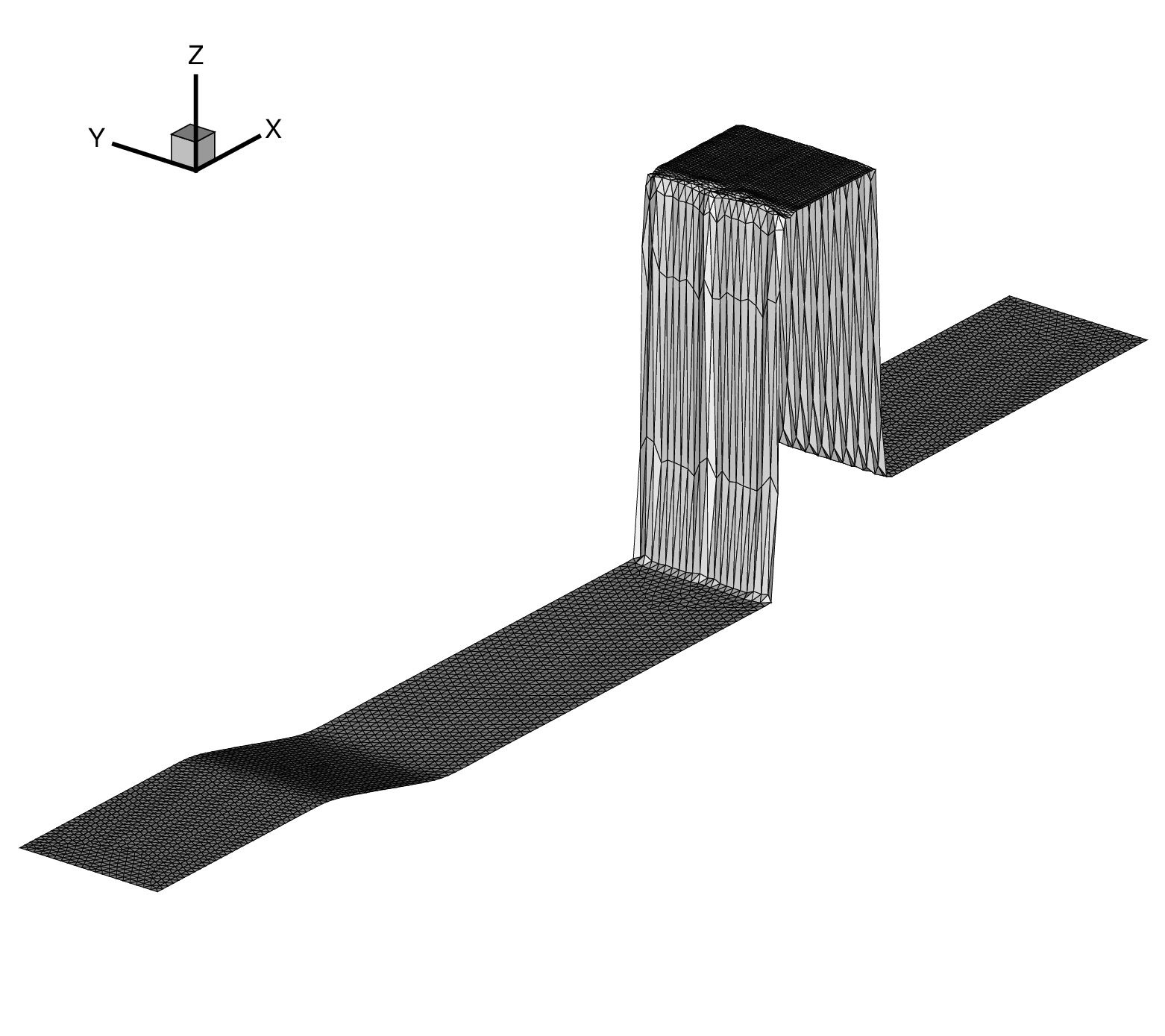}  &           
\includegraphics[width=0.47\textwidth]{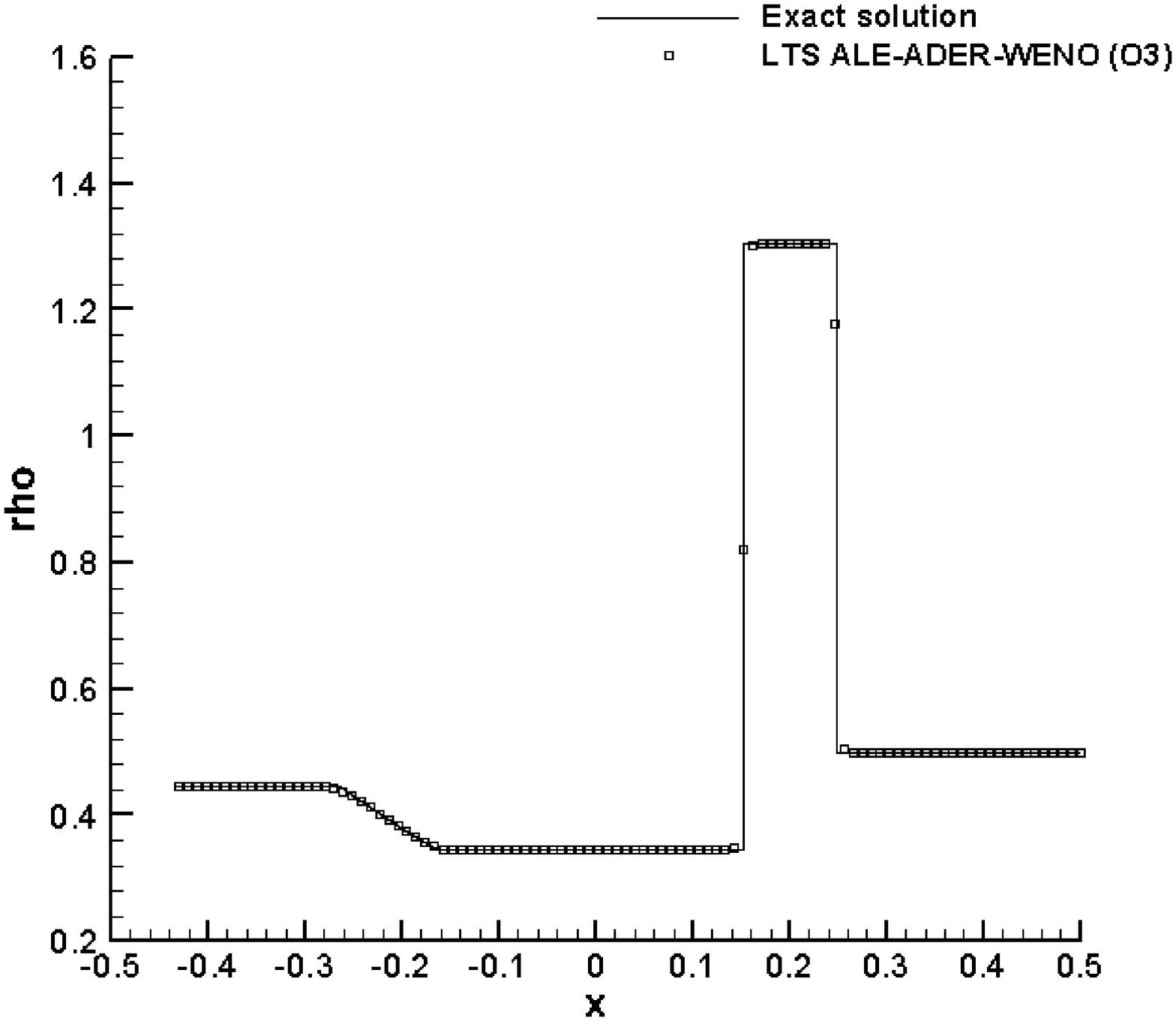} \\
\includegraphics[width=0.47\textwidth]{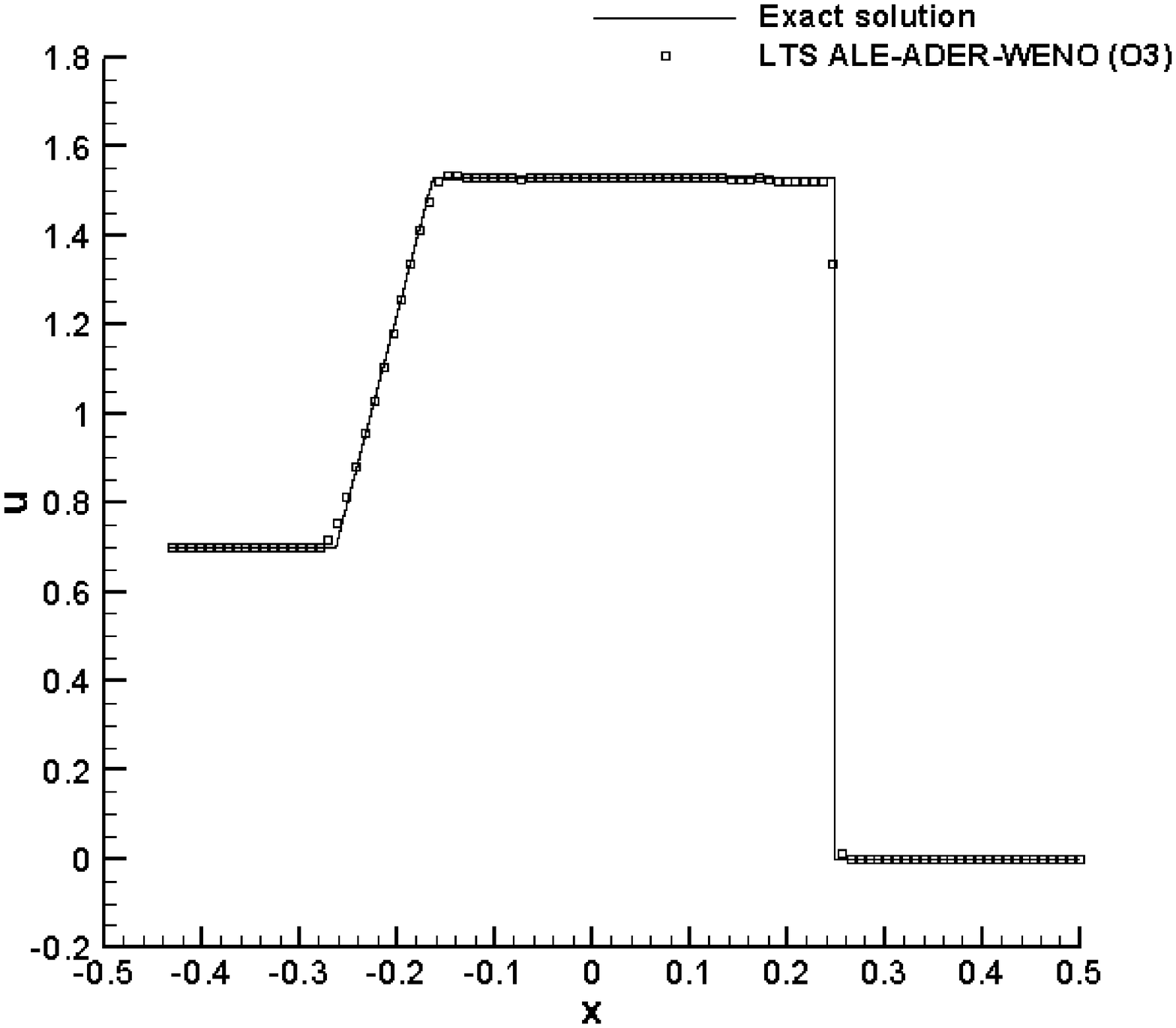}  &           
\includegraphics[width=0.47\textwidth]{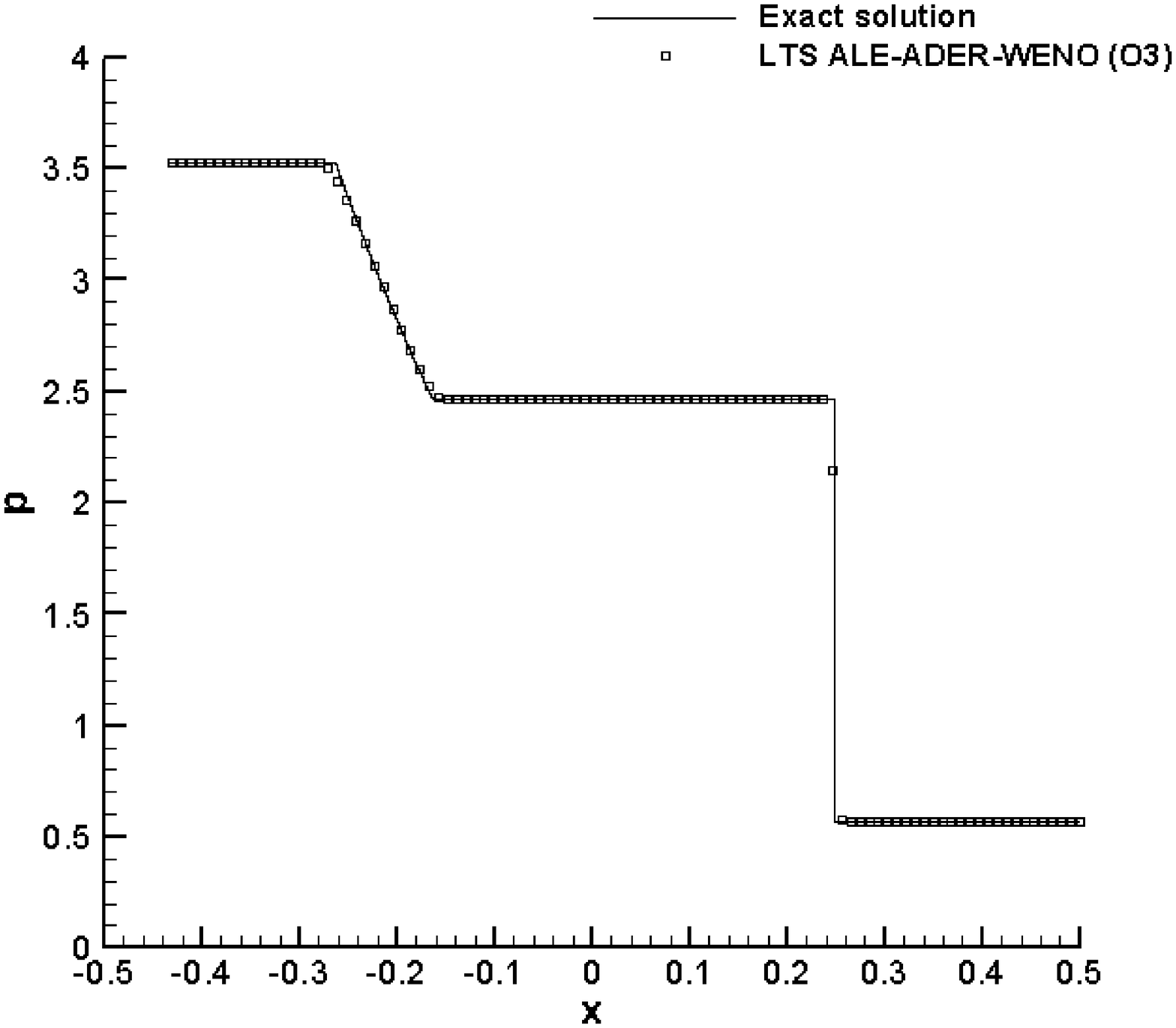} \\    
\end{tabular} 
\caption{Comparison between exact and third order accurate numerical solution for the Lax shock tube problem RP2. Density (top right), velocity (bottom left) and pressure (bottom right) distribution are shown as well as a 3D view of the density solution at the final time $t_f=0.1$ (top left).} 
\label{fig.RP2}
\end{center}
\end{figure}

The exact solution is computed with the exact Riemann solver presented in \cite{ToroBook}. We use the third order version of our Lagrangian ADER-WENO schemes with LTS using the Osher-type 
numerical flux to obtain the results depicted in Figures \ref{fig.RP1}-\ref{fig.RP2}, where a comparison between the exact and the numerical solution is shown. We observe an excellent 
resolution of the contact wave with only one intermediate point for both RP1 and RP2, and a very good agreement with the analytical solution can also be noticed for density, as well as 
for pressure and for the horizontal velocity component. Table \ref{tab.RPdata} aims at showing the computational efficiency of the LTS algorithm w.r.t. the Lagrangian ADER-WENO schemes 
with global time stepping (GTS) presented in \cite{Lagrange2D}. In order to give a fair comparison between LTS and GTS schemes, the efficiency is not measured in terms of computational 
time, which may depend on the machine hardware or on the algorithm implementation, but rather we count the total number of element updates needed to reach the final time of the simulation, 
as done in \cite{ALELTS1D}. Hence, looking at Table \ref{tab.RPdata}, we notice that the Lagrangian algorithm with global time stepping requires a total number of element updates that is 
a factor of 3-4 times larger than the one of our new Lagrangian scheme with LTS. 

%\begin{figure}[!htbp]
%\begin{center}
%\begin{tabular}{cc} 
%\includegraphics[width=0.47\textwidth]{./RP3-LTS-rho3D.eps}  &           
%\includegraphics[width=0.47\textwidth]{./RP3-LTS-rho.eps} \\
%\includegraphics[width=0.47\textwidth]{./RP3-LTS-u.eps}  &           
%\includegraphics[width=0.47\textwidth]{./RP3-LTS-p.eps} \\    
%\end{tabular} 
%\caption{Comparison between exact and third order accurate numerical solution for the Riemann problem RP3. Density (top right), velocity (bottom left) and pressure (bottom right) distribution are shown as well as a 3D view of the density solution at the final time $t_f=0.012$.} 
%\label{fig.RP3}
%\end{center}
%\end{figure}

\begin{table}[!htbp]
  \caption{Comparison of the computational efficiency between GTS and LTS algorithm in terms of  the total number of element updates for RP1 and RP2. A third order scheme has been adopted.}
	\begin{center}
		\begin{tabular}{cccc}
		\hline
		         & \multicolumn{3}{c}{Number of element updates} 		\\
		Case     &      GTS                &     LTS                 &  GTS/LTS \\
		\hline
		RP1      & $10.120404 \cdot 10^6 $ & $3.257847 \cdot 10^6 $ & 3.11 \\
		RP2      & $23.349964 \cdot 10^6 $ & $5.020780 \cdot 10^6 $ & 4.45 \\
		\hline
		\end{tabular}
	\end{center}
	\label{tab.RPdata}
\end{table}

\subsection{Two-dimensional explosion problems} 
\label{sec.EP2D}
Circular explosion problems can be regarded as the two-dimensional extension of Riemann problems. The initial domain $\Omega(0)=\left\{ \mathbf{x} : \left\| \mathbf{x} \right\| < R_o \right\} $ 
is given  by the unit circle of radius $R_o=1$. A circle of radius $R=0.5$ separates two different states that define the initial conditions reported in Table \ref{tab.iniEP2D} in terms of 
primitive variables $\U=(\rho,u,v,p)$. The two states are addressed here as the \textit{inner} state $\U_i$ and the \textit{outer} state $\U_o$, respectively. Transmissive boundary conditions have 
been imposed on the external boundary and we set $\gamma=1.4$. EP1 corresponds to the initial data of the classical Sod shock tube problem RP1, while EP2 is taken from \cite{ToroBook}. 
In both cases we use the same computational mesh $m_1$, with a characteristic mesh size of $h=1/100$ for $r\leq R$ and $h=1/50$ for $r>R$, hence obtaining a total number $N_E=43756$ of triangles.   

\begin{table}  
\caption{Initial conditions for the two-dimensional explosion problems EP1 and the EP2 with $t_f$ denoting the final time of the simulation.}  
\begin{center} 
\begin{small}
\renewcommand{\arraystretch}{1.0}
\begin{tabular}{cccccccccc} 
\hline
 Case   &  $\rho_i$ &  $u_i$   &  $v_i$  &  $p_i$  &  $\rho_o$  &  $u_o$   &  $v_o$  &  $p_o$ &  $t_f$    \\
\hline		
EP1 &  1.0     &  0.0   &  0.0  &  1.0     &  0.125   &  0.0   &  0.0  &  0.1  &  0.2    \\
EP2 &  1.0     &  0.0   &  0.0  &  1000.0  &  1.0     &  0.0   &  0.0  &  0.01 &  0.012    \\
\hline   
\end{tabular}
\end{small}
\end{center}
\label{tab.iniEP2D}
\end{table}

\begin{figure}[!htbp]
\begin{center}
\begin{tabular}{cc} 
\includegraphics[width=0.47\textwidth]{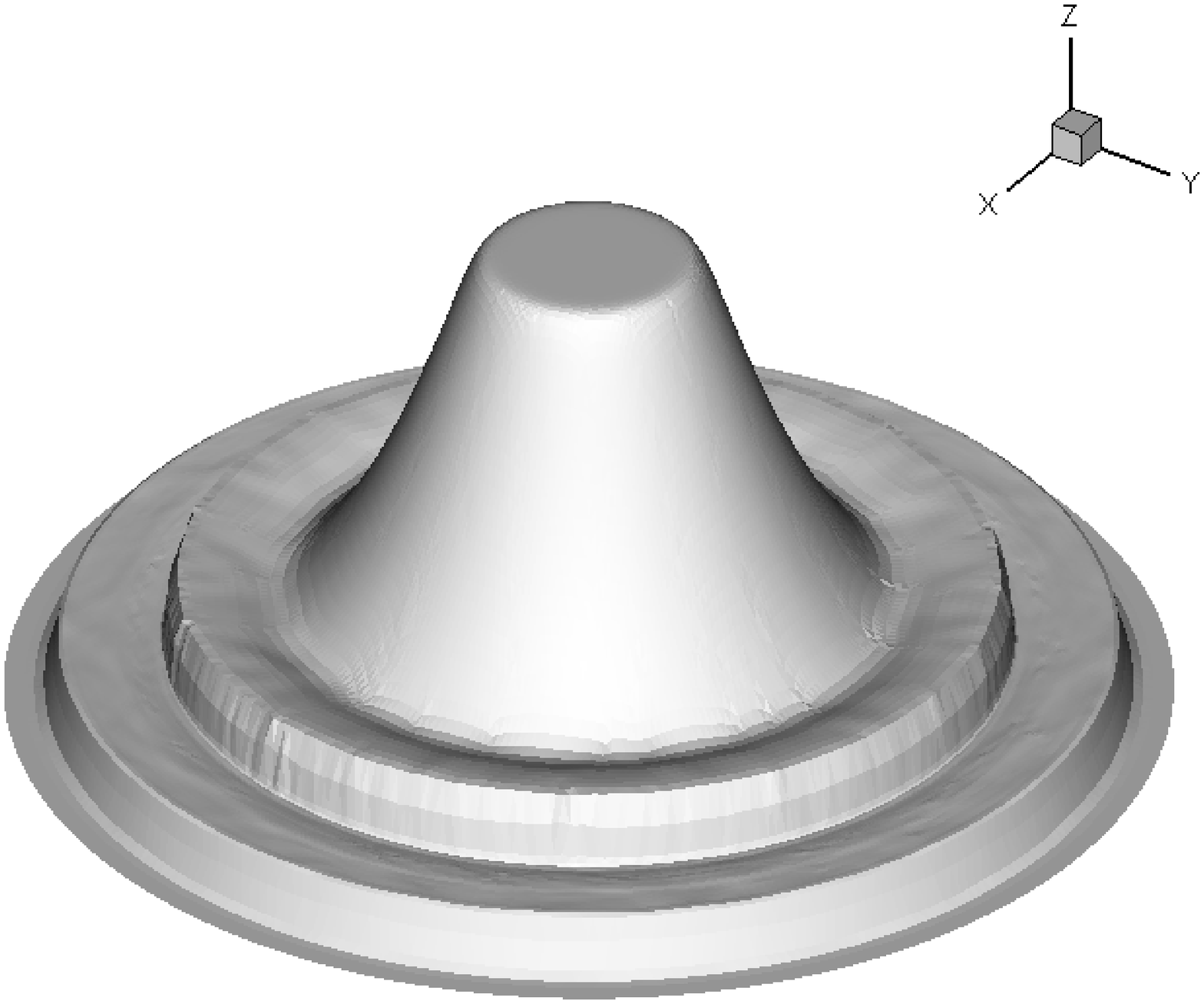}  &           
\includegraphics[width=0.47\textwidth]{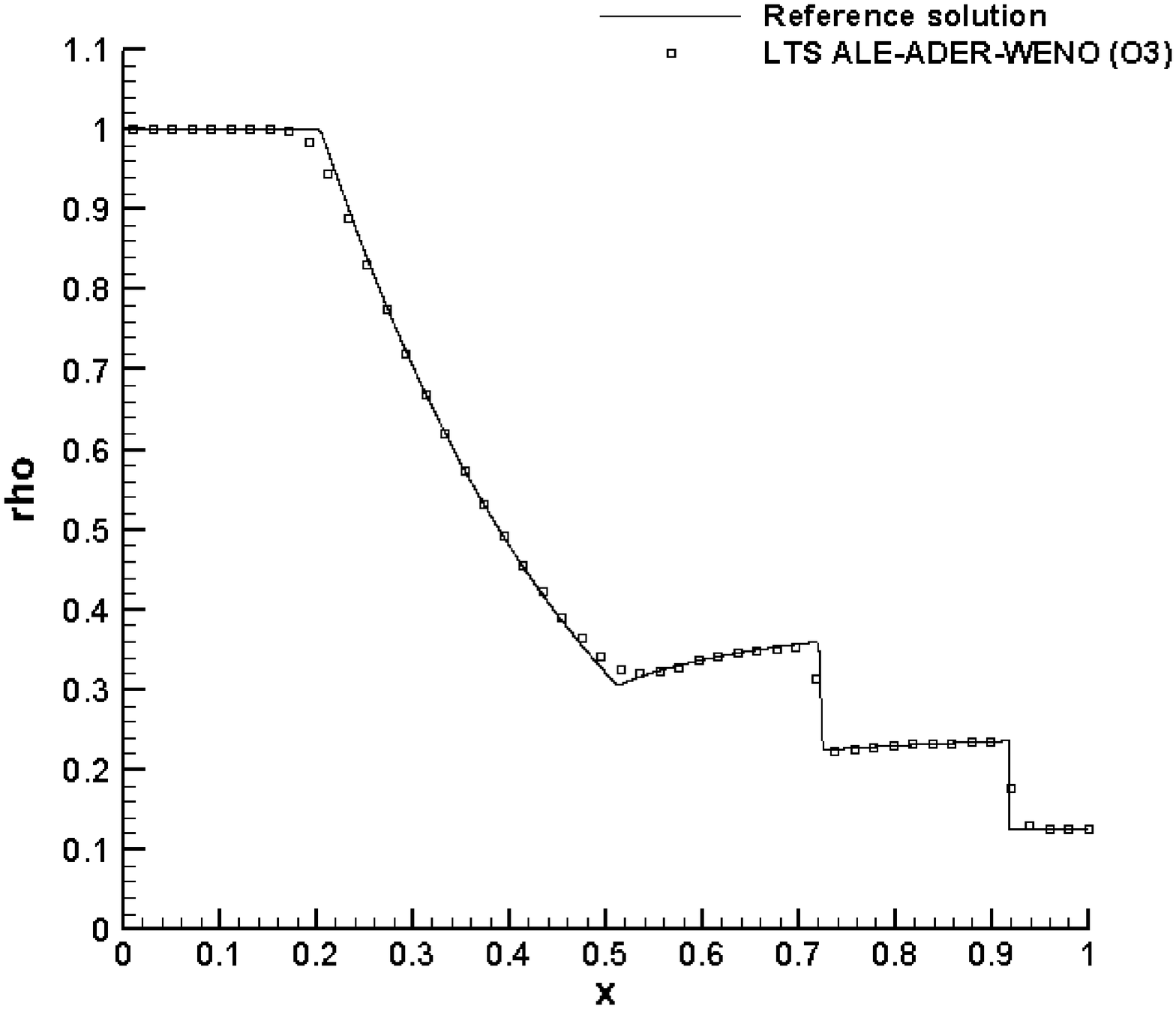} \\
\includegraphics[width=0.47\textwidth]{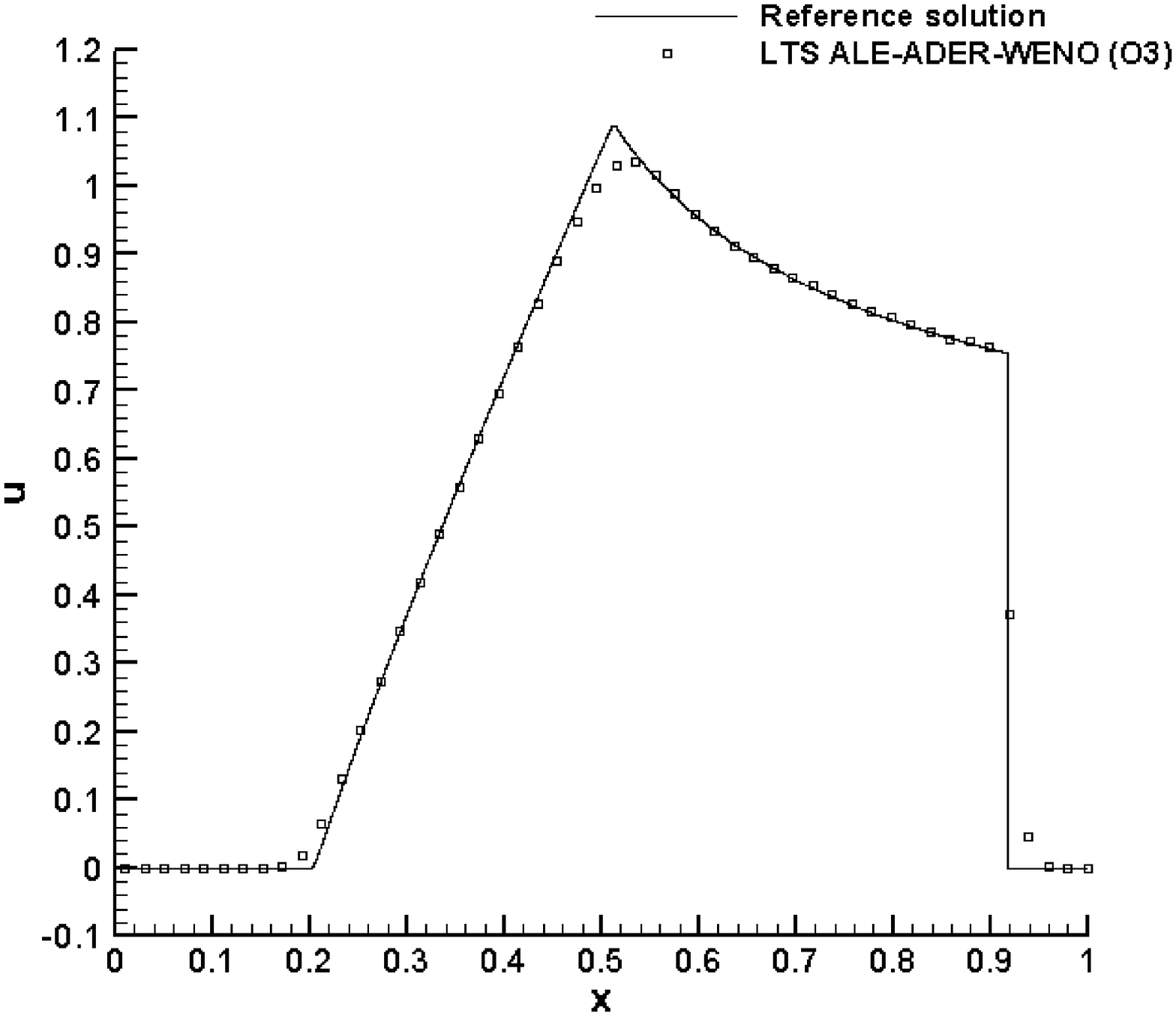}  &           
\includegraphics[width=0.47\textwidth]{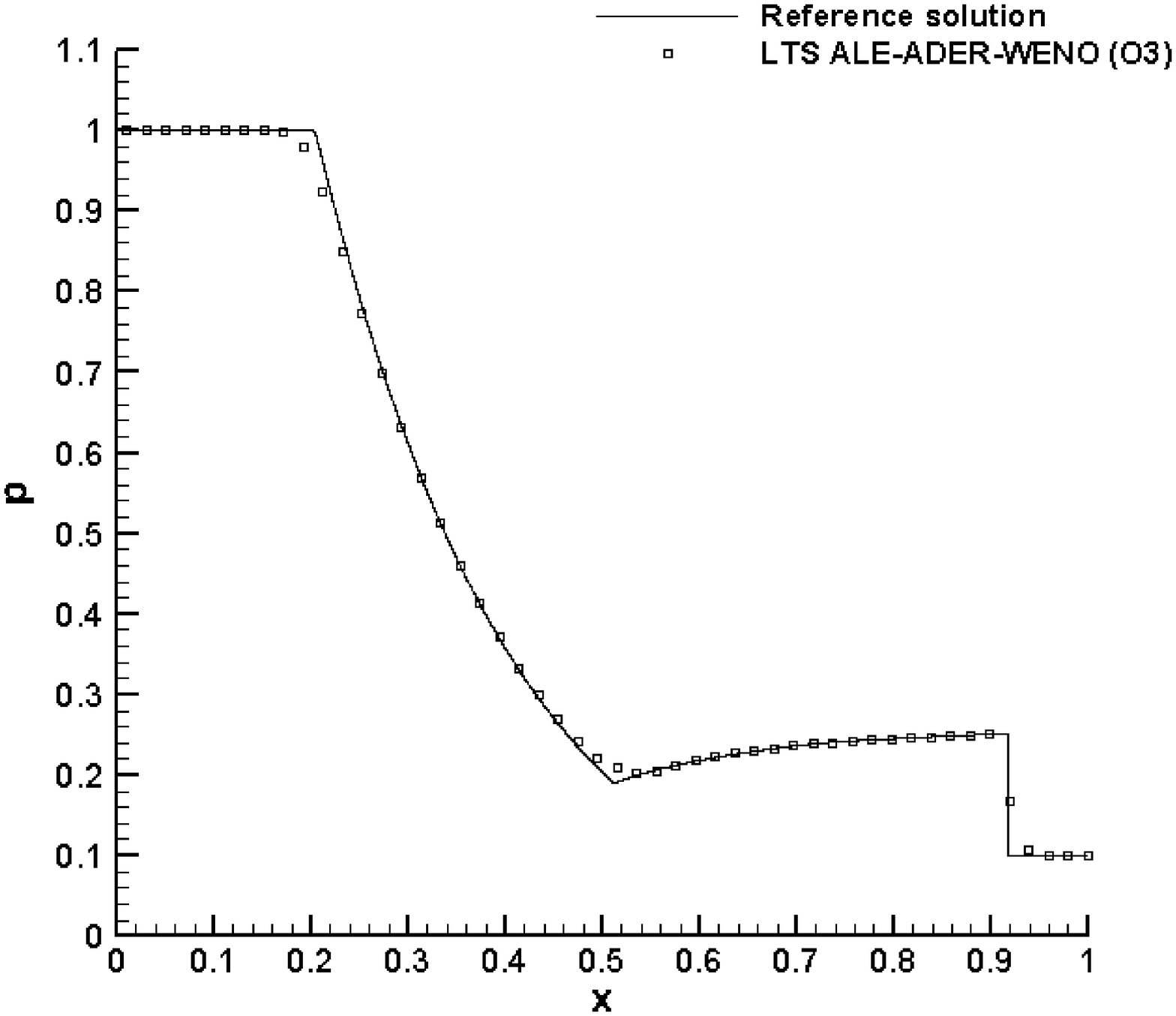} \\    
\end{tabular} 
\caption{Comparison between reference and third order accurate numerical solution for the explosion problem EP1. Density (top right), velocity (bottom left) and pressure (bottom right) distribution are shown as well as a 3D view of the density solution at the final time $t_f=0.25$ (top left).} 
\label{fig.EP1}
\end{center}
\end{figure}

\begin{figure}[!htbp]
\begin{center}
\begin{tabular}{cc} 
\includegraphics[width=0.47\textwidth]{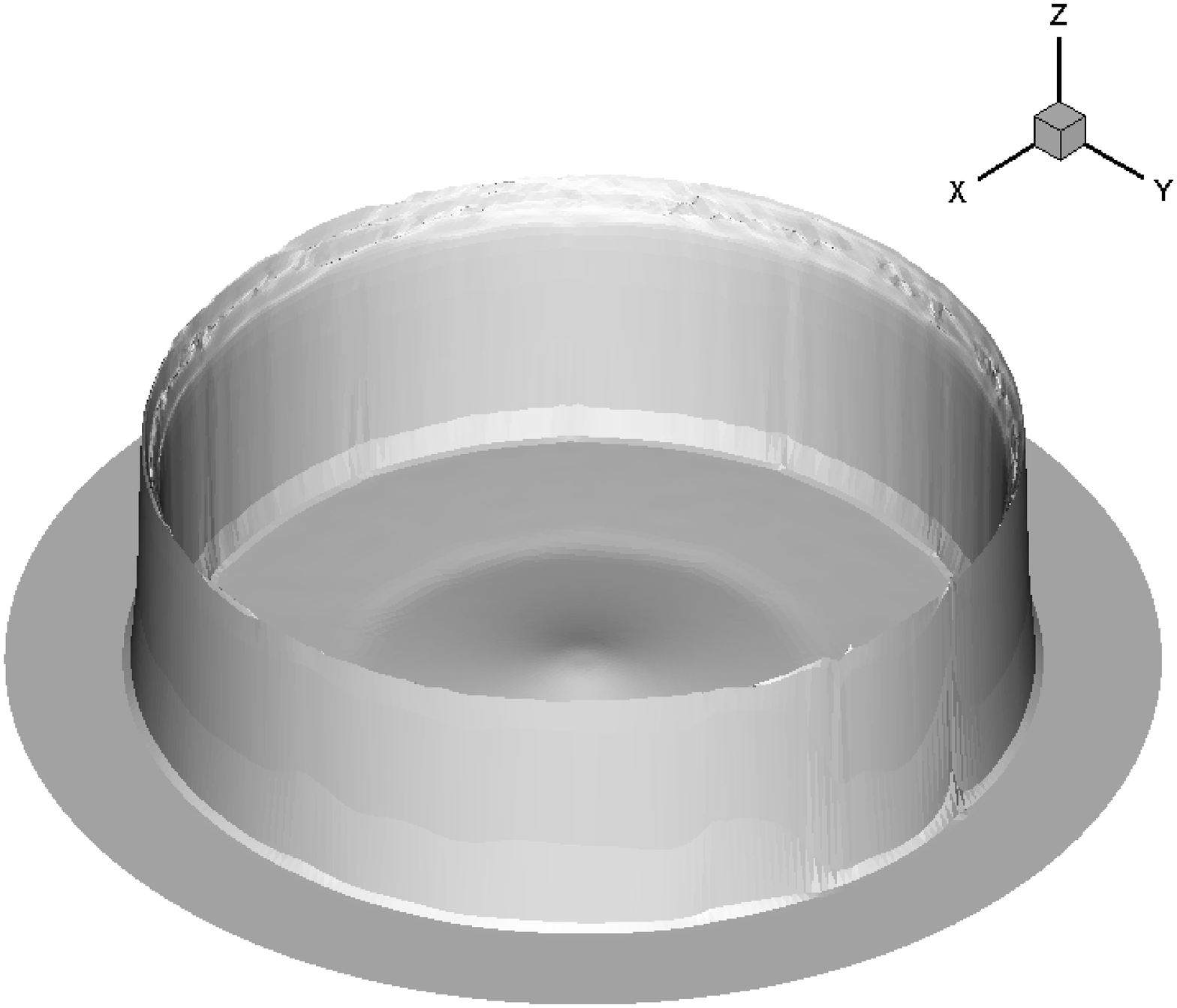}  &           
\includegraphics[width=0.47\textwidth]{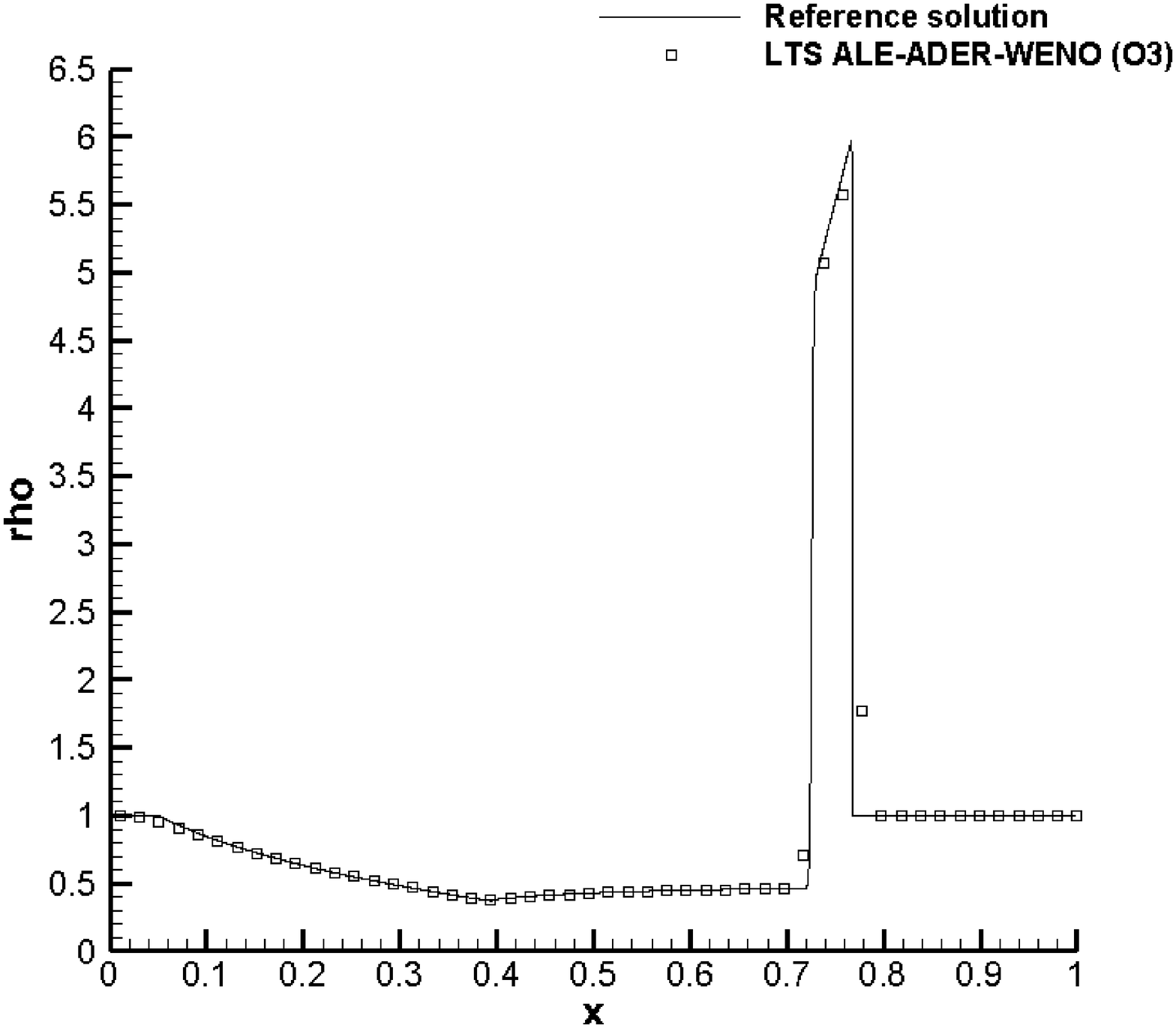} \\
\includegraphics[width=0.47\textwidth]{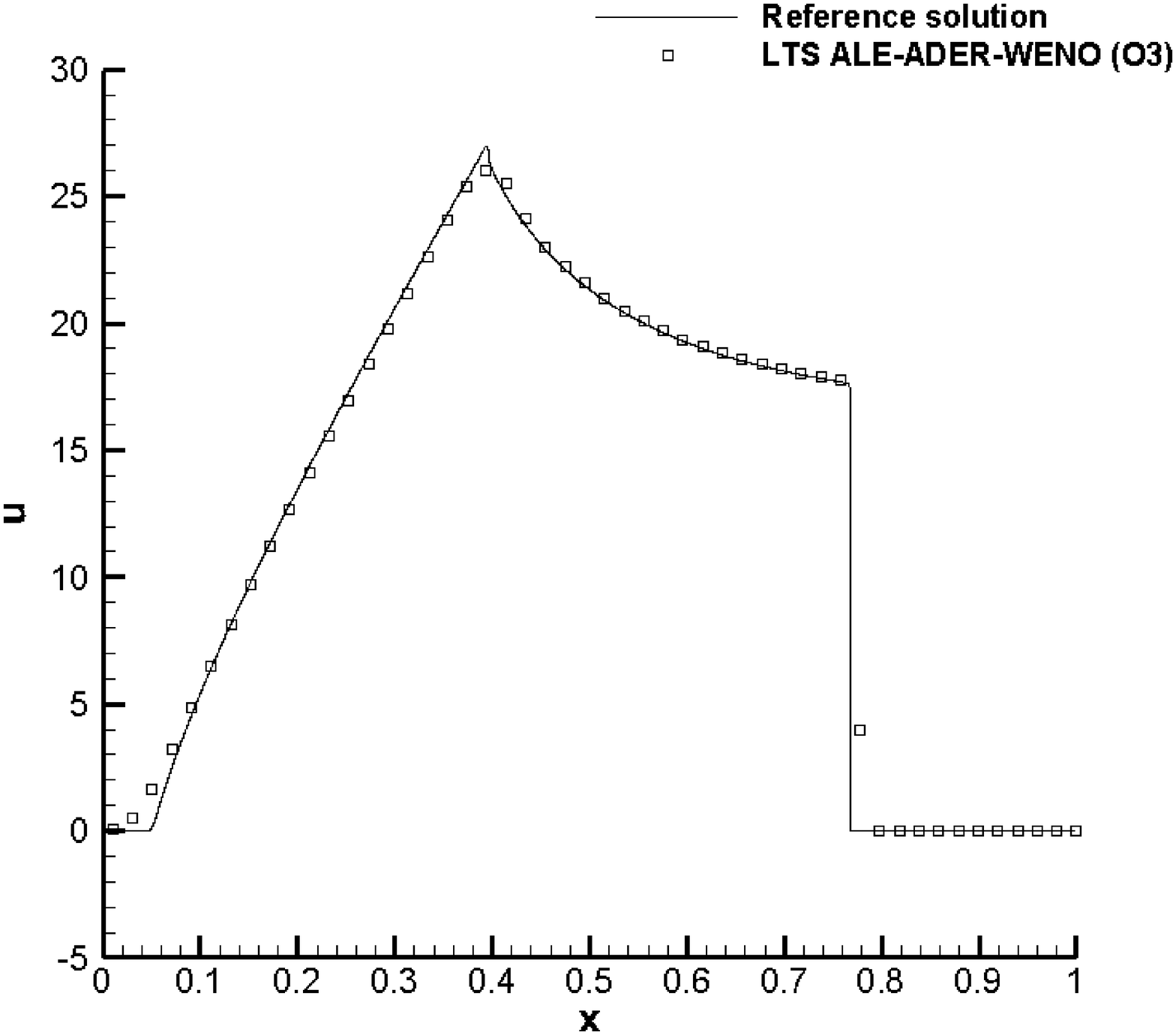}  &           
\includegraphics[width=0.47\textwidth]{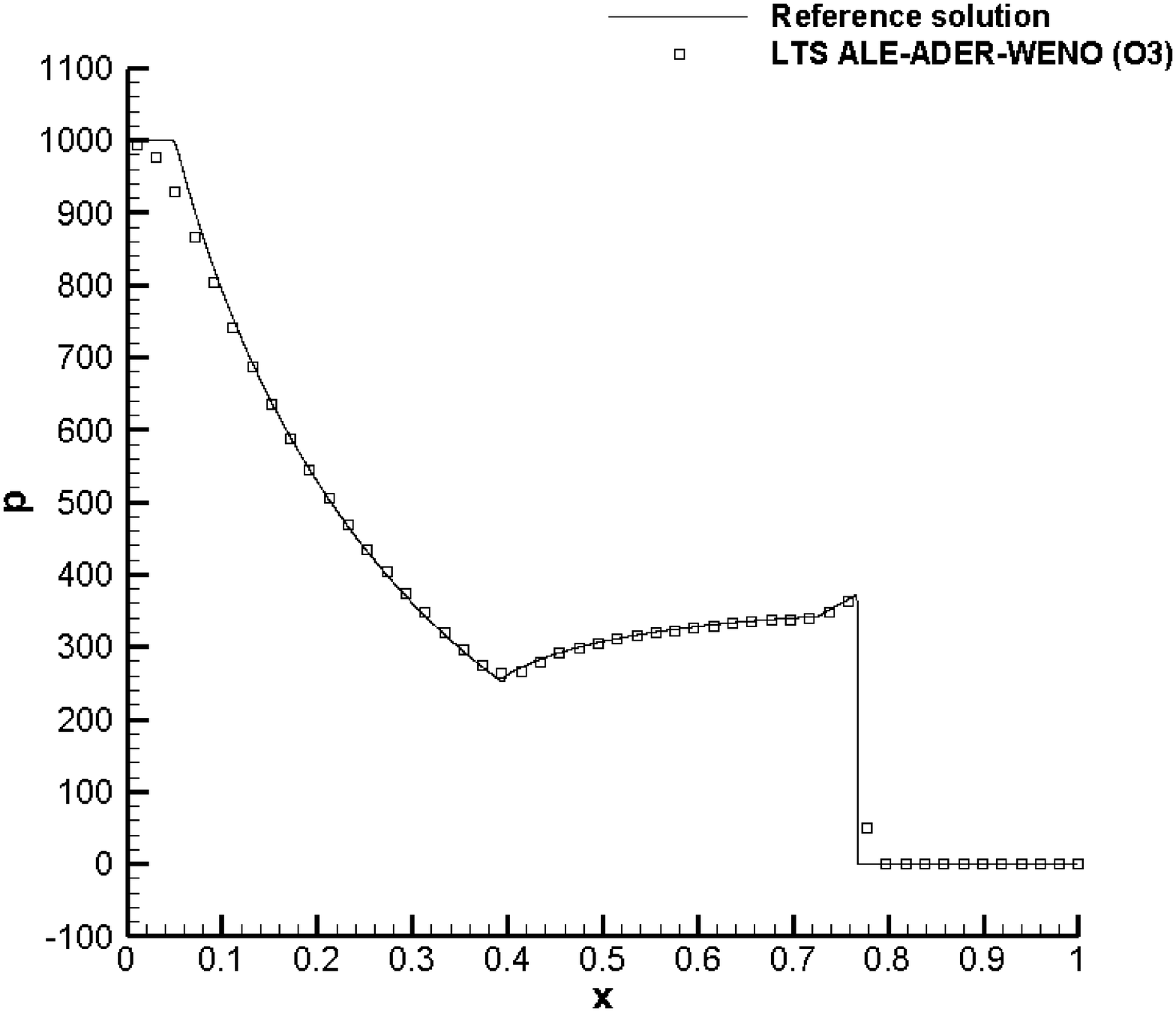} \\    
\end{tabular} 
\caption{Comparison between reference and third order accurate numerical solution for the explosion problem EP2. Density (top right), velocity (bottom left) and pressure (bottom right) distribution are shown as well as a 3D view of the density solution at the final time $t_f=0.012$ (top left).} 
\label{fig.EP3}
\end{center}
\end{figure}

\begin{figure}[!htbp]
\begin{center}
\begin{tabular}{cc} 
\includegraphics[width=0.47\textwidth]{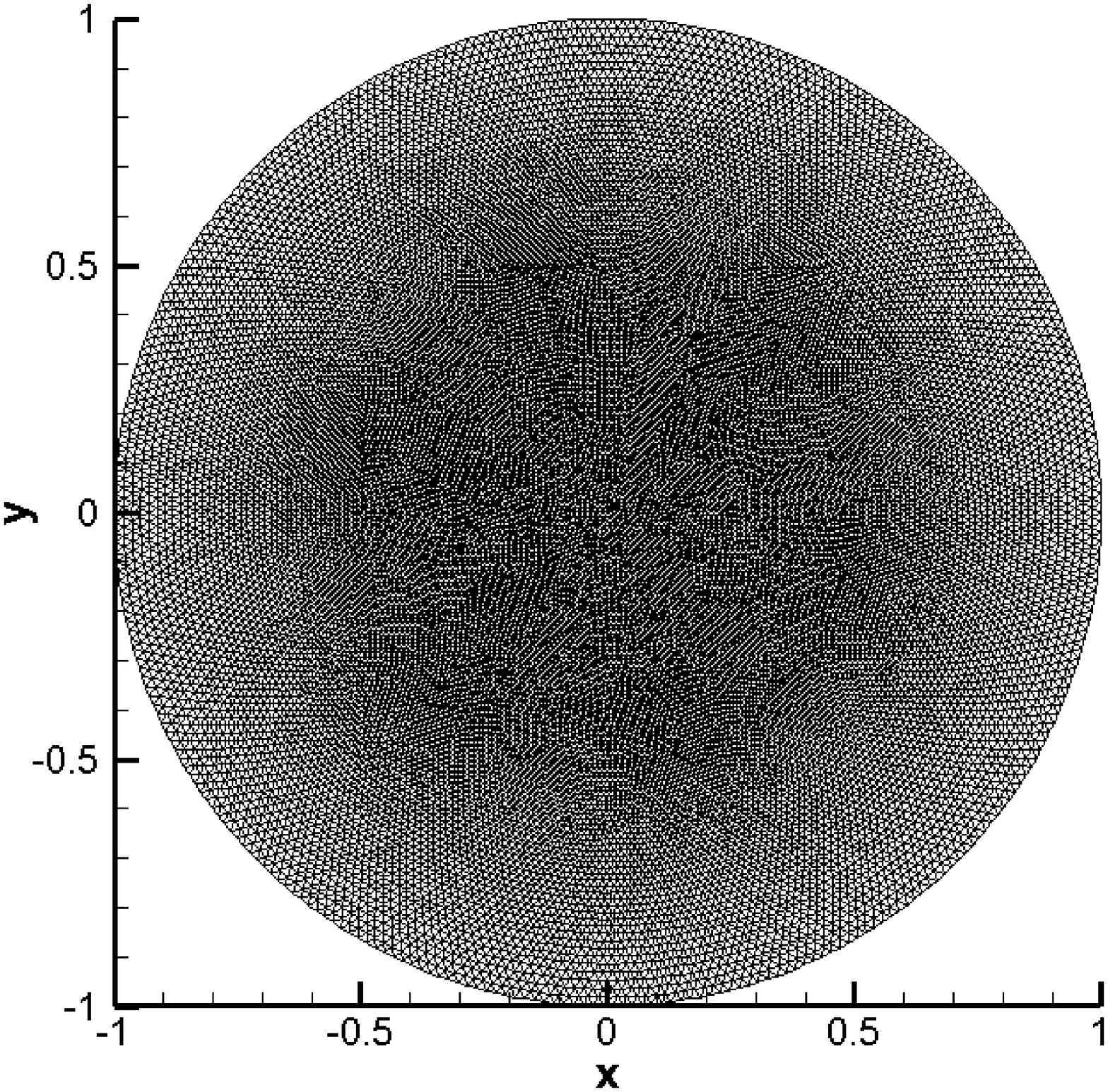}  &              
\includegraphics[width=0.47\textwidth]{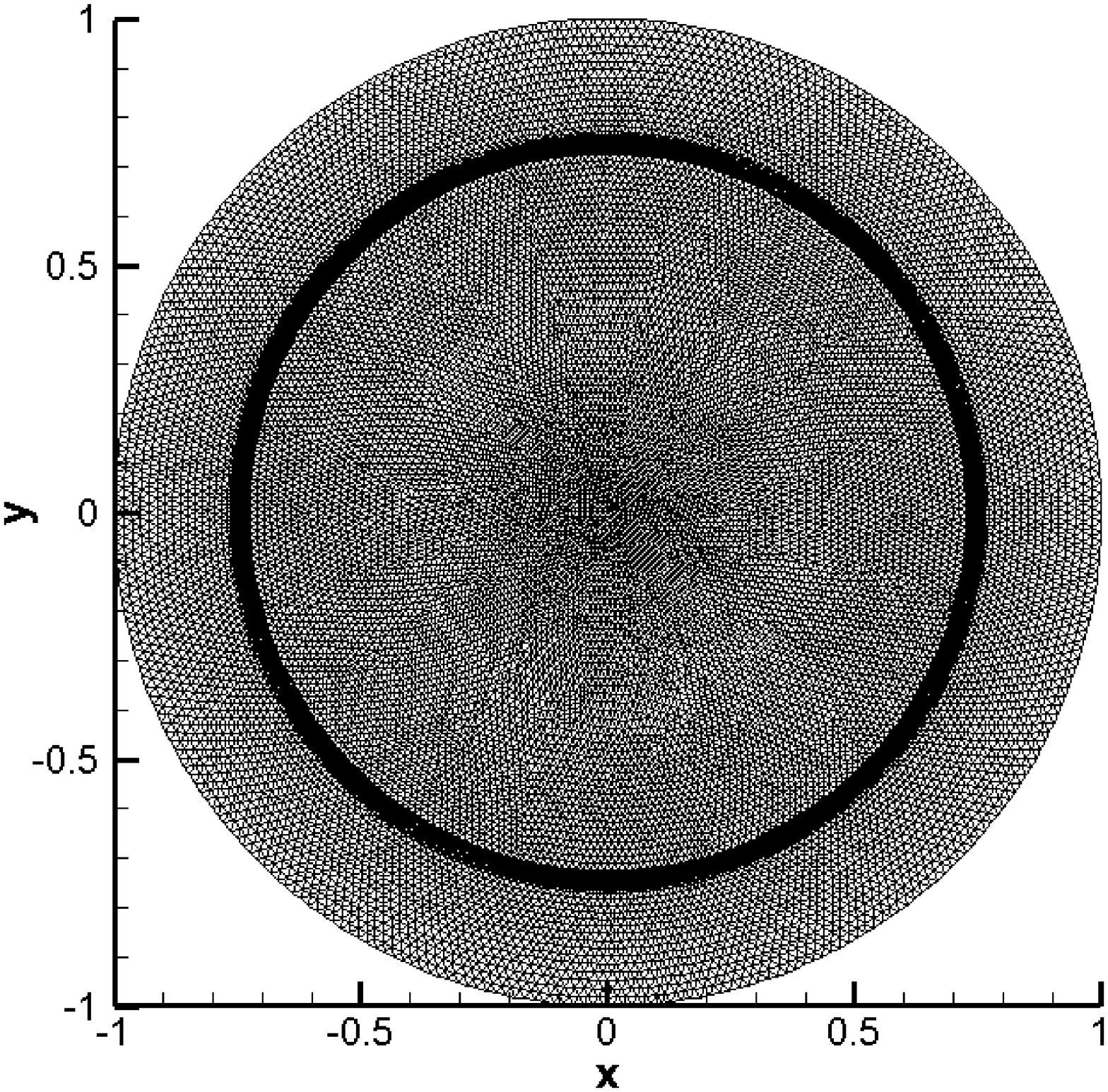}  \\            
\end{tabular} 
\caption{Initial (left) and final (right) mesh configuration for the explosion problem EP2. The strong shock generates a high compression of those elements which follow the wave.} 
\label{fig.EP3-grid}
\end{center}
\end{figure}

As proposed in \cite{Lagrange2D} a suitable reference solution can be obtained simplifying the two-dimensional Euler equations to a one-dimensional system with geometric source terms \cite{ToroBook}, which reads
\begin{equation}
\label{inhomEuler}
\Q_t + \F(\Q)_r = \S(\Q),
\end{equation}  
with
\begin{equation}
\label{matrixInhomEuler}
\Q = \left(\begin{array}{c} \rho \\ \rho u \\ \rho E \end{array}\right), \quad \F = \left(\begin{array}{c} \rho u \\ \rho u^2 + p \\  u(\rho E+p) \end{array}\right), \quad \S = -\frac{1}{r}\left(\begin{array}{c} \rho u \\ \rho u^2 \\  u(\rho E+p) \end{array}\right)\,.
\end{equation}
Here $r$ and $u$ represent the radial direction and the radial velocity, respectively.
As a result, a proper reference solution is obtained after solving the inhomogeneous system of equations \eqref{inhomEuler}-\eqref{matrixInhomEuler} on a one-dimensional mesh of 15000 points in the radial interval $r \in [0;1]$ using a classical second order TVD scheme \cite{ToroBook} with a Rusanov-type numerical flux. Third order accurate Lagrangian ADER-WENO schemes with time accurate local time stepping have been used together with the Osher-type numerical flux \eqref{eqn.osher} to compute the explosion problems EP1 and EP2. Figures \ref{fig.EP1}-\ref{fig.EP3} show a comparison between the  numerical solution obtained with the Lagrangian LTS scheme and the 1D reference solution. As for the Riemann problems presented in the previous section, one can appreciate the very good resolution of the contact wave in the density distribution and a good agreement with the reference solution is achieved also for horizontal velocity and pressure. We point out that EP2 is more challenging than EP1 because it involves a strong shock wave which causes a high compression of some elements in the mesh, as clearly depicted in Figure \ref{fig.EP3-grid}. By using the LTS approach we can avoid that those small triangles dictate the timestep for the entire mesh, hence allowing the other control volumes to reach the end of the simulation 
much faster and with a lower number of element updates, as highlighted in Table \ref{tab.EPdata}. 

\begin{figure}[!htbp]
\begin{center}
\begin{tabular}{cc} 
\includegraphics[width=0.47\textwidth]{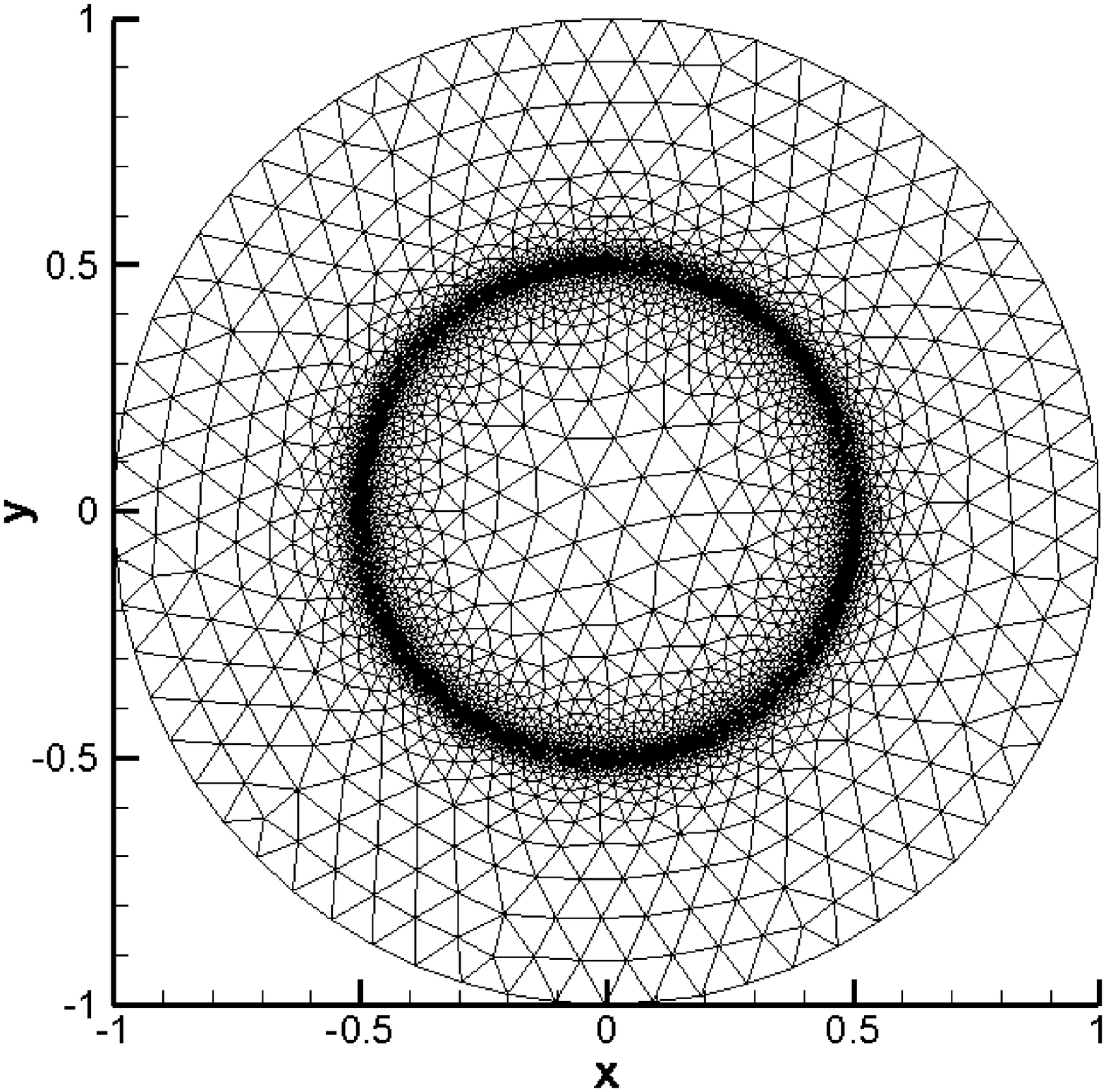}  &           
\includegraphics[width=0.47\textwidth]{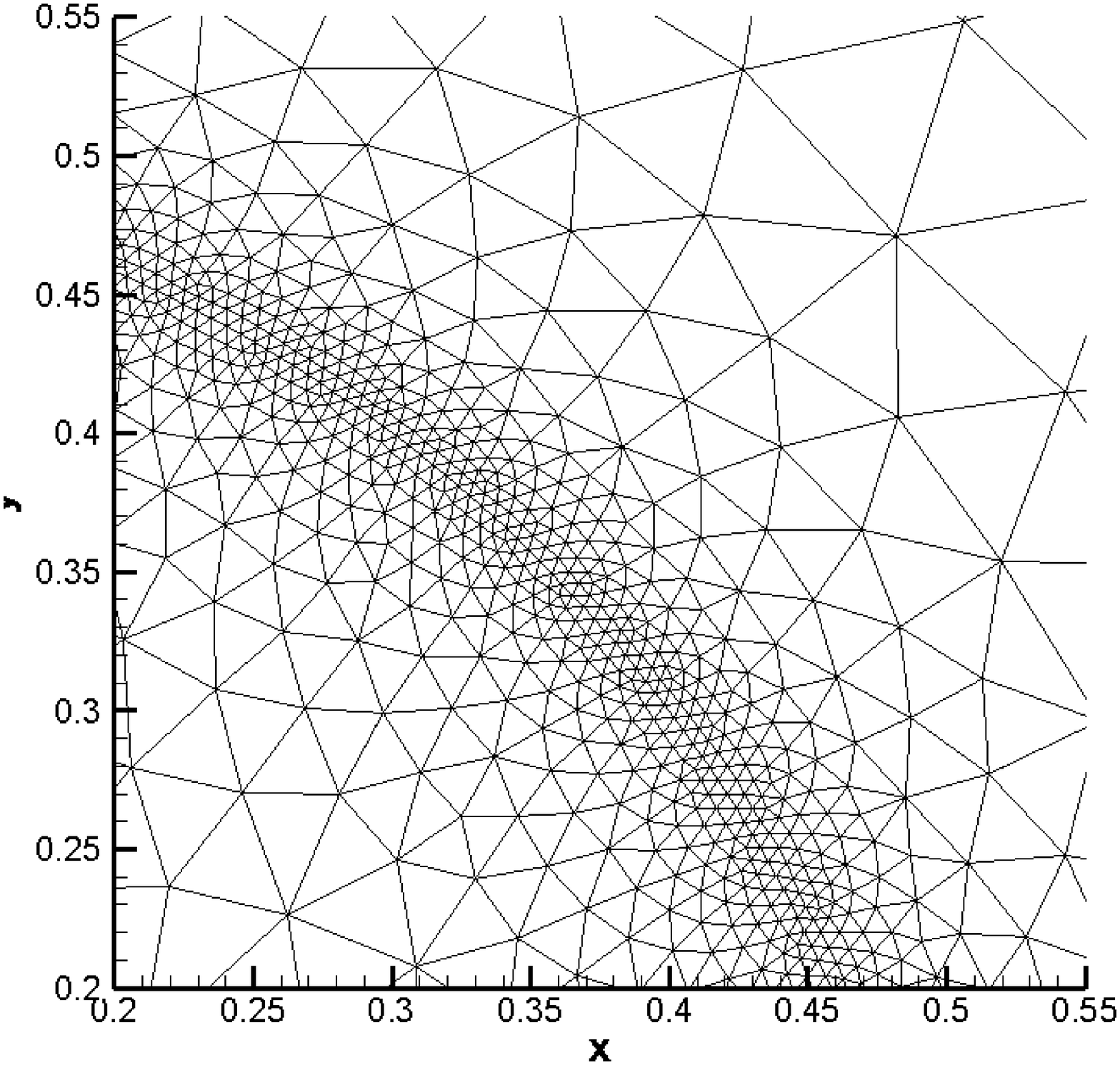} \\   
\includegraphics[width=0.47\textwidth]{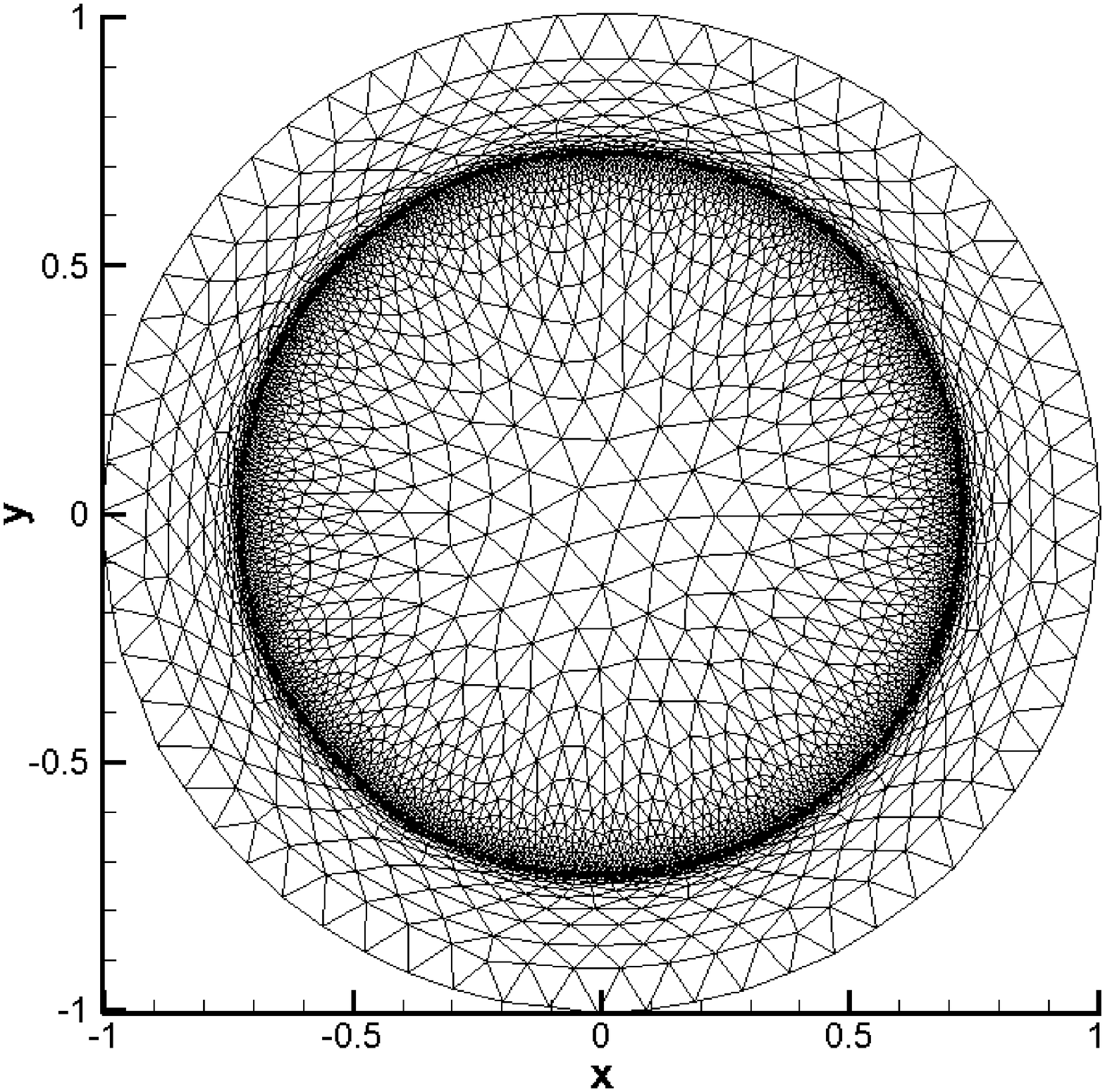}  &           
\includegraphics[width=0.47\textwidth]{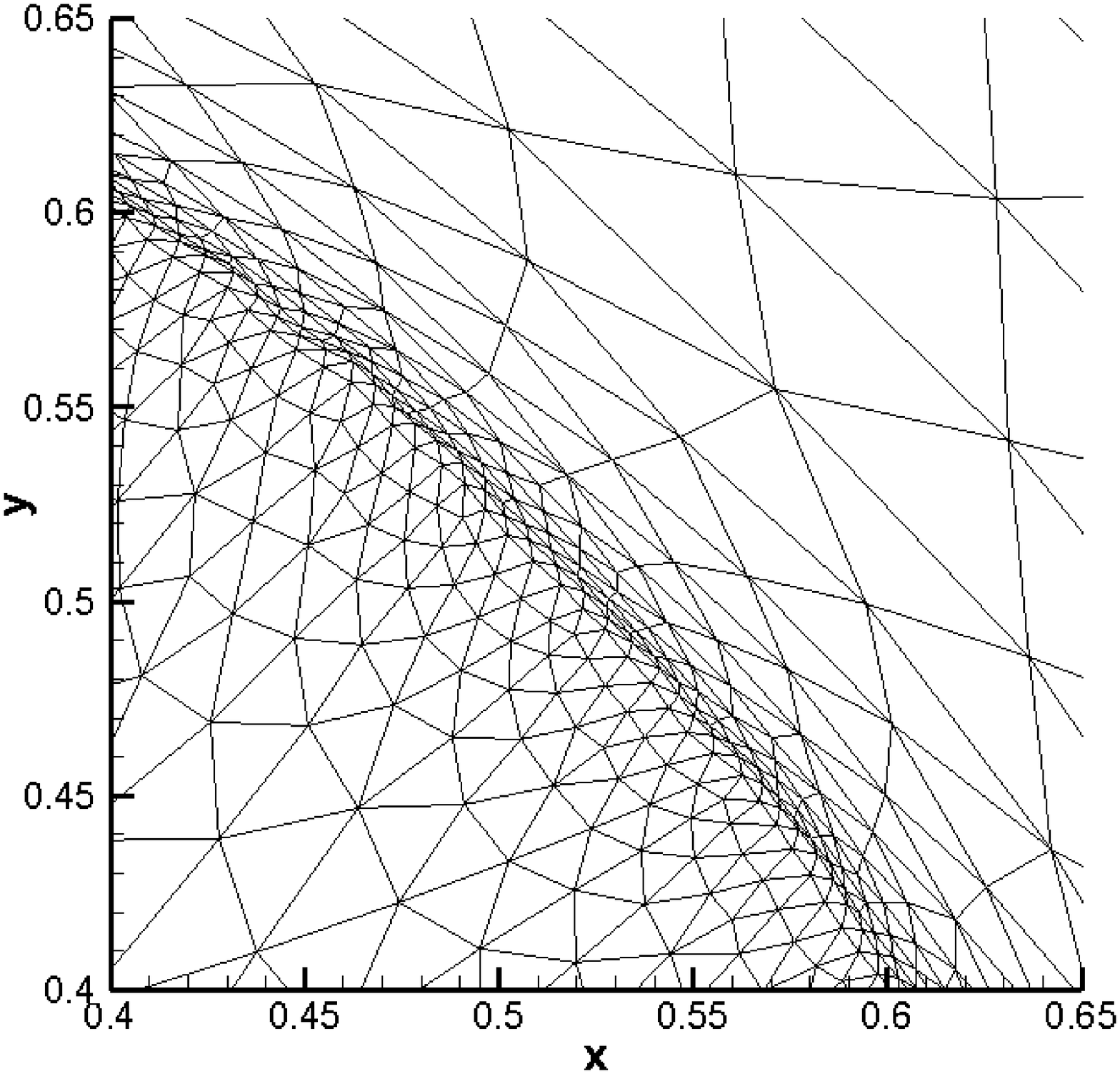} \\ 
\end{tabular} 
\caption{Initial (top) and final (bottom) configuration of the grid $m_2$ for the explosion problem EP1$^\ast$. A zoom of the mesh configuration across the contact wave is shown on the right at time $t=0.0$ (top) and $t=0.25$ (bottom).} 
\label{fig.EP1sizing-grid}
\end{center}
\end{figure}

The numerical simulation of many important phenomena arising in science and engineering typically requires the use of non-uniform computational grids with small elements 
clustered in some portions of the computational domain. In such circumstances, the use of a classical global time stepping algorithm would slow down the computation severely, 
since the smallest element of the mesh reduces the admissible timestep for the entire grid. 
Within the Eulerian framework on Cartesian grids, such a problem can be conveniently circumvented by resorting to Adaptive-Mesh-Refinement (AMR)  
with local time stepping, see e.g. \cite{Berger-Oliger1984,Berger-Colella1989,MuletAMR1,MuletAMR2,Burger2012,Dumbser2012b,Dumbser2014,Zanotti2013d}. 
There, the mesh is forced to refine only when and where this is needed, while it is recoarsened as soon as the chosen refinement criterion is no longer satisfied. 
An alternative option consists of preparing the computational mesh with a \textit{local static} refinement, which will remain fixed during the evolution if an 
Eulerian approach is adopted, while it will respond to the dynamics of the fluid if a Lagrangian framework is adopted, like in the present paper. In both cases, 
a local time stepping algorithm would make a huge difference in terms of computational efficiency, avoiding large control volumes to be slowed down by very small ones. 
Motivated by these considerations, we have run a modified version of the explosion problem EP1, denoted as EP1$^\ast$, which uses the same initial conditions of the former, 
apart for the mesh $m_2$, which has been built with a local mesh refinement around the initial location of the discontinuity, i.e. at $R=0.5$.  
More specifically, the mesh size is $h=1/100$ in the refined zone, and it grows with a growth rate of $s=1.5$ until $h=1/10$, which is used in the rest of the domain. 
Figure \ref{fig.EP1sizing-grid} shows the initial and the final configuration of the computational grid $m_2$ as well as a zoom onto the discontinuity. 
%The numerical results depicted in Figure \ref{fig.EP1sizing} have been collected with the fourth order version of the LTS Lagrangian ADER-WENO algorithm 
%with the Rusanov-type numerical flux \eqref{eqn.rusanov}. A good resolution at the contact discontinuity has been achieved, whereas the shock wave is not well resolved because the mesh $m_2$ is very coarse there. 
With this test problem we want to make another case for adopting the LTS approach rather than the classical GTS algorithm, and the advantages of the former can be easily 
deduced by looking at Table \ref{tab.EPdata}, where we compare the total number of element updates for each explosion problem needed to reach the final time of 
the simulation. 

%\begin{figure}[!htbp]
%\begin{center}
%\begin{tabular}{cc} 
%\includegraphics[width=0.47\textwidth]{./EP1sizing-LTS-rho3D.eps}  &           
%\includegraphics[width=0.47\textwidth]{./EP1sizing-LTS-rho.eps} \\    
%\end{tabular} 
%\caption{3D view of the density distribution (left) and comparison between reference and fourth order accurate numerical solution (right) for the explosion problem EP1$^\ast$ at the final time $t_f=0.25$ using the mesh $m_2$.} 
%\label{fig.EP1sizing}
%\end{center}
%\end{figure}

\begin{table}[!htbp]
  \caption{Comparison of the computational efficiency between GTS and LTS algorithm using the total number of element updates for EP1, EP2 and EP1$^*$.}
	\begin{center}
		\begin{tabular}{cccc}
		\hline
		         & \multicolumn{3}{c}{Number of element updates} 		\\
		Case     &      GTS                &     LTS                 &  GTS/LTS \\
		\hline
		EP1      & $  30.804224 \cdot 10^6 $ & $12.206887 \cdot 10^6 $ & 2.52 \\
		EP2      & $ 181.412376 \cdot 10^6 $ & $38.274477 \cdot 10^6 $ & 4.74 \\
		EP1$^*$  & $  22.171520 \cdot 10^6 $ & $ 8.313123 \cdot 10^6 $ & 2.67 \\
		\hline
		\end{tabular}
	\end{center}
	\label{tab.EPdata}
\end{table}

\subsection{The Kidder problem} 
\label{sec.Kidder}
The Kidder problem is a classical benchmark problem for Lagrangian algorithms. It has been widely used in the literature \cite{Maire2009,Despres2009} in order to assure that no spurious entropy is produced by the Lagrangian scheme. This test case was first designed by Kidder in \cite{Kidder1976} and it consists of an isentropic compression of a portion of a shell filled with an ideal gas. The shell $\Omega(0)$ is initially bounded by $r_i(t) \leq r \leq r_e(t)$, where $r=\sqrt{x^2+y^2}$ represents the general radial coordinate while $r_i(t),r_e(t)$ denote the time-dependent internal and external radius, respectively. The perfect gas is initially assigned with the following vector of primitive variables $\U_0$:
\begin{equation}
\U_0 = \left( \begin{array}{c} \rho_0(r) \\ u_0(r) \\ v_0(r) \\ p_0(r) \end{array}  \right) = \left( \begin{array}{c}  \left(\frac{r_{e,0}^2-r^2}{r_{e,0}^2-r_{i,0}^2}\rho_{i,0}^{\gamma-1}+\frac{r^2-r_{i,0}^2}{r_{e,0}^2-r_{e,0}^2}\rho_{e,0}^{\gamma-1}\right)^{\frac{1}{\gamma-1}} \\ 0 \\ 0 \\ s_0\rho_0(r)^\gamma \end{array}  \right), 
\label{eq:KidderIC}
\end{equation}
where $\rho_{i,0}=1$ and $\rho_{e,0}=2$ are the initial values of density at the internal and the external frontier, respectively. According to \cite{Maire2009} the ratio of specific heats is $\gamma=2$ and the initial entropy distribution $s_0$ is assumed to be uniform, i.e. $s_0= \frac{p_0}{\rho_0^\gamma} = 1$.

\begin{figure}[!htbp]
\begin{center}
\begin{tabular}{cc} 
\includegraphics[width=0.47\textwidth]{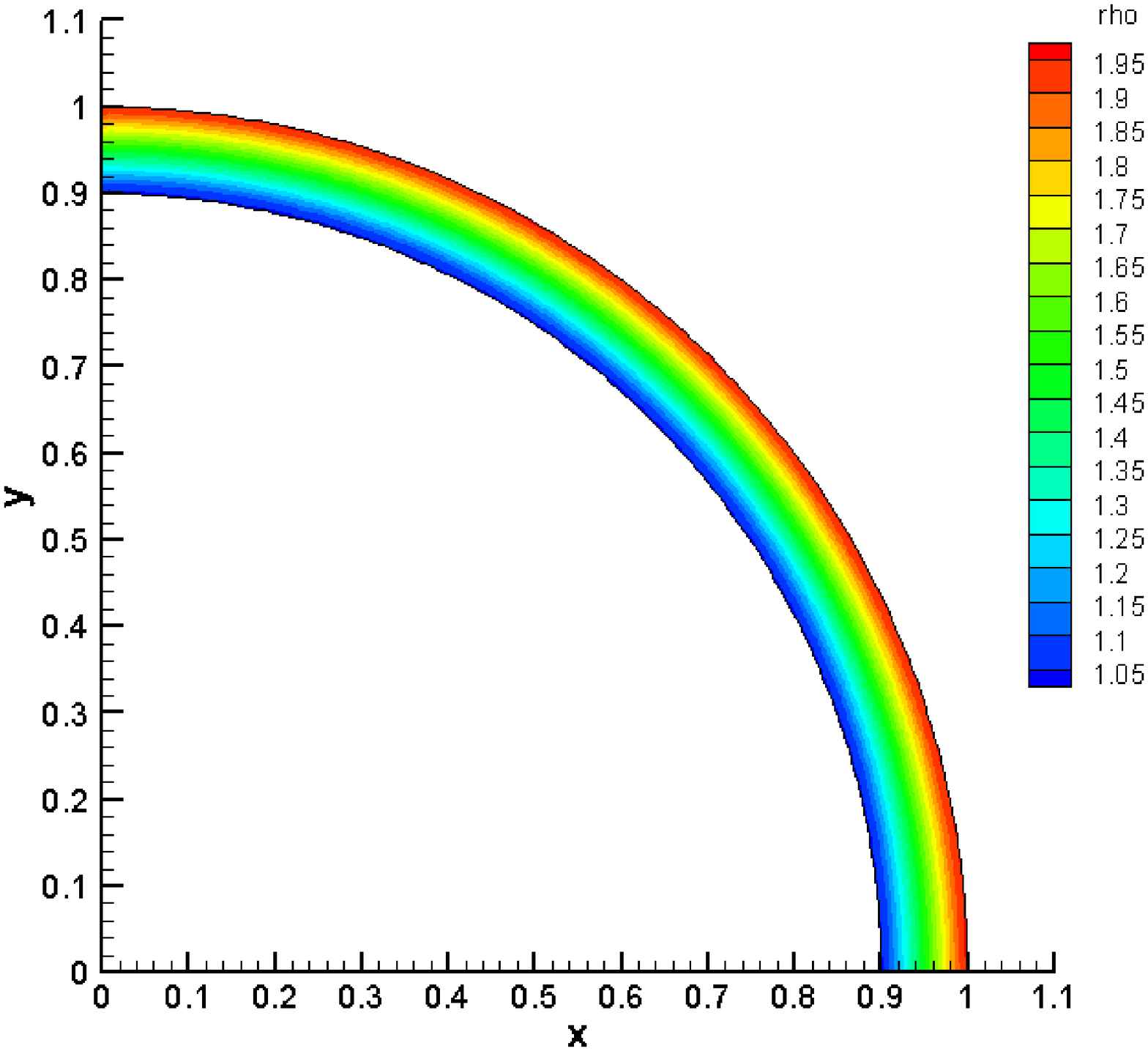}  &           
\includegraphics[width=0.47\textwidth]{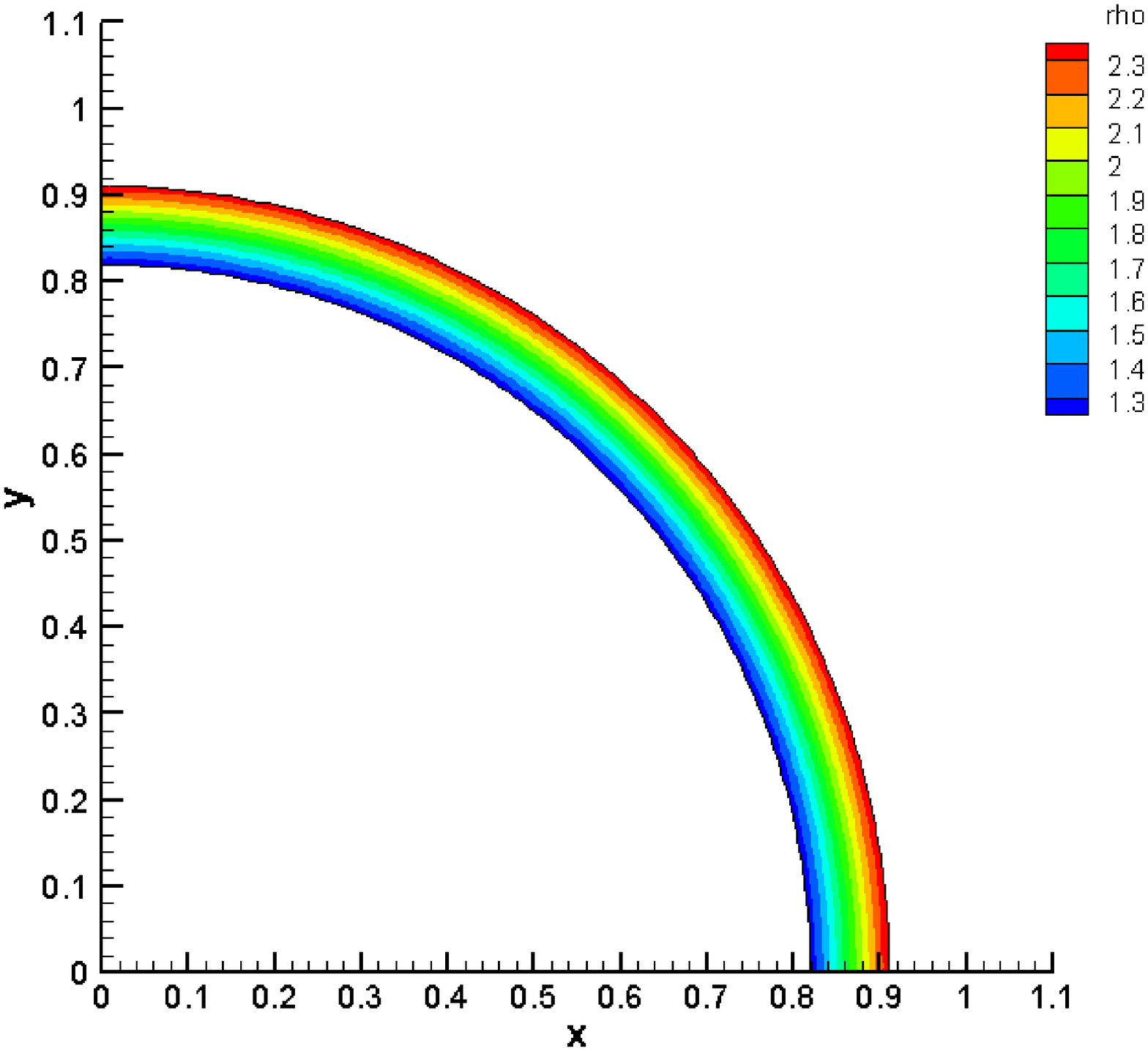} \\   
\includegraphics[width=0.47\textwidth]{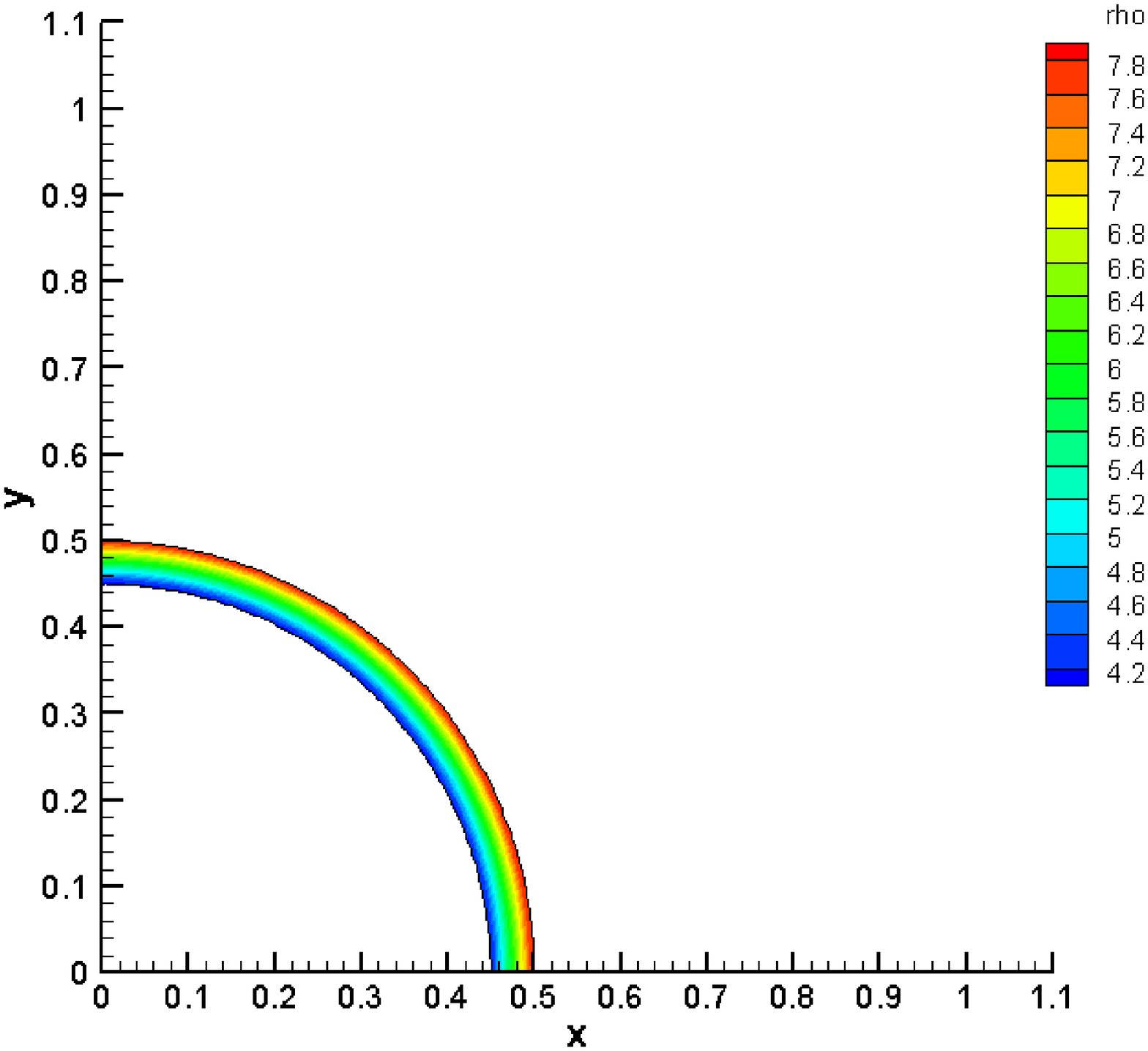}  &           
\includegraphics[width=0.47\textwidth]{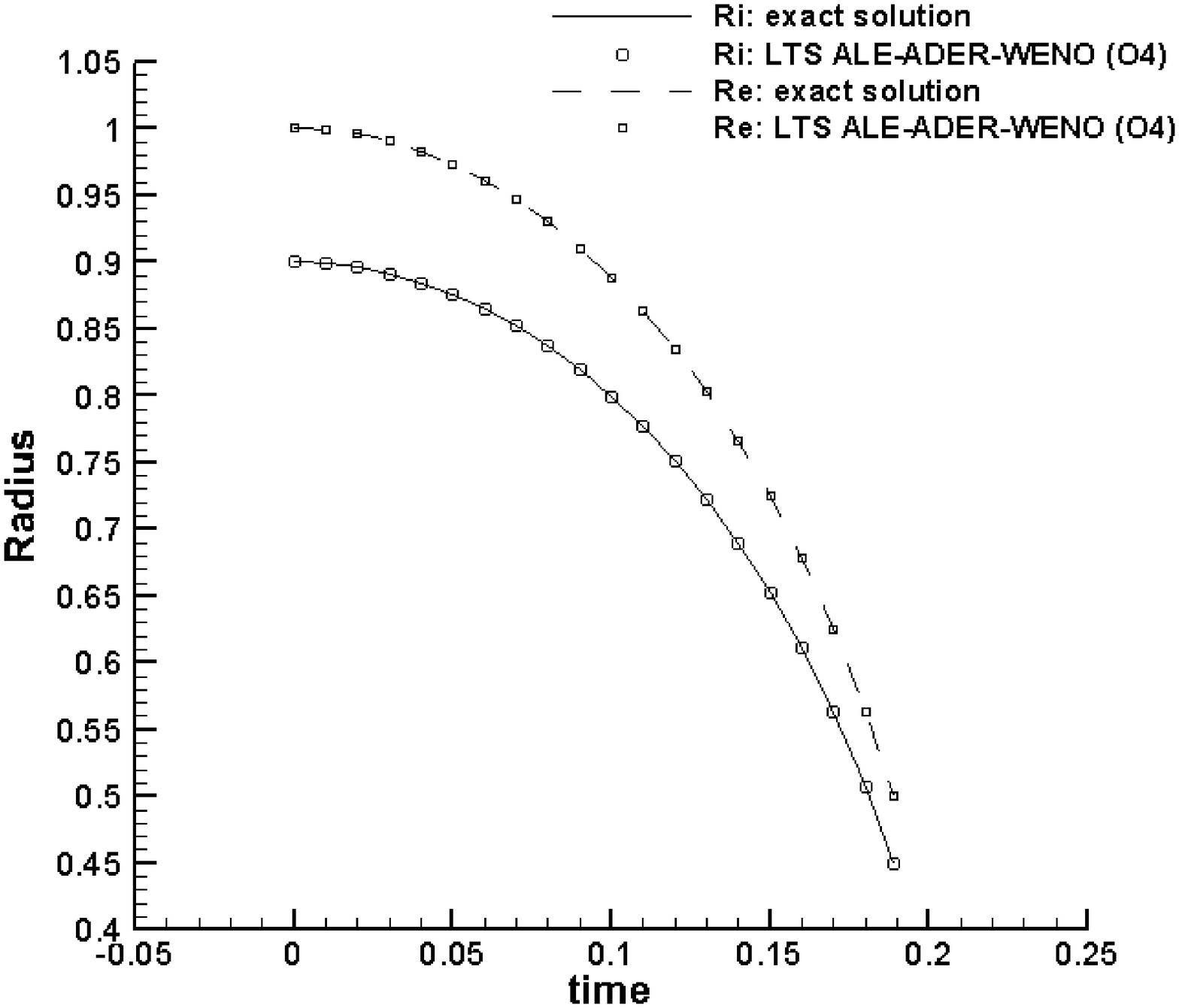} \\ 
\end{tabular} 
\caption{Fourth order accurate density distribution for the Kidder problem at the initial time $t=0.0$ (top left), at $t=0.9$ (top right) and at the final time $t=t_f$ (bottom left). The evolution of the internal and external radius of the shell is also shown (bottom right) and compared with the analytical solution.}
\label{fig.Kidder}
\end{center}
\end{figure}

Sliding wall boundary conditions are imposed on the horizontal and vertical edges that bound the portion of the shell, while on the internal and on the external frontier we set a space-time dependent state, which is assigned according to the exact solution $R(r,t)$ \cite{Kidder1976}. The analytical solution for the Kidder problem is given at the general time $t$ for a fluid particle initially located at radius $r$ as a function of the radius and of the homothety rate $h(t)$, i.e.
\begin{equation}
  R(r,t) = h(t)r, \qquad h(t) = \sqrt{1-\frac{t^2}{\tau^2}},
\label{eqKidderEx}
\end{equation}
where $\tau$ is the focalisation time 
\begin{equation}
\tau = \sqrt{\frac{\gamma-1}{2}\frac{(r_{e,0}^2-r_{i,0}^2)}{c_{e,0}^2-c_{i,0}^2}}
\end{equation}
with $c_{i,e}=\sqrt{\gamma\frac{p_{i,e}}{\rho_{i,e}}}$ representing the internal and external sound speeds. Following \cite{Despres2009,Maire2009}, the final time of the simulation is chosen to be $t_f=\frac{\sqrt{3}}{2}\tau$, so that the compression rate is $h(t_f)=0.5$ and the exact solution is given by the shell located within the interval $0.45 \leq R \leq 0.5$. We use a fourth order accurate  version of our new Lagrangian ADER-WENO scheme with LTS using the Osher-type numerical flux \eqref{eqn.osher}. The results are depicted in Figure \ref{fig.Kidder}, which shows the numerical solution for density at three different output times $t=0.0$, $t=0.9$ and $t=t_f$. Moreover, the evolution of the internal and external radius of the shell has been monitored during the simulation and Table \ref{tab.Kidder} reports the absolute error $|err|$ of the frontier positions, which is defined as the difference between the analytical and the numerical location of the internal and external radius at the final time. 

\begin{table}[!htbp]
  \caption{Absolute error for the internal and external radius location between exact $R_{ex}$ and numerical $R_{num}$ solution.}
	\begin{center}
		\begin{tabular}{cccc}
		\hline
		              		  			& $R_{ex}$ 		& $R_{num}$  & $|err|$     \\
		\hline
		\textit{Internal radius}	& 0.45000000 	& 0.44996063  & 3.94E-05 \\
		\textit{External radius}	& 0.50000000 	& 0.49930053  & 6.99E-04 \\
		\hline
		\end{tabular}
	\end{center}
	\label{tab.Kidder}
\end{table}

\subsection{The Saltzman problem} 
\label{sec.Saltzman}
\begin{figure}[!htbp]
\begin{center}
\begin{tabular}{cc} 
\includegraphics[width=0.47\textwidth]{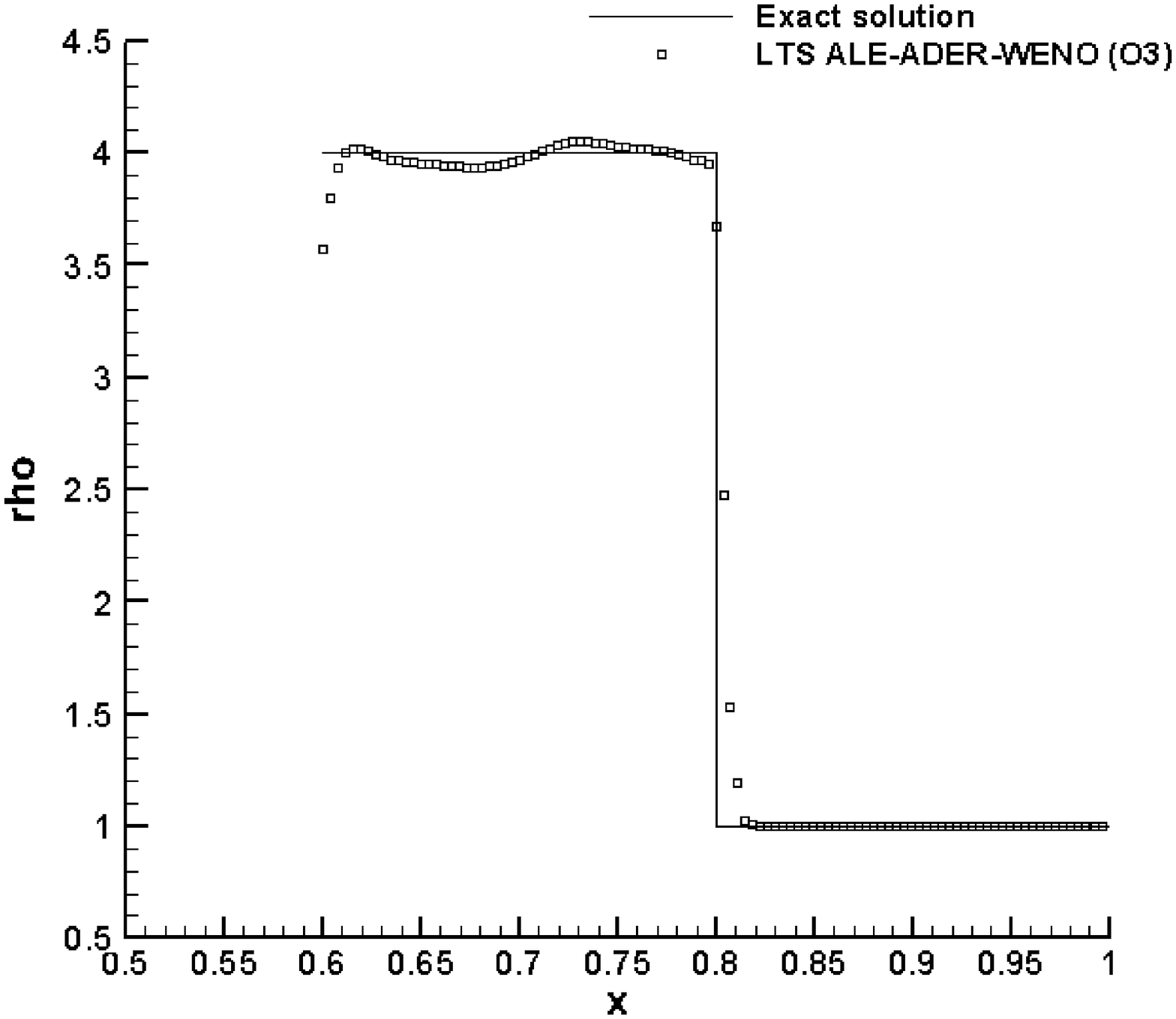}  &           
\includegraphics[width=0.47\textwidth]{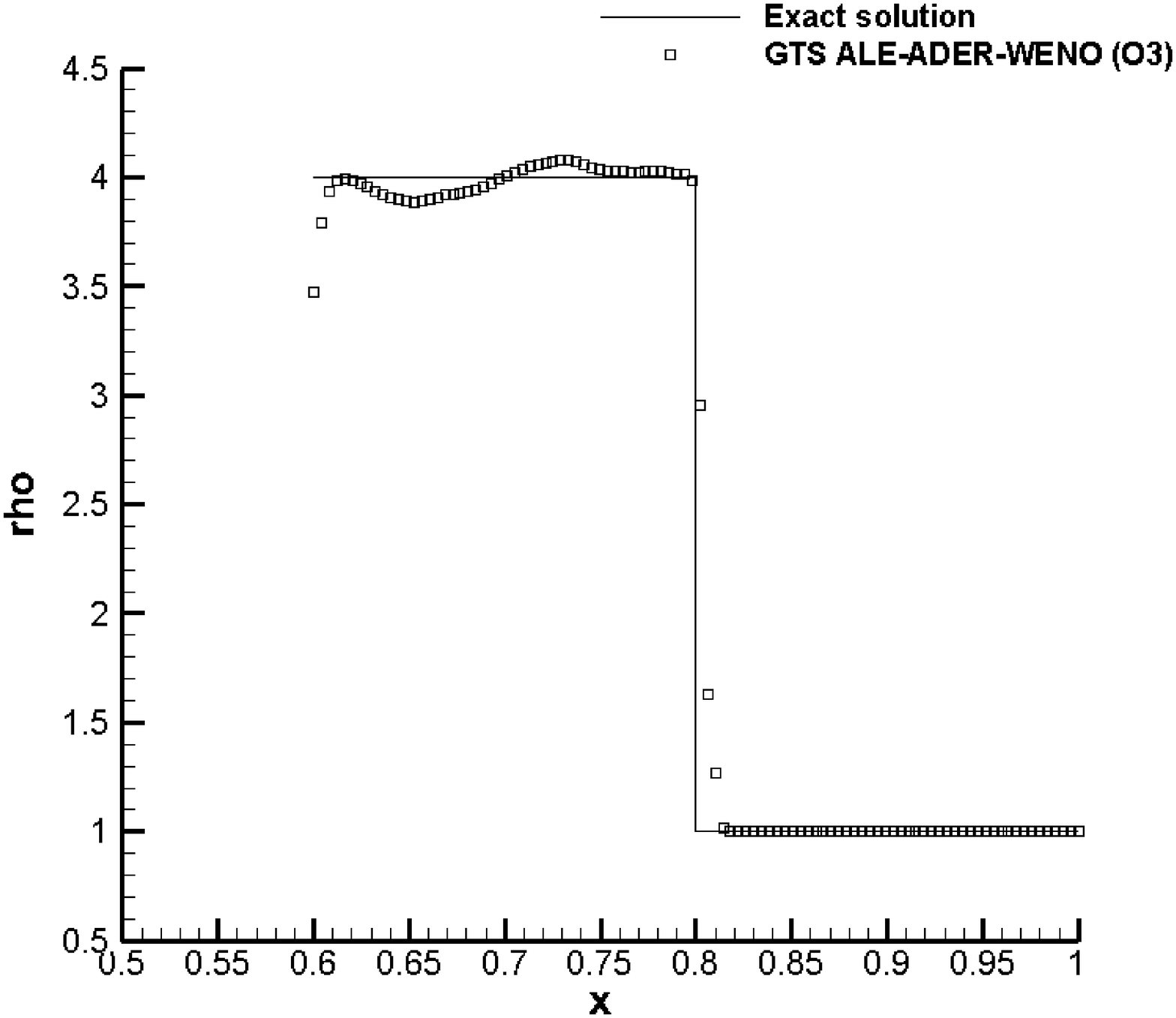} \\  
\includegraphics[width=0.47\textwidth]{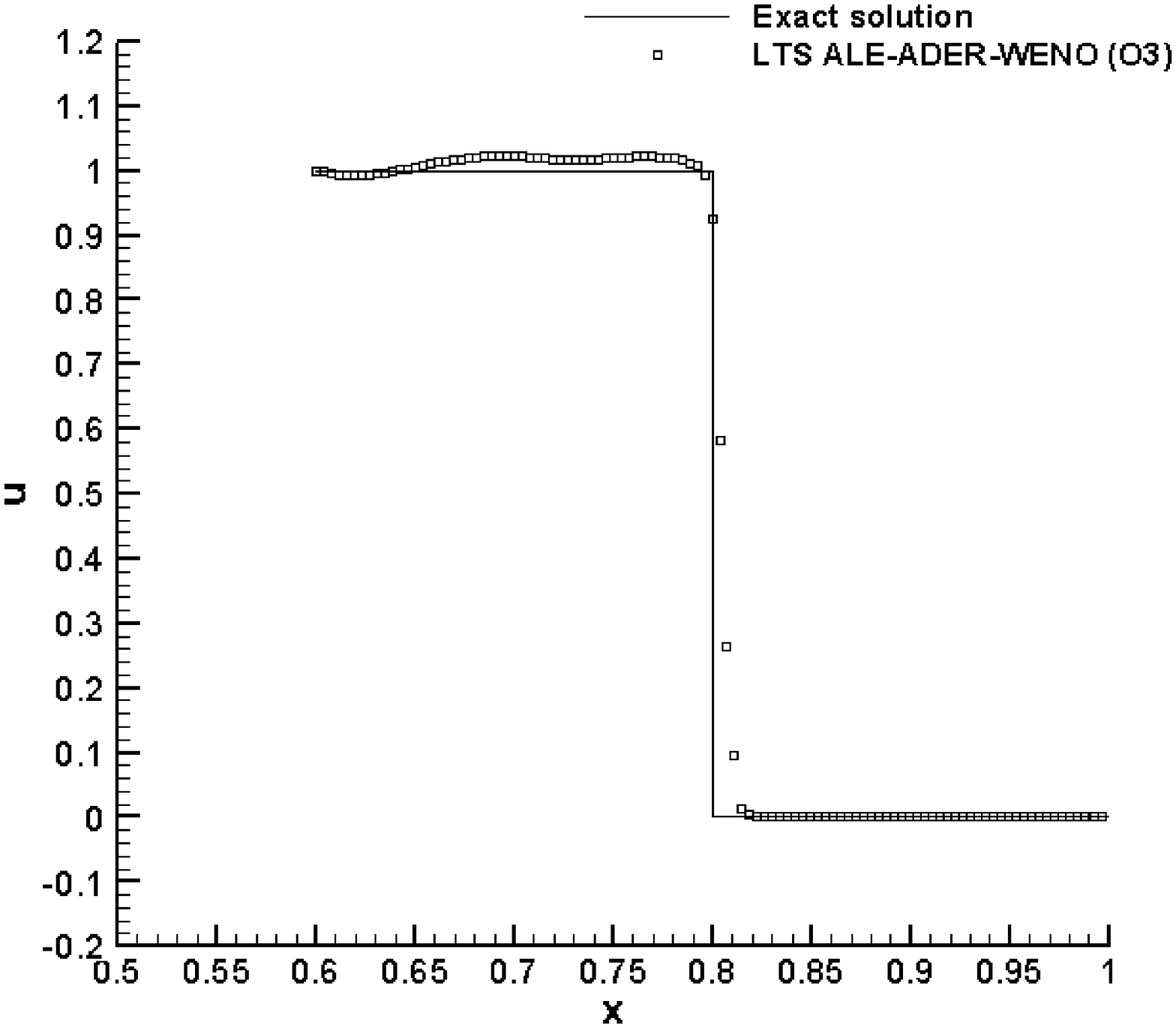}  &           
\includegraphics[width=0.47\textwidth]{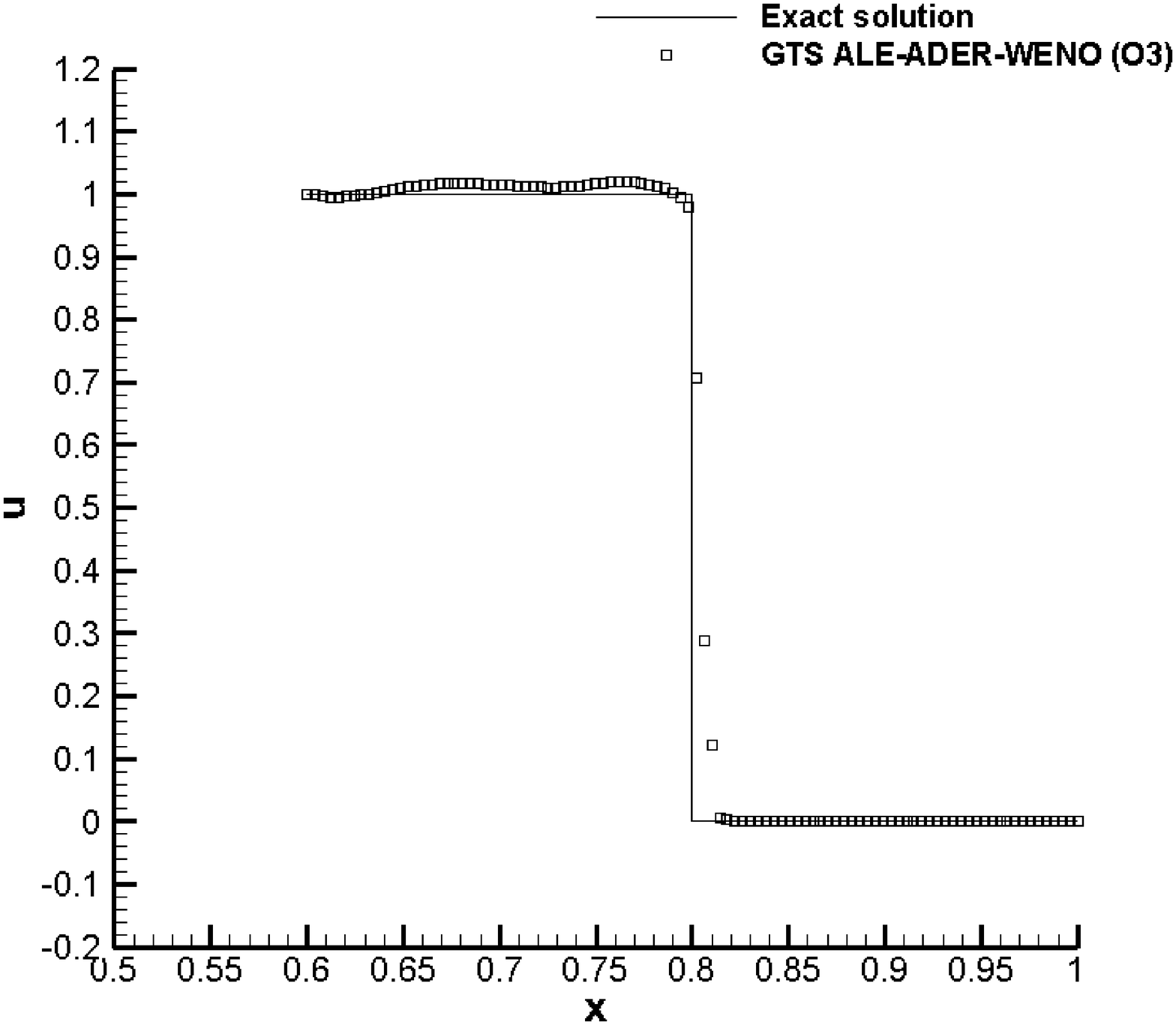} \\ 
\end{tabular} 
\caption{Third order accurate numerical solution for the Saltzman problem at the final time $t_f=0.6$. 
{\em Left panels}: solution obtained with LTS.
{\em Right panels}: solution obtained with GTS.
}
\label{fig.Saltz2D}
\end{center}
\end{figure}

\begin{figure}[!htbp]
\begin{center}
\begin{tabular}{cc}            
\includegraphics[width=0.94\textwidth]{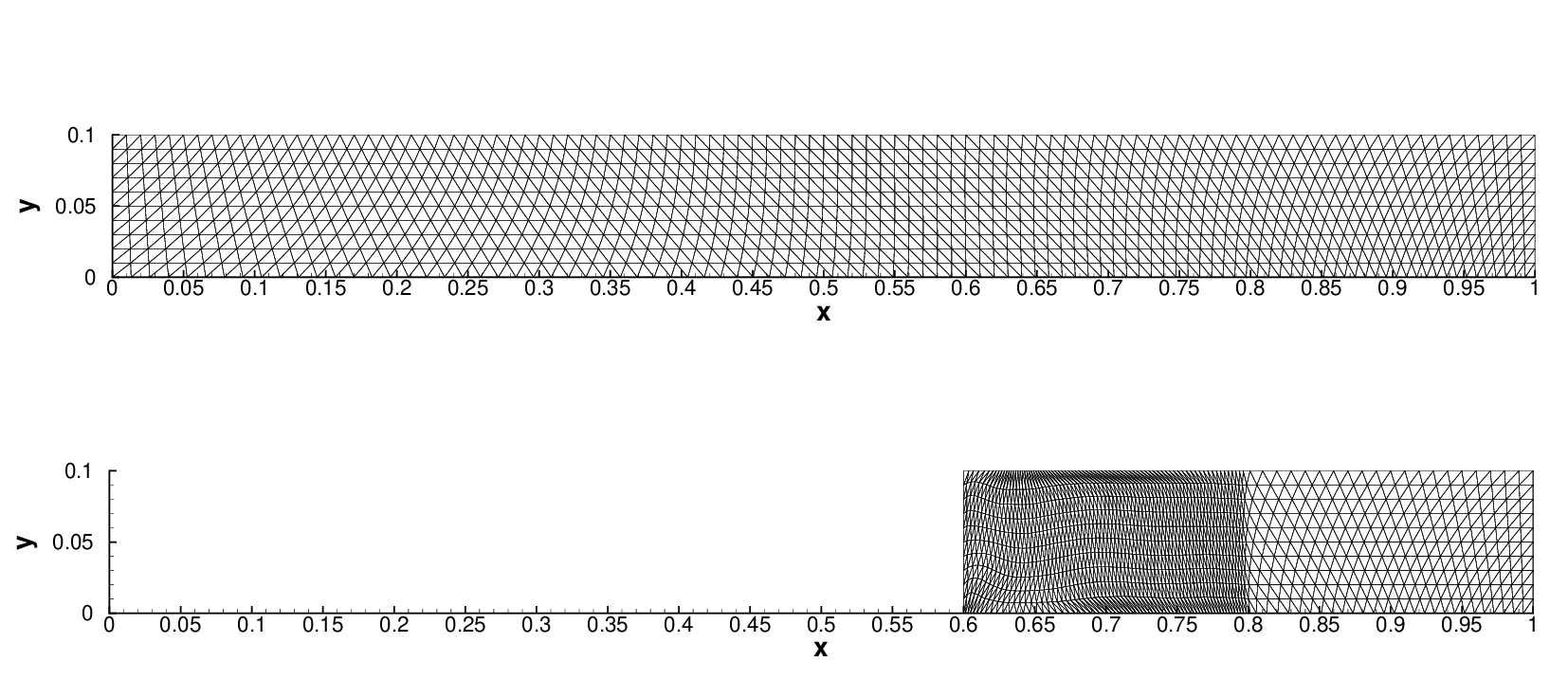}
\end{tabular} 
\caption{Initial and final mesh configuration for the Saltzman problem using the new LTS algorithm.}
\label{fig.Saltz2Dgrid}
\end{center}
\end{figure}

Another classical test case for Lagrangian gas dynamics is the Saltzman problem, which was presented for the first time by Dukowicz et al. in \cite{SaltzmanOrg} for a two-dimensional Cartesian grid that has been skewed in such a way that no element edges are aligned with the main fluid flow. It is a very challenging test problem 
against which any Lagrangian scheme ought to be validated \cite{Maire2009,chengshu2}.
It involves a strong shock wave caused by a piston that is moving along the main direction $x$ of the initial rectangular domain $\Omega(0)=[0;1]\times[0;0.1]$, which is initially discretized by $100 \times 10$ \textit{square} elements. According to \cite{Maire2009,chengshu2}, each element is then split into two \textit{right triangles}, so that we obtain a total number of elements of $N_E=2 \cdot 100 \times 10 = 2000$, and finally the following mapping is applied in order to skew the mesh: 
\begin{equation}
  x' = x + \left( 0.1 - y \right) \sin(\pi x) \qquad y' = y,
\end{equation}
where $\mathbf{x}=(x,y)$ represents the coordinate vector of the uniform grid, while $\mathbf{x'}=(x',y')$ are the final skewed coordinates. As done in \cite{chengshu2} the fluid is initially at rest and is assigned an internal energy $e_0=10^{-4}$ and a density $\rho_0=1$, hence the initial condition in terms of conserved variables reads $\Q_{0} = \left( \rho_0, \rho u_0, \rho v_0, \rho E_0 \right) = \left( 1, 0, 0, 10^{-4} \right)$. According to \cite{chengshu2}, the ratio of specific heats is set to $\gamma = \frac{5}{3}$ and the final time is assumed to be $t_f=0.6$, while the piston is moving with velocity $\mathbf{v}_p = (1,0)$ towards the right boundary of the domain. Moving slip wall boundary condition is imposed on the piston, whereas fixed slip wall boundaries have been set on the remaining sides of the domain. As fully explained in \cite{Lagrange2D,ToroBook}, the exact solution $\Q_{ex}(\x,t)$ is computed by solving a one-dimensional Riemann problem and at the final time $t_f$ it is given by
\begin{equation}
  \Q_{ex}(\x,t_f) = \left\{ \begin{array}{ccc} \left( 4, 1, 0, 2.5     \right) & \textnormal{ if } & x \leq x_f, \\
                                               \left( 1, 0, 0, 10^{-4} \right) & \textnormal{ if } & x > x_f,        
                      \end{array}  \right. 
\end{equation}
where $x_f=0.8$ denotes the final shock location. The piston is moving very fast, so that the fluid next to the piston is highly compressed and elements there must typically obey a severe CFL condition.
In practice, 
we have to start the simulation with CFL=0.1, hence using very small and \textit{global} timesteps. After time $t=0.01$ the numerical scheme proceeds with the new time accurate \textit{local time stepping} algorithm described in this article. We have used the third order version of our LTS Lagrangian ADER-WENO schemes and the very robust Rusanov-type numerical flux \eqref{eqn.rusanov}. Figure \ref{fig.Saltz2D} shows a comparison between the exact and the numerical solution for density and horizontal velocity at the final time of the simulation for both the LTS and the GTS version of our algorithm, while the initial and the final mesh configurations are depicted in Figure \ref{fig.Saltz2Dgrid}. An overall good agreement of the numerical solution with the exact solution can be observed and the decrease of the density which occurs near the  piston is due to the well known \textit{wall-heating problem}, see \cite{toro.anomalies.2002}. Furthermore we point out that the results obtained with the LTS scheme given in the left column of Figure \ref{fig.Saltz2D} do not differ very much from the numerical solution obtained with global time stepping (GTS) shown in the right column of Figure \ref{fig.Saltz2D}.  

\section{Conclusions}
\label{sec.concl}
In this article we have presented a high order Lagrangian finite volume schemes with time-accurate local time stepping (LTS) on moving unstructured triangular meshes. The numerical scheme is derived from \cite{Lagrange2D,LagrangeNC}, where a classical global time stepping approach is adopted, and from the recently developed one-dimensional high order Lagrangian LTS numerical scheme \cite{ALELTS1D}.  In our approach the WENO reconstruction technique is used to achieve high order of accuracy in space, while high order of accuracy in time is obtained via the local space-time Galerkin predictor. The new algorithm illustrated in this article is based on a non-conforming mesh in time, with hanging nodes that are continuously moving and in principle never match the same time level, unless either an 
intermediate output time or the final time of the simulation is reached. As a consequence, the reconstruction is carried out locally, i.e. within each control volume, using a virtual geometry and a  virtual set of cell averages of the surrounding elements that are both computed using the high order space-time predictor solution. In order to develop a fully conservative numerical scheme, the 
fluxes are evaluated relying on memory variables, which allow to record all fluxes accumulated in the past within each control volume. Unlike the one-dimensional version of the algorithm presented 
in \cite{ALELTS1D}, on two-dimensional unstructured meshes we need also to compute additional fluxes over triangular space-time sub-surfaces, whenever an element performs the update timestep. This
additional computational and algorithmic complexity is due to the increased complexity of the topology of a 2D mesh, which consists in control volumes, edges and nodes. 
By construction, our scheme is conservative and automatically satisfies the geometric conservation law (GCL) due to the integration over a closed space-time control volume. 

The algorithm has been applied to the Euler equations of compressible gas dynamics, solving a set of canonical test problems and benchmarks for Lagrangian schemes. Furthermore convergence rates 
up to  fourth order of accuracy in space and in time have been shown. 

Further work may contain the extension of the presented LTS algorithm to three space dimensions and non-conservative hyperbolic balance laws as well as the implementation of a proper treatment 
for stiff source terms, hence allowing the scheme to be applied to more complex systems of equations. 

\section*{Acknowledgments}
The presented research has been financed by the European Research Council (ERC) under the European Union's Seventh Framework 
Programme (FP7/2007-2013) with the research project \textit{STiMulUs}, ERC Grant agreement no. 278267. The authors acknowledge 
PRACE for awarding us access to the SuperMUC supercomputer of the Leibniz Rechenzentrum (LRZ) in Munich, Germany.

\bibliography{Lagrange2D_LTS}
\bibliographystyle{plain}

\end{document}